\documentclass[a4paper]{amsart}
\usepackage{amssymb} 
\usepackage[T1]{fontenc}
\usepackage[latin1]{inputenc}
\usepackage{amsfonts}
\usepackage{amsxtra}
\usepackage{ae}
\usepackage[all]{xy}
\usepackage{enumerate}

\include{diagram}

\newcommand*{\ket}{\rangle}
\newcommand*{\bra}{\langle}

\newcommand*{\Cov}{\mathfrak{Cov}}

\newcommand*{\D}{\mathcal{D}}
\newcommand*{\E}{\mathcal{E}}
\newcommand*{\cotimes}{\hat{\otimes}}
\newcommand*{\SHom}{\mathfrak{Hom}}
\newcommand*{\LSMod}{\textsf{-}\mathfrak{Mod}}
\newcommand*{\ev}{ev}

\DeclareMathOperator{\pro}{pro}

\DeclareMathOperator{\Aut}{Aut}
\DeclareMathOperator{\Hom}{Hom}
\DeclareMathOperator{\Mor}{Mor}
\DeclareMathOperator{\LMod}{-Mod}
\DeclareMathOperator{\End}{End}
\DeclareMathOperator{\id}{id}
\DeclareMathOperator{\Ad}{ad}

\DeclareMathOperator{\Tr}{tr}

\newenvironment{bnum}
{\begin{list}{}
    {\setlength{\labelwidth}{15pt}
     \setlength{\leftmargin}{\labelwidth}
    }
}
{\end{list}}

\numberwithin{equation}{section}
\newtheorem{theorem}{Theorem}[section]
\newtheorem{prop}[theorem]{Proposition}
\newtheorem{lemma}[theorem]{Lemma}
\newtheorem{cor}[theorem]{Corollary}
\newtheorem{definition}[theorem]{Definition}

\begin{document}

\title{Equivariant periodic cyclic homology}
\author{Christian Voigt}
\address{Mathematisches Institut\\
         Westfälische Wilhelms-Universität Münster\\
         Einsteinstr.\ 62\\
         48149 Münster\\
         Germany
}
\email{cvoigt@math.uni-muenster.de}


\thanks{This research was supported by the EU-Network \emph{Quantum
    Spaces and Noncommutative Geometry} (Contract HPRN-CT-2002-00280)
  and the \emph{Deutsche Forschungsgemeinschaft} (SFB 478).}

\begin{abstract}
We define and study equivariant periodic cyclic homology for locally compact 
groups. This can be viewed as a noncommutative generalization of equivariant 
de Rham cohomology. Although the construction resembles the Cuntz-Quillen 
approach to ordinary cyclic homology, a completely new feature in the equivariant setting 
is the fact that the basic ingredient in the theory is not a complex in the usual sense. 
As a consequence, in the equivariant context only the periodic cyclic theory can 
be defined in complete generality. Our definition recovers particular cases studied 
previously by various authors. 
We prove that bivariant equivariant periodic cyclic homology 
is homotopy invariant, stable and satisfies excision in both variables. Moreover 
we construct the exterior product which generalizes the obvious composition product. 
Finally we prove a Green-Julg theorem in cyclic homology for compact groups and 
the dual result for discrete groups. 
\end{abstract}

\maketitle

\section{Introduction}

In the general framework of noncommutative geometry cyclic homology plays the role of de Rham 
cohomology ~\cite{Connes2}. It was introduced by Connes ~\cite{Connes1} as the target of the noncommutative 
Chern character. Besides cyclic cohomology itself Connes also defined periodic cyclic cohomology. The latter is 
particularly important because it is the periodic theory that gives de Rham cohomology in the commutative case. \\
In this paper we develop a general framework in which cyclic homology can be extended to the equivariant 
context. Special cases of our theory have been defined and studied by various authors ~\cite{BG}, 
~\cite{Brylinski1}, ~\cite{Brylinski2}, ~\cite{Bues1}, ~\cite{Bues2}, ~\cite{KKL1}, ~\cite{KKL2}. However, all these 
approaches are limited to actions of compact Lie groups or even finite groups. Hence a substantial open problem was 
how to treat non-compact groups. Even for compact Lie groups an important open question was how to give a correct 
definition of equivariant cyclic cohomology (in contrast to homology) apart from the case of finite groups. \\
In this paper we define and study bivariant equivariant periodic cyclic homology $ HP^G_*(A,B) $ for locally compact 
groups $ G $. Throughout we work in the setting of bornological vector spaces and use the theory of 
smooth representations of locally compact groups on bornological vector spaces developped by Meyer ~\cite{Meyersmoothrep}. 
In this way we can treat many interesting examples of group actions on algebras in a unified fashion. 
In particular we obtain a theory which applies to discrete groups and totally disconnected groups as 
well as to Lie groups. \\
The construction of the theory follows the Cuntz-Quillen approach to cyclic homology based on the $ X $-complex 
~\cite{CQ1}, ~\cite{CQ2}, ~\cite{CQ3}, ~\cite{CQ4}. In fact a certain part 
of the Cuntz-Quillen machinery can be carried over to the equivariant situation without change. However, a 
new feature in the equivariant theory is the fact that the basic objects are not complexes in the sense of homological 
algebra. More precisely, we define an equivariant version $ X_G $ of the $ X $-complex but the differential 
$ \partial $ in $ X_G $ does not satisfy $ \partial^2 = 0 $ in general. To describe this behaviour we introduce 
the notion of a paracomplex. It turns out that in order to obtain ordinary complexes 
it is crucial to work in the bivariant setting from the very beginning. 
Although many tools from homological algebra are not available for paracomplexes, the resulting 
theory is computable to some extent. We point out that the occurence of 
paracomplexes is the reason why we only define and study the periodic theory 
$ HP^G_* $. It seems to be unclear how ordinary equivariant cyclic homology $ HC^G_* $ can be defined 
correctly in general apart from the case of compact groups. \\
An important ingredient in the definition of $ HP^G_* $ is the algebra 
$ \mathcal{K}_G $ which can be viewed as a certain subalgebra of the algebra of compact operators on the regular 
representation $ L^2(G) $. 
For instance, if $ G $ is discrete the elements of $ \mathcal{K}_G $ are simply finite matrices indexed by $ G $. 
The ordinary Hochschild homology and cyclic homology of this algebra are rather trivial. 
However, in the equivariant setting $ \mathcal{K}_G $ carries homological information of the group $ G $ 
if it is viewed as a $ G $-algebra equipped with the action induced from the regular representation. 
This resembles the properties of the total space $ EG $ of the universal principal bundle over the 
classifying space $ BG $. As a topological space $ EG $ is contractible, but its equivariant cohomology 
is the group cohomology of $ G $. Moreover, in the classical theory an arbitrary 
action of $ G $ on a space $ X $ can be turned into a free action by 
replacing $ X $ with the $ G $-space $ EG \times X $. In our theory 
tensoring with the algebra $ \mathcal{K}_G $ is used to associate to an arbitrary 
$ G $-algebra another $ G $-algebra which is projective as a $ G $-module. \\
Let us now explain how the text is organized. In section \ref{secborn} we review basic definitions and 
results from the theory of bornological vector spaces and the theory of smooth representations of locally compact groups.  
After this we introduce the category of covariant modules in section \ref{seccov} and discuss the 
natural symmetry operator on this category. Covariant modules constitute 
the appropriate framework for studying equivariant cyclic homology. 
In section \ref{secprocat} we review some facts about pro-categories. Since the work of Cuntz and Quillen ~\cite{CQ4} it is 
known that periodic cyclic homology is most naturally defined for pro-algebras. The same holds true in the 
equivariant situation where one has to consider pro-$ G $-algebras. 
We introduce the pro-categories needed in our framework and fix some notation. 
In section \ref{secpara} we define paracomplexes 
and paramixed complexes. As explained above, paracomplexes play an important role in our theory. \\
After these preparations we define and study quasifree pro-$ G $-algebras in section \ref{secproalg}. This discussion extends 
in a straightforward way the theory of quasifree algebras introduced by Cuntz and Quillen. Next we define 
equivariant differential forms for pro-$ G $-algebras in section \ref{secdiffform} and show that one obtains 
paramixed complexes in this way. Equivariant differential forms are used to construct the equivariant 
$ X $-complex $ X_G(A) $ for a pro-$ G $-algebra $ A $ in section \ref{secX}. As mentioned before this leads to a paracomplex. 
We show that the paracomplexes obtained from the equivariant $ X $-complex and 
from the Hodge tower associated to equivariant differential forms are homotopy equivalent. In this way we 
generalize one of the main results of Cuntz and Quillen to the equivariant setting. The proof from the 
nonequivariant situation has to be modified because there is no spectral decomposition of the Karoubi 
operator available in the equivariant context.  
In section \ref{secHPdef} we give the definition of bivariant equivariant periodic cyclic homology 
$ HP^G_*(A,B) $ for pro-$ G $-algebras $ A $ and $ B $. We show that $ HP^G_* $ is homotopy 
invariant with respect to smooth equivariant homotopies and stable in a natural sense in both variables in the subsequent sections.
Moreover we prove that $ HP^G_* $ satisfies excision in both variables. 
This shows on a formal level that $ HP^G_* $ shares important properties 
with equivariant $ KK $-theory ~\cite{Kasparov1}, ~\cite{Kasparov2}. 
In section \ref{secext} we construct the exterior product for equivariant periodic cyclic homology. Again, the properties of this 
product are parallel to the situation in $ KK $-theory. \\
After these general considerations we explain in section \ref{seclie} how our definition is related to previous constructions 
in the literature. In particular we discuss the example of a compact Lie group $ G $ acting smoothly 
on a compact manifold $ M $. In this case the equivariant cyclic homology of $ C^\infty(M) $ has been
computed by Block and Getzler ~\cite{BG}. This example is illuminating since 
it exhibits the relations between equivariant cyclic homology and the classical Cartan model of equivariant cohomology ~\cite{Cartan1}, ~\cite{Cartan2}. 
In fact, one may think of equivariant cyclic homology as a noncommutative version of the Cartan model. \\
Finally, we prove a homological version of the Green-Julg theorem 
$ HP^G_*(\mathbb{C},A) \cong HP_*(A \rtimes G) $ for compact groups in section \ref{secGJ} and 
the dual result $ HP^G_*(A, \mathbb{C}) \cong HP^*(A \rtimes G) $ for discrete 
groups in section \ref{secDGJ}. Again this is analogous to the situation in $ KK $-theory. \\
We do not treat the construction of a Chern character from equivariant $ K $-theory into  
equivariant cyclic homology in this paper. Let us remark that for compact Lie groups and finite groups partial 
Chern characters have been defined before ~\cite{BG}, ~\cite{KKL2}. \\
This paper is based on the main part of my thesis ~\cite{Voigtthesis} which was written under the 
direction of Prof. Dr. J. Cuntz. 

\section{Bornological vector spaces and smooth representations}
\label{secborn}

In this section we recall some basic results of the 
theory of bornological vector spaces and smooth representations of locally compact groups. For more information we refer to 
~\cite{H-L1}, ~\cite{H-L2}, ~\cite{Meyer}, ~\cite{Meyersmoothrep}, \cite{Meyerborntop}. \\
A convex bornology on a complex vector space $ V $ is a collection of subsets $ \mathfrak{S}(V) $ of $ V $ satisfying 
some axioms. The elements $ S \in  \mathfrak{S}(V) $ are called the small subsets of the bornology. 
The motivating example of a bornology is given by the collection of bounded subsets in a locally convex vector space. 
A bornological vector space is a vector space $ V $ together with a convex bornology $  \mathfrak{S}(V) $ on $ V $. 
A linear map $ f: V \rightarrow W $ between bornological vector spaces is called bounded if it maps small sets 
to small sets. The space of bounded linear maps from $ V $ to $ W $ is 
denoted by $ \Hom(V,W) $. Recall that a subset $ S $ of a complex vector space is called a disk if it is circled and convex. 
The disked hull $ S^\Diamond $ is the circled convex hull of $ S $. If $ S $ is a small subset in 
a bornological vector space then $ S^\Diamond $ is again small. To a disk $ S \subset V $ one associates 
the semi-normed space $ \bra S \ket $ which is defined as the linear span of $ S $ endowed with the 
semi-norm $ \| \cdot \|_S $ given by the Minkowski functional. The disk $ S $ is called norming if $ \bra S \ket $ is a 
normed space and completant if $ \bra S \ket $ is a Banach space. 
A bornological vector space is called separated if all disks 
$ S \in \mathfrak{S} $ are norming. It is called complete if each 
$ S \in \mathfrak{S} $ is contained in a completant small disk $ T \in \mathfrak{S} $. A complete bornological 
vector space is always separated. \\ 
We will usually only work with complete bornological vector spaces. To any 
bornological vector space $ V $ one can associate a complete bornological vector 
space $ V^c $ and a bounded linear map $ \natural: V \rightarrow V^c $ such that composition with $ \natural $ induces a 
bijective correspondence between bounded 
linear maps $ V^c \rightarrow W $ with complete 
target $ W $ and bounded linear maps $ V \rightarrow W $. 
In the category of complete bornological vector spaces direct sums, direct 
products, projective limits and inductive limits exist. In all these cases one has characterizations 
by universal properties. Moreover there exists a natural tensor product which is universal for bounded 
bilinear maps. \\
A complete bornological algebra is a complete bornological 
vector space $ A $ with an associative multiplication given as a bounded linear map 
$ m: A \hat{\otimes} A \rightarrow A $. A homomorphism between 
complete bornological algebras is a bounded linear map $ f: A \rightarrow B $ 
which is compatible with multiplication. Remark that complete bornological algebras are not assumed to have a unit. Even if 
$ A $ and $ B $ are unital a homomorphisms $ f: A \rightarrow B $ need not preserve the unit of $ A $. A homomorphism 
$ f: A \rightarrow B $ between unital bornological algebras  satisfying $ f(1) = 1 $ will be called a unital homomorphism. \\
We denote the unitarization of a complete bornological algebra $ A $ by $ A^+ $. 
It is the complete bornological algebra with underlying vector space 
$ A \oplus \mathbb{C} $ and multiplication defined by 
$ (a, \alpha)\cdot (b,\beta) = (ab + \alpha b + \beta a, \alpha \beta) $. If 
$ f: A \rightarrow B $ 
is a homomorphism between complete bornological algebras there exists a 
unique extension to a unital homomorphism $ f^+: A^+ \rightarrow B^+ $. \\
Let us discuss briefly the definition of a module over a complete bornological algebra $ A $. A left $ A $-module is a complete 
bornological vector space $ M $ together with a bounded linear map 
$ \lambda: A \hat{\otimes} M \rightarrow M $ satisfying the 
axiom $ \lambda(\id \hat{\otimes} \lambda) = \lambda(m \hat{\otimes} \id) $ for 
an action. A homomorphisms $ f: M \rightarrow N $ of $ A $-modules 
is a bounded linear map commuting with the action of $ A $. 
We denote by $ \Hom_A(M,N) $ the space of all $ A $-module homomorphisms. 
Let $ V $ be any complete bornological vector 
space. An $ A $-module of the form $ M = A^+ \hat{\otimes} V $ with 
action given by left multiplication is called the free $ A $-module over $ V $. 
If an $ A $-module $ P $ is a direct summand in a free $ A $-module it is called projective. 
Projective modules are characterized by the following property. If $ P $ is projective and $ f: M \rightarrow N $ a surjective 
$ A $-module homomorphism with a bounded linear splitting $ s: N \rightarrow M $ then any $ A $-module homomorphism 
$ g: P \rightarrow N $ can be lifted to an $ A $-module homomorphism 
$ h: P \rightarrow M $ such that $ fh = g $. \\
In a similar way one can define and study right $ A $-modules and $ A $-bimodules. 
We can also work in the unital category starting with 
a unital complete bornological algebra $ A $. A unitary module $ M $ over a unital complete bornological algebra $ A $ is an $ A $-module such that 
$ \lambda(1 \otimes m) = m $ for 
all $ m \in M $. In the category of unitary modules the modules of 
the form $ A \hat{\otimes} V $ where $ V $ is a complete bornological vector 
space are free. Projective modules are again direct summands of free modules and can 
be characterized by a lifting property as before. \\
Let us briefly discuss the most relevant examples of bornological vector spaces. 

\subsection*{Fine spaces}

Let $ V $ be an arbitrary complex vector space. The fine bornology 
$ \mathfrak{Fine}(V) $ 
is the smallest possible bornology on $ V $. This means that $ S \subset V $ is 
contained in $ \mathfrak{Fine}(V) $ iff it is a bounded subset of a finite dimensional 
subspace of $ V $. 
It follows immediately from the definitions that 
all linear maps $ f: V \rightarrow W $ from a fine space $ V $ into any 
bornological space $ W $ are bounded. In particular we obtain 
a fully faithful functor $ \mathfrak{Fine} $ from the category of complex vector 
spaces into the category of complete bornological vector spaces. 
This embedding is compatible with tensor products. If $ V_1 $ and 
$ V_2 $ are fine spaces the completed bornological tensor product 
$ V_1 \hat{\otimes} V_2 $ is the algebraic tensor product 
$ V_1 \otimes V_2 $ equipped with the fine bornology. 
In particular every algebra $ A $ over the complex numbers can be viewed as a complete bornological 
algebra with the fine bornology. \\
Since the completed bornological tensor product is compatible with 
direct sums we see that $ V_1 \hat{\otimes} V_2 $ is as a vector space 
simply the algebraic tensor product $ V_1 \otimes V_2 $ provided $ V_1 $ or 
$ V_2 $ is a fine space. However, the bornology on the tensor product is in 
general not the fine bornology. 

\subsection*{Locally convex spaces}

The most important examples of bornological vector spaces are obtained 
from locally convex vector spaces. If $ V $ is any locally convex vector space one 
can associate two natural bornologies $ \mathfrak{Bound}(V) $ and 
$ \mathfrak{Comp}(V) $ to $ V $ which are called the bounded bornology and the 
precompact bornology, respectively. \\
The elements in $ \mathfrak{Bound}(V) $ are by definition the 
bounded subsets of $ V $. Equipped with the bornology $ \mathfrak{Bound}(V) $ the space $ V $ is separated 
if its topology is Hausdorff and 
complete if the topology of $ V $ is sequentially complete. \\
The bornology $ \mathfrak{Comp}(V) $ consists of all precompact subsets 
of $ V $. This means that $ S \in \mathfrak{Comp}(V) $ iff 
for all neighborhoods $ U $ of the origin there is a finite subset $ F \subset V $ 
such that $ S \subset F + U $. If $ V $ is complete then 
$ S \subset V $ is precompact iff its closure is compact. 
Equipped with the bornology $ \mathfrak{Comp}(V) $ the space $ V $ is separated if the 
topology of $ V $ is Hausdorff and complete if $ V $ is a complete topological 
vector space.  

\subsection*{Fr\'echet spaces}

In the case of Fr\'echet spaces the properties of the bounded bornology and the 
precompact bornology can be described more in detail. 
Let $ V $ and $ W $ be Fr\'echet spaces endowed both with the bounded or 
the precompact bornology. A linear map $ f: V \rightarrow W $ is bounded if and only if it is continuous. This is 
due to the fact that a linear map between metrizable topological spaces is continuous iff it is sequentially continuous. 
Hence the functors $ \mathfrak{Bound} $ and $ \mathfrak{Comp} $ 
from the category of Fr\'echet spaces into the category of complete bornological vector spaces are fully faithful. \\
The following theorem describes the completed bornological tensor product 
of Fr\'echet spaces with the precompact bornology and is proved in ~\cite{Meyer}. 
\begin{theorem} \label{Frechettensor}
Let $ V $ and $ W $ be Fr\'echet spaces and let 
$ V \hat{\otimes}_\pi W $ be their completed projective tensor product. 
Then there is a natural isomorphism
\begin{align*}
&(V,\mathfrak{Comp}) \hat{\otimes} (W,\mathfrak{Comp}) \cong 
(V \hat{\otimes}_\pi W,\mathfrak{Comp})
\end{align*}
of complete bornological vector spaces. 
\end{theorem}

\subsection*{LF-spaces}

More generally we can consider LF-spaces. A locally convex 
vector space $ V $ is an LF-space if there exists an increasing sequence 
of subspaces $ V_n \subset V $ with union equal to $ V $ such that each $ V_n $ is a Fr\'echet space in the subspace topology and $ V $ carries the 
corresponding inductive limit topology. A linear map 
$ V \rightarrow W $ from the LF-space $ V $ into an arbitrary locally convex 
space $ W $ is continuous iff its restriction to the subspaces $ V_n $ is continuous 
for all $ n $. From the definition of the inductive limit topology it follows that a bounded subset of an LF-space $ V $ is contained in 
a Fr\'echet subspace $ V_n $. If $ V_1 $ and $ V_2 $ are LF-spaces endowed 
with the bounded or the precompact bornology a bilinear map 
$ b: V_1 \times V_2 \rightarrow W $ is bounded iff it is separately continuous. 
This implies that an LF-space equipped with a separately continuous 
multiplication becomes a complete bornological algebra with respect to the 
bounded or the precompact bornology. \\
The following description of tensor products of LF-spaces can also be found in ~\cite{Meyer}. 
\begin{theorem}\label{LFtensor}
Let $ V $ and $ W $ be nuclear LF-spaces endowed with the bounded bornology. 
Then $ V \hat{\otimes} W $ is isomorphic to the inductive tensor product 
$ V \hat{\otimes}_\iota W $ endowed with the bounded bornology. 
\end{theorem}
Next we review the basic theory of smooth representations of locally compact groups on bornological 
vector spaces ~\cite{Meyersmoothrep}. In the sequel integration of functions on a locally 
compact group is always understood with respect to a fixed left Haar measure. \\
A representation of a locally compact group $ G $ on a complete bornological vector space $ V $ is a group homomorphism $ \pi: G \rightarrow 
\Aut(V) $ where $ \Aut(V) $ denotes the group of bounded linear automorphisms of $ V $. Let $ F(G,V) $ 
be the vector space of all functions from $ G $ to $ V $. The space  $ F(G,V) $ is simply 
the direct product of copies of the space $ V $ taken over the set $ G $. To a representation 
$ \pi: G \rightarrow \Aut(V) $ we associate the linear map $ [\pi]: V \rightarrow F(G,V) $ defined by 
$ [\pi](v)(t) = \pi(t)(v) $. 
\begin{definition} Let $ G $ be a locally compact group and let $ V $ be a complete bornological vector space. 
A representation $ \pi $ of $ G $ on $ V $ is smooth if $ [\pi] $ defines a bounded linear map from 
$ V $ into $ \E(G,V) $. A smooth representation is also called a $ G $-module. 
A bounded linear map $ f: V \rightarrow W $ between $ G $-modules is called equivariant if 
$ f(s \cdot v) = s \cdot f(v) $ for all $ v \in V $ and $ s \in G $.
\end{definition}
Here $ \E(G,V) $ denotes the space of smooth functions on $ G $ with values in $ V $. Smoothness has its 
usual meaning if $ G $ is a Lie group and $ V $ is a Banach space. If $ G $ is discrete any function from $ G $ to $ V $ is 
smooth. It follows that every representation of a discrete group is smooth. If $ G $ is totally disconnected and 
$ V $ is a fine space then a function from $ G $ to $ V $ is smooth iff it is locally constant. Hence 
for totally disconnected groups and fine spaces one recovers the ordinary theory of smooth representations on complex vector spaces. 
For the general definition of the space $ \E(G,V) $ and more information we refer to ~\cite{Meyersmoothrep}. \\
We denote by $ G \LMod $ the category of $ G $-modules and equivariant 
linear maps. The direct sum of a family of $ G $-modules is again a $ G $-module. 
The tensor product $ V \hat{\otimes} W $ of two $ G $-modules becomes 
a $ G $-module using the diagonal action $ s \cdot (v \otimes w) =
s \cdot v \otimes s \cdot w $ for $ v \in V $ and $ w \in W $. 
For every group the trivial one-dimensional $ G $-module $ \mathbb{C} $
is a unit with respect to the tensor product. In this way $ G\LMod $ becomes 
an additive monoidal category. \\  
Let $ \D(G) $ be the space of smooth functions with compact support on $ G $.  
For a Lie group $ G $ this is the space of smooth functions with compact support on $ G $ in the usual sense. 
If $ G $ is totally disconnected we obtain the space of locally constant functions on $ G $ with compact support. 
The group $ G $ acts on $ \D(G) $ by left translations 
$$
(s \cdot f)(t) = f(s^{-1}t)
$$
and $ \D(G) $ becomes a $ G $-module in this way. \\
A $ G $-module is called projective if it has the lifting property with respect to equivariant surjections 
$ M \rightarrow N $ of $ G $-modules with bounded linear splitting $ N \rightarrow M $. 
\begin{lemma} \label{omegaprojective}
Let $ V $ be any $ G $-module. Then the $ G $-module $ \D(G) \cotimes V $ is projective.
\end{lemma}
\proof 
We use a standard argument ~\cite{Blanc}. Let $ \pi: M \rightarrow N $ be a surjective equivariant map with a 
bounded linear splitting $ \sigma $. Morever let $ \phi: D(G) \cotimes V \rightarrow M $ be any equivariant linear map.
Choose a function $ \chi \in \D(G) $ such that 
$$
\int_G \chi(s) ds = 1 
$$
and define 
\begin{equation*}
f_s(t) = f(t) \chi(t^{-1}s)
\end{equation*}
for every $ f \in \D(G) $ and $ s \in G $. Then one computes 
\begin{equation*}
\int_G f_s(t) ds = f(t)
\end{equation*}
and $ t \cdot (f_{t^{-1}s}) = (t \cdot f)_s $ for all $ f \in \D(G) $ and $ s,t \in G $. 
We set   
\begin{equation*}
\psi(f \otimes v) = 
\int_G t \cdot \sigma \phi(t^{-1} \cdot (f_t \otimes v)) dt. 
\end{equation*}
Since we have $ t^{-1} \cdot(f_t) = ({t^{-1}} \cdot f)_e $
the integral is well-defined. It is easy to check that $ \psi $ extends to an equivariant linear
map $ \D(G) \cotimes V \rightarrow M $. Finally  we have 
\begin{align*}
\pi \psi(f \otimes v) = 
\int_G t \cdot \pi\sigma \phi(t^{-1} \cdot (f_t \otimes v)) dt = \int_G \phi(f_t \otimes v) dt = \phi(f \otimes v)
\end{align*}
using that $ \pi $ and $ \phi $ are equivariant. This yields the assertion. \qed \\
Next we specify the class of $ G $-algebras we are going to work with. Expressed in the language of 
category theory our definition amounts 
to saying that a $ G $-algebra is an algebra in the monoidal category $ G \LMod $. 
\begin{definition}
Let $ G $ be a locally compact group. A $ G $-algebra is a complete 
bornological algebra $ A $ which is at the same time a $ G $-module such that the multiplication satisfies
\begin{equation*}
s \cdot (xy) = (s \cdot x) (s \cdot y)
\end{equation*}
for all $ x,y \in A $ and $ s \in G $. An equivariant homomorphism
$ f: A \rightarrow B $ between $ G $-algebras is an
algebra homomorphism which is equivariant. 
\end{definition} 
If $ A $ is unital we say that $ A $ is a unital $ G $-algebra if
$ s \cdot 1 = 1 $ for all $ s \in G $. 
The unitarisation $ A^+ $ of a $ G $-algebra $ A $ is a unital $ G $-algebra in 
a natural way.
We will occasionally also speak of an action of 
$ G $ on $ A $ to express that $ A $ is a $ G $-algebra. \\
There is a natural way to enlarge any $ G $-algebra to a $ G $-algebra where all 
group elements act by inner automorphisms. This is the crossed product
construction which we study next.
\begin{definition}
Let $ G $ be a locally compact group and let $ A $ be a $ G $-algebra. The crossed
product $ A \rtimes G $ of $ A $ by $ G $ is 
$ A \hat{\otimes} \D(G) = \D(G,A) $ with multiplication given by
\begin{equation*}
(f * g)(t) = \int_G f(s) s \cdot g(s^{-1}t) ds
\end{equation*}
for $ f, g \in \D(G,A) $. 
\end{definition}
It is easy to check that $ A \rtimes G $ is a complete bornological algebra. 
If we consider the case $ A = \mathbb{C} $ with the trivial action we 
obtain by definition the smooth group algebra $ \D(G) $ of $ G $. If $ G $ is discrete this is simply the 
complex group ring $ \mathbb{C}G $ endowed with the fine bornology. \\
In general the crossed product does not posses a unit, the algebra 
$ A \rtimes G $ is unital if $ A $ has a unit and $ G $ is discrete. 
We want to show that the crossed product $ A \rtimes G $ still has an approximate identity 
whenever $ A $ has one. Let us first recall from ~\cite{Meyersmoothrep} the concept of an 
approximate identity. A complete bornological algebra $ A $ is said to have an approximate identity 
if for any bornologically compact subset $ S \subset A $ there exists 
a sequence $ (u_n)_{n \in \mathbb{N}} $ in $ A $ such that 
$ u_n \cdot a - a $ and $ a \cdot u_n - a $ converge to zero uniformly for 
$ a \in S $. A subset of a bornological vector space $ V $ is bornologically compact 
if it is a compact subset of the Banach space $ \bra T \ket $ for some 
completant small disk $ T \subset V $. Uniform convergence means that there 
exists a completant small disk $ T \subset A $ such  that the sequences $ u_n \cdot a - a $ and 
$ a \cdot u_n - a $ converge uniformly to zero in the Banach space $ \bra T \ket $. \\
An $ A $-module $ M $ over a bornological algebra $ A $ with approximate identity is called 
nondegenerate if the module action $ A \cotimes M \rightarrow M $ 
is a bornological quotient map. This is equivalent to saying that the natural map 
$ A \cotimes_A M \rightarrow M $ is a bornological isomorphism ~\cite{Meyersmoothrep}. \\
Given a smooth representation $ \pi $ of $ G $ on $ V $ one defines a $ \D(G) $-module structure on 
$ V $ by setting 
$$
f \cdot v = \int_G f(t) \, t \cdot v dt. 
$$
It is shown in ~\cite{Meyersmoothrep} that the smooth group algebra $ \D(G) $ has an approximate 
identity and that the previous construction defines an isomorphism between the category of smooth representations of $ G $ 
and the category of nondegenerate $ \D(G) $-modules for every locally compact group $ G $. \\
We have the following extension of proposition 4.3 in ~\cite{Meyersmoothrep}. 
\begin{prop}
Let $ G $ be a locally compact group and let $ A $ be a $ G $-algebra with approximate identity. 
Then the crossed product $ A \rtimes G $ has an approximate identity. 
\end{prop}
\proof The idea is to combine the approximate identity of $ A $ with an approximate 
identity for $ \D(G) $, the latter being constructed in \cite{Meyersmoothrep}. 
In the sequel we will view elements of $ A $ and $ \D(G) $ as left and right multipliers 
of the crossed product $ A \rtimes G $ in the obvious way. Let $ S \subset A \rtimes G $ be a 
bornologically compact subset. Right multiplication of $ \D(G) $ on $ A \rtimes G $ does not 
involve $ A $. Let us consider left multiplication. 
Since $ A $ is a smooth representation, the left action of $ G $ on $ A \rtimes G $ is smooth. 
Hence there exists a bounded linear splitting $ \sigma: A \rtimes G \rightarrow \D(G) \cotimes (A \rtimes G) $ 
for the left action of $ \D(G) $ on the crossed product. Clearly the image $ \sigma(S) $ of $ S $ is again 
bornologically compact. 
Using Grothendieck's result about compact subsets of the projective tensor product of Fr\'echet spaces 
\cite{Grothendieck} we see that $ S $ is contained in the completant disked hull of 
$ R_l \otimes C_l $ for bornologically compact subsets $ R_l \subset A $ and $ C_l \subset \D(G) $. 
Similarly, $ \sigma(S) $ is contained in the completant disked hull of 
$ C_r \otimes R_r $ for bornologically compact subsets $ C_r \subset \D(G) $ and $ R_r \subset A \rtimes G $. 
Hence we obtain a sequence $ (h_n)_{n \in \mathbb{N}} $ in $ \D(G) $ such that $ f \cdot h_n - f $ 
and $ h_n \cdot \sigma(f) - \sigma(f) $ converge to zero uniformly for $ f \in S $. After applying 
the multiplication map $ \D(G) \cotimes (A \rtimes G) \rightarrow A \rtimes G $ we see that 
$ h_n \cdot f - f $ converges uniformly to zero in $ A \rtimes G $.  \\
Left multiplication of $ A $ on $ A \rtimes G $ does not involve $ \D(G) $. For right 
multiplication the explicit formula is $ (f \cdot a)(t) = f(t) (t \cdot a) $ for 
$ f \in \D(G,A) $ and $ a \in A $. Let $ \phi: A \rtimes G \rightarrow A \cotimes \D(G) = \D(G,A) $ 
be the isomorphism given by $ \phi(f)(t) = t^{-1} \cdot f(t) $. Then the right 
action of $ A $ on $ A \rtimes G $ corresponds under the map $ \phi $ to the trivial right action 
$ (f \cdot a)(t) = f(t) a $ on $ A \cotimes \D(G) $. 
As above we choose a sequence $ (a_n)_{n \in \mathbb{N}} $ in 
$ A $ such that $ a_n \cdot f - f $ and $ \phi(f) \cdot a_n - \phi(f) $ converge 
uniformly to zero for all $ f \in S $. Then $ f \cdot a_n - f $ converges 
uniformly to zero in $ A \rtimes G $ for all $ f \in S $. 
Define $ u_n = a_n \otimes h_n \in A \rtimes G $. Using the 
equations
$$
u_n \cdot f - f = a_n \cdot (h_n \cdot f - f) + (a_n \cdot f - f)
$$ 
and 
$$
f \cdot u_n - f = (f \cdot a_n - f)h_n + (f\cdot h_n - f)
$$
we see that $ u_n \cdot f - f $ and $ f \cdot u_n - f $ converge to zero uniformly for $ f \in S $. 
Hence $ A \rtimes G $ has an approximate identity. \qed 
\begin{definition}
A covariant representation of a $ G $-algebra $ A $ with approximate identity is 
a complete bornological vector space $ M $ which is both a $ G $-module and a 
nondegenerate $ A $-module such that 
$$
s \cdot (a \cdot m) = (s \cdot a) \cdot (a \cdot m)
$$
for all $ s \in G, f \in A $ and $ m \in M $. A bounded linear map $ f: M \rightarrow N $ between 
covariant representations is covariant if it is $ A $-linear and equivariant. 
\end{definition}
Clearly covariant representations of a $ G $-algebra $ A $ and covariant maps form a category. 
The next result shows that this category is closely relate to the crossed product construction. 
\begin{prop}\label{covrepprop} Let $ A $ be a $ G $-algebra with an approximate identity. Then the 
category of nondegenerate $ A \rtimes G $-modules is isomorphic to the category of covariant representations of $ A $. 
\end{prop}
\proof Let $ M \cong (A\rtimes G) \cotimes_{A \rtimes G} M $ be a nondegenerate $ A \rtimes G $-module. 
Then we obtain a representation of $ G $ and an $ A $-module structure on $ M $ by 
letting act $ s \in G $ and $ a \in A $ as left multipliers on $ A \rtimes G $. 
Since the action of $ G $ on $ A \rtimes G $ is smooth we have natural isomorphisms
$$ 
\D(G) \cotimes_{\D(G)} M \cong \D(G) \cotimes_{\D(G)} (A \rtimes G) \cotimes_{A \rtimes G} M \cong 
(A \rtimes G) \cotimes_{A \rtimes G} M \cong M 
$$ 
for the integrated 
form of this representation of $ G $ and it follows that $ M $ becomes 
a $ G $-module. Moreover we have 
$$
A \cotimes_A M \cong  A \cotimes_A (A \rtimes G) \cotimes_{A \rtimes G} \cotimes M \cong 
(A \rtimes G) \cotimes_{A \rtimes G} M \cong M 
$$
in a natural way using the fact that multiplication induces an isomorphism $ A \cotimes_A A \cong A $ 
due to the existence of an approximate identity for $ A $. It follows that  $ M $ is a 
nondegenerate $ A $-module. In this way $ M $ becomes a covariant representation. \\
Conversely, assume that $ M $ is a covariant representation of $ A $. Then we obtain 
an $ A \rtimes G $-module structure on $ M $ by setting 
$$
f \cdot m = \int_G f(t) (t \cdot m) dt 
$$
for $ f \in \D(G,A) $. The module structure $ \mu: (A \rtimes G) \cotimes M \rightarrow M $ 
can be decomposed as 
\begin{equation*} 
(A \rtimes G) \cotimes M \;= 
 \xymatrix{
     A \cotimes \D(G) \cotimes M \; \ar@{->}[r]^{\quad \; \; \; \id \cotimes \mu_G} &
         A \cotimes M  \ar@{->}[r]^{\;\;\mu_A} & M   
     }
\end{equation*}
where $ \mu_G: \D(G) \cotimes M \rightarrow M $ and $ \mu_A: A \cotimes M \rightarrow M $ are 
the given module structures. Since $ M $ is a $ G $-module the map $ \mu_G $ has a bounded linear 
splitting. Hence the first arrow is a bornological quotient map. Moreover $ \mu_A $ is a bornological 
quotient map since $ M $ is a nondegenerate $ A $-module. It follows that $ M $ is a nondegenerate 
$ A \rtimes G $-module. \\
The previous constructions are compatible with morphisms and it is easy to see that they are 
inverse to each other. This yields the assertion. \qed \\
Let us have a look at some basic examples of $ G $-algebras and the associated 
crossed products. In particular the algebra 
$ \mathcal{K}_G $ introduced below will play an important role in our theory. 

\subsection*{Trivial actions}

The simplest example of a $ G $-algebra is the algebra of complex
numbers with the trivial $ G $-action. More generally one can equip any
complete bornological algebra $ A $ with the trivial action to obtain a 
$ G $-algebra. The corresponding crossed product algebra $ A \rtimes G $ is
simply a tensor product,
\begin{equation*}
A \rtimes G \cong A \hat{\otimes} \D(G).
\end{equation*}
This explains why one may view crossed products in general as
twisted tensor products.

\subsection*{Commutative algebras}

Let $ M $ be a smooth manifold on which the Lie group $ G $ acts
smoothly  and let $ C^\infty_c(M) $ be the LF-algebra of 
compactly supported smooth functions on $ M $. 
Then we get an action of $ G $ on $ A = C^\infty_c(M) $ by defining
\begin{equation*}
(s \cdot f)(x) = f(s^{-1} \cdot x)
\end{equation*}
for all $ s \in G $ and $ f \in A $. This algebra is
unital if $ M $ is compact and $ G $ is discrete.
The associated crossed product $ A \rtimes G $ may be described as 
the smooth convolution algebra
of the translation groupoid $ M \rtimes G $ associated to the action
of $ G $ on $ M $.

\subsection*{Algebras associated to representations of $ G $}

Let $ V $ and $ W $ be $ G $-modules and let 
$ b: W \times V \rightarrow \mathbb{C} $ be an equivariant bounded bilinear map. Then 
$ l(b) = V \cotimes W $ is a $ G $-algebra with the multiplication 
\begin{equation*}
(v_1 \otimes w_1) \cdot (v_2 \otimes w_2) = v_1 \otimes b(w_1, v_2) w_2 
\end{equation*}
and the diagonal $ G $-action. \\ 
In the case $ V = W $ we have a natural homomorphism $ l(b) \rightarrow \End(V) $ 
given by 
$$ 
\iota(v \otimes w)(u) = v\; b(w, u).
$$
If we equip $ \End(V) $ with the representation of $ G $ defined by 
the formula
\begin{equation*}
(s \cdot T)(u) = s \cdot T(s^{-1} \cdot u)
\end{equation*}
for $ s \in G $ and $ u \in V  $ the homomorphism $ \iota $ becomes equivariant. \\
A basic example is given by the left regular representation on $ \D(G) $.  
We set $ V = W = \D(G) $ and consider the pairing 
\begin{equation*}
b(f, g) = \int_G f(t) g(t) dt.
\end{equation*}
The corresponding $ G $-algebra will be denoted by $ \mathcal{K}_G $. 
Elements in $ \mathcal{K}_G $ can be viewed as kernels $ k \in \D(G \times G) $ 
of integral operators acting on $ \D(G) $ by 
\begin{equation*}
(k f)(s) = \int_G k(s,t) f(t) dt. 
\end{equation*}
Finally observe that by lemma \ref{omegaprojective} the tensor product $ V \hat{\otimes} \mathcal{K}_G $ is a 
projective $ G $-module for every $ G $-module $ V $. 

\section{Covariant modules} 
\label{seccov}

In this section we introduce the notion of a covariant modules which plays an important role 
in equivariant cyclic homology. \\
Let $ G $ be a locally compact group. Then $ G $ can be viewed as a $ G $-space using the adjoint action. 
This induces an action of $ G $ on $ \D(G) $ viewed as a commutative algebra with pointwise multiplication. 
The resulting $ G $-algebra will be denoted by $ \mathcal{O}_G $ in order to distinguish it from 
the smooth group algebra of $ G $. Explicitly we have $ (t \cdot f)(s) = f(t^{-1} s t) $ for
$ f \in \mathcal{O}_G $ and $ s \in G $. It is evident that the algebra $ \mathcal{O}_G $ has an approximate identity. 
Remark that $ \mathcal{O}_G $ is unital iff the group $ G $ is compact. \\
We are interested in covariant representations of this particular $ G $-algebra
and give the following explicit definition. 
\begin{definition} Let $ G $ be a locally compact group. A (smooth) $ G $-covariant
module is a complete bornological vector space $ M $ which is 
both a nondegenerate $ \mathcal{O}_G $-module and a $ G $-module such that
\begin{equation*}
s \cdot (f \cdot m) = (s \cdot f) \cdot (s \cdot m)
\end{equation*}
for all $ s \in G, f \in \mathcal{O}_G $ and $ m \in M $. A bounded linear map
$ \phi: M \rightarrow N $ between covariant modules is called covariant
if it is $ \mathcal{O}_G $-linear and equivariant.
\end{definition}
We remark that covariant modules may be thought of as spaces of global sections of equivariant sheaves over 
$ G $ viewed as a $ G $-space with the adjoint action. Moreover, due to proposition \ref{covrepprop} a covariant module 
is the same thing as a nondegenerate module over the crossed product $ \mathcal{O}_G \rtimes G $. In the sequel 
we will also write $ \Cov(G) $ for the crossed product $ \mathcal{O}_G \rtimes G $. \\  
Usually we will not mention the group explicitly in our terminology and simply speak of covariant modules and covariant maps. 
The category of covariant modules and covariant maps will be denoted by $ G\LSMod $ 
and we will write $ \SHom_G(M,N) $ for the space of covariant maps between covariant modules $ M $ and $ N $. 
In addition we let $ \SHom(M,N) $ be the collection of maps that are only $ \mathcal{O}_G $-linear. \\
A basic example of a covariant module is the algebra
$ \mathcal{O}_G $ itself. More generally, let $ V $ be a $ G $-module.
We obtain an associated covariant module by
considering $ \mathcal{O}_G \hat{\otimes} V $ with the diagonal $ G $-action
and the obvious $ \mathcal{O}_G $-module structure given by multiplication. In the case 
$ V = \D(G) $ we obtain just $ \Cov(G) $ viewed as a left module over itself. 
If $ V $ is any $ G $-module then $ \Cov(G) \cotimes V $ 
becomes a covariant module by the diagonal action of $ G $ and left multiplication of $ \mathcal{O}_G $. \\
Let us consider the covariant module $ \Cov(G) $. We can view elements in $ \Cov(G) $ as 
smooth functions with compact support on $ G \times G $ where the first variable corresponds to $ \mathcal{O}_G $ and 
the second variable corresponds to $ \D(G) $. The multiplication in the crossed product becomes 
$$
(f \cdot g)(s,t) = \int_G f(s,r)g(r^{-1}sr, r^{-1} t) dr 
$$
in this picture.  
\begin{lemma} \label{Toplemma} The bounded linear map $ T: \Cov(G) \rightarrow \Cov(G) $ defined by 
$$
T(f)(s,t) = f(s,st)
$$
is an isomorphism of $ \Cov(G) $-bimodules. 
\end{lemma}
\proof It is clear that $ T $ is a bounded linear isomorphism with inverse given by 
$ T^{-1}(f)(s,t) = f(s,s^{-1}t) $. We compute 
\begin{align*}
(f \cdot T(g))(s,t) &= \int_G f(s,r) T(g)(r^{-1}sr, r^{-1}t) dr \\
&= \int_G f(s,r) g(r^{-1}sr, r^{-1}st) dr = T(f \cdot g)(s,t) \\
&= \int_G f(s, sr) g(r^{-1}sr, r^{-1}t) dr = (T(f) \cdot g)(s,t)
\end{align*}
for $ f, g \in \Cov(G) $. This proves the assertion. \qed \\ 
Now consider an arbitrary covariant module $ M $. Since $ \Cov(G) $ has an approximate 
identity we have a natural isomorphism $ M \cong \Cov(G) \cotimes_{\Cov(G)} M $. 
Let us define $ T: M \rightarrow M $ by 
$$
T(f \otimes m) = T(f) \otimes m 
$$
for $ f \otimes m \in \Cov(G) \otimes_{\Cov(G)} M $. It follows from lemma \ref{Toplemma} that this definition makes 
sense. The operator $ T $ has the following fundamental properties.  
\begin{prop} \label{covparaadd} The operator $ T: M \rightarrow M  $ is a covariant isomorphism for all covariant modules $ M $. 
If $ \phi: M \rightarrow N $ is any covariant map between covariant modules then we have $ T \phi = \phi T $. Hence $ T $ defines a natural 
isomorphism $ T: \id \rightarrow \id $ of the identity functor $ \id: G\LSMod \rightarrow G\LSMod $. 
\end{prop}
\proof It is clear from lemma \ref{Toplemma} that $ T: M \rightarrow M $ is a covariant isomorphism for all $ M $. 
Using the fact that $ M $ and $ N $ are nondegenerate $ \Cov(G) $-modules the equation $ T \phi = \phi T $ follows easily 
after identifying $ \phi $ with the covariant map $ \id \cotimes \phi: \Cov(G)\cotimes_{\Cov(G)} M 
\rightarrow \Cov(G)\cotimes_{\Cov(G)} N $. The last statement is just a reformulation of the first two assertions. \qed \\
We conclude this section by exhibiting certain projective objects in the category of covariant modules. 
A covariant module $ P $ is projective if for every covariant map 
$ \pi: M \rightarrow N $ with a bounded linear 
splitting $ \sigma: N \rightarrow M $ between covariant modules and every covariant map 
$ \phi: P \rightarrow N $ 
there exists a covariant map $ \psi: P \rightarrow M $ such that 
$ \pi \psi = \phi $. 
\begin{lemma} \label{covmodproj} Let $ V $ be any $ G $-module. Then the covariant module $ \Cov(G) \cotimes V $ is projective. 
\end{lemma}
\proof Let $ \pi: M \rightarrow N $ be a surjective covariant map  
with bounded linear splitting $ \sigma: N \rightarrow M $ and let  
$ \phi: \Cov(G) \cotimes V \rightarrow N $ be any covariant map. Moreover let $ (\chi_j)_{j \in J} $ 
be a partition of unity for $ G $ with $ \chi_k \in \D(G) $ for all $ k $ such that 
$ \sum_{j \in J} \chi_j^2 = 1 $. We define a bounded linear map $ \eta: \Cov(G) \cotimes V \rightarrow  M $ as follows. For 
$ f \otimes g \otimes v \in \mathcal{O}_G \otimes \D(G) \otimes V $ set 
\begin{equation*}
\eta(f \otimes g \otimes v) = 
\sum_{j \in J} (f \chi_j)\cdot \sigma \phi(\chi_j \otimes g \otimes v) 
\end{equation*}
and observe that the sum is actually finite since the support of $ f $ is compact 
for every $ f \in \mathcal{O}_G $. It is easy to check that $ \eta $ extends to 
the completion $ \Cov(G) \cotimes V $. Moreover it follows from the definitions that $ \eta $ is $ \mathcal{O}_G $-linear 
and that we have $ \pi \eta = \phi $. \\
With the same notation as in the proof of lemma \ref{omegaprojective} we set  
\begin{equation*}
\psi(f \otimes g \otimes v) = 
\int_G t \cdot \eta(t^{-1} \cdot (f \otimes g_t \otimes v)) dt
\end{equation*}
for an element $ f \otimes g \otimes v \in \mathcal{O}_G \otimes \D(G) \otimes V $. 
One checks that $ \psi $ extends to a bounded linear 
map $ \Cov(G) \cotimes V \rightarrow M $. Moreover $ \psi $ is $ \mathcal{O}_G $-linear and 
equivariant. Finally  
one computes $ \pi \psi = \phi $
using that $ \pi \eta = \phi $ is covariant. This yields the assertion. \qed 

\section{Projective systems}
\label{secprocat}

The most natural way to define equivariant periodic cyclic homology is to work in the category
of pro-$ G $-algebras. This means that we have to consider projective systems of 
$ G $-modules and covariant modules. In this section we review these notions 
and fix our notation. \\
To any additive category $ \mathcal{C} $ one associates the pro-category $ \pro(\mathcal{C}) $ of projective 
systems over $ \mathcal{C} $ as follows. A projective system over $ \mathcal{C} $ consists of a directed
index set $ I $, objects $ V_i $ for all $ i \in I $ and 
morphisms $ p_{ij}: V_j \rightarrow V_i $ for
all $ j \geq i $. The morphisms are assumed to satisfy
$ p_{ij} p_{jk} = p_{ik} $ if $ k \geq j \geq i $. These conditions
are equivalent to saying that we have a contravariant functor from
the small category $ I $ to $ \mathcal{C} $. The class of
objects of $ \pro(\mathcal{C}) $ consists by definition of all
projective systems over $ \mathcal{C} $. The space of morphisms between projective systems
$ (V_i)_{i \in I} $ and $ (W_j)_{j \in J} $ is defined by 
\begin{equation*}
\Mor((V_i),(W_j)) = \varprojlim_{j} \varinjlim_i
\Mor_\mathcal{C}(V_i,W_j)
\end{equation*}
where the limits are taken in the category of abelian groups. Of 
course one has to check that the composition of morphisms can be defined
in a consistent way. We refer to ~\cite{AM} for further details. \\
It is useful to study pro-objects by comparing them to constant pro-objects. A constant pro-object is by definition a pro-object 
where the index set consists only of one element. If $ V = (V_i)_{i \in I} $ is any pro-object a morphism 
$ V \rightarrow C $ with constant range $ C $ is given by a 
morphism $ V_i \rightarrow C $ for some $ i $. \\ 
In the category $ \pro(\mathcal{C}) $ projective limits always exist. 
This is due to the fact that a projective system of pro-objects $ (V_j)_{j \in J} $ 
can be identified naturally with a pro-object. \\
Since there are finite direct sums in $ \mathcal{C} $ we also have finite direct sums in $ \pro(\mathcal{C}) $. Explicitly, 
the direct sum of $ V = (V_i)_{i \in I} $ and 
$ W =  (W_j)_{j \in J} $ is given by 
\begin{equation*}
(V_i)_{i \in I} \oplus (W_j)_{j \in J} =
(V_i \oplus W_j)_{(i,j) \in I \times J}
\end{equation*}
where the index set $ I \times J $ is ordered using the product
ordering. The structure maps of this projective system are obtained 
by taking direct sums of the structure maps of 
$ (V_i)_{i \in I} $ and $ (W_j)_{j \in J} $. With this notion of direct sums the category $ \pro(\mathcal{C}) $ becomes an additive category. \\
If we apply these general constructions to the category of $ G $-modules we obtain 
the category of pro-$ G $-modules.
A morphism in $ \pro(G\LMod) $ will be called an equivariant linear
map. Similarly we have the category of covariant pro-modules as the pro-category 
of $ G\LSMod $. Morphisms in $ \pro(G\LSMod) $ will be called covariant maps. \\
Let us come back to the general situation. Assume in addition that 
$ \mathcal{C} $ is monoidal such that the tensor product functor $ \mathcal{C} \times \mathcal{C} \rightarrow \mathcal{C} $ is bilinear. 
In this case we define the tensor product $ V \otimes W $ for pro-objects 
$ V = (V_i)_{i \in I} $ and $ W =  (W_j)_{j \in J} $ by
\begin{equation*}
(V_i)_{i \in I} \otimes (W_j)_{j \in J} =
(V_i \otimes W_j)_{(i,j) \in I \times J}
\end{equation*}
where again $ I \times J $ is ordered using the product
ordering. The structure maps are obtained 
by tensoring the structure maps of $ (V_i)_{i \in I} $ and $ (W_j)_{j \in J} $. 
Observe that any morphism $ f: V \otimes W \rightarrow C $ with constant range $ C $ factors through $ V_i \otimes W_j $ for some 
$ i \in I, j \in J $. This means that we can write $ f $ in the form 
$ f = g(f_V \otimes f_W) $ where $ f_V: V \rightarrow C_V $ and 
$ f_W: W \rightarrow C_W $ are morphisms with constant range and 
$ g: C_V \otimes C_W \rightarrow W $ is a morphism of constant 
pro-objects. \\
Equipped with this tensor product the category $ \pro(\mathcal{C}) $ is 
additive monoidal and we obtain a natural faithful additive monoidal functor $ \mathcal{C} \rightarrow \pro(\mathcal{C}) $. \\
The existence of a tensor product in $ \pro(\mathcal{C}) $ yields a natural notion 
of algebras and algebra homomorphisms in this category. Such algebras will be
called pro-algebras and their homomorphism will be called pro-algebra homomorphisms. Moreover we can 
consider pro-modules for pro-algebras and their homomorphisms. \\
The category $ G\LMod $ is monoidal in the sense explained above. To indicate that 
we use completed bornological tensor products in $ G\LMod $ we will denote the tensor product of two pro-$ G $-modules $ V $ and 
$ W $ by $ V \hat{\otimes} W $. \\
In order to fix terminology we give the following definition. 
\begin{definition} A pro-$ G $-algebra $ A $ is an algebra in the category 
$ \pro(G \LMod) $. An algebra homomorphism $ f: A \rightarrow B $ in 
$ \pro(G \LMod) $ is called an equivariant homomorphism of 
pro-$ G $-algebras. 
\end{definition} 
Occasionally we will consider unital pro-$ G $-algebras. The unitarisation 
$ A^+ $ of a pro-$ G $-algebra $ A $ is defined in the same way as for 
$ G $-algebras. \\
We also include a short discussion of extensions. Let again $ \mathcal{C} $ be 
any additive category and let $ K, E $ and $ Q $ be objects in 
$ \pro(\mathcal{C}) $. A (strict) extension is a diagram of the form 
\begin{equation*}
   \xymatrix{
     K \ar[r]^{\iota} & E \ar@/^/@{.>}[l]^{\rho} \ar[r]^{\pi} & Q \ar@/^/@{.>}[l]^{\sigma}
     }
\end{equation*}
in $ \pro(\mathcal{C}) $ such that $ \rho \iota = \id $, $ \pi \sigma = \id $ and $ \iota \rho + \sigma \pi = \id $. 
In other words we require that $ E $ decomposes into a direct sum of 
$ K $ and $ Q $. We will frequently omit the splitting $ \sigma $ and the retraction $ \rho $ in our notation and 
write simply 
\begin{equation*}
   \xymatrix{
     K \;\; \ar@{>->}[r]^{\iota} & E \ar@{->>}[r]^{\pi} & Q 
     }
\end{equation*}
or $ (\iota,\pi): 0 \rightarrow K \rightarrow E \rightarrow Q 
\rightarrow 0 $ for an extension. \\
Let us give the following definition in the situation $ \mathcal{C} = \pro(G\LMod) $.
\begin{definition}
Let $ K, E $ and $ Q $ be pro-$ G $-algebras. An extension of pro-$ G $-algebras is an extension 
\begin{equation*}
   \xymatrix{
     K \;\; \ar@{>->}[r]^{\iota} & E \ar@{->>}[r]^{\pi} & Q
     }
\end{equation*}
in $ \pro(G \LMod) $ where $ \iota $ and $ \pi $ are equivariant algebra homomorphisms. 
\end{definition}
Later we will need the concept of relatively projective pro-$ G $-modules and covariant pro-modules. 
A pro-$ G $-module $ P $ is called relatively projective if for every equivariant linear map $ \pi: M \rightarrow N $ 
of pro-$ G $-modules with pro-linear section $ N \rightarrow M $ and every equivariant linear map $ \phi: P \rightarrow N $ 
there exists an equivariant linear map $ \psi: P \rightarrow M $ such that $ \pi \psi = \phi $. Similarly a covariant 
pro-module is called relatively projective if it has the lifting property with respect to covariant maps 
between covariant pro-modules having a pro-linear section. 
The following lemma gives a simple criterion for relative projectivity. 
\begin{lemma}\label{relproj} Let $ V $ be a pro-$ G $-module. Then $ \D(G) \cotimes V $ is a relatively projective pro-$ G $-module 
and $ \Cov(G) \cotimes V $ is a relatively projective covariant pro-module. 
\end{lemma}
\proof This follows from the fact that the constructions in the proofs of lemma \ref{omegaprojective} and lemma \ref{covmodproj} are natural. \qed \\
Working with pro-$ G $-modules or covariant pro-modules may seem 
somewhat difficult because there are no longer concrete elements 
to manipulate with. Nevertheless we will write down explicit formulas 
involving ``elements'' in subsequent sections. This can be justified 
by noticing that these formulas are concrete expressions for identities between abstractly defined morphisms. 

\section{Paracomplexes}
\label{secpara}

In this section we introduce the concept of a paramixed complex. Our terminology is motivated from ~\cite{GJ1} but it is 
slightly different. The related notion of a paracyclic module is well-known in the study of the cyclic 
homology of crossed products and smooth groupoids \cite{FT}, \cite{GJ1}, \cite{Nistor1}, \cite{Crainic}. \\
Whereas cyclic modules and mixed complexes are fundamental concepts in cyclic homology, paracyclic modules 
are mainly regarded as a tool in computations. However, 
in the equivariant situation the point of view has to be changed drastically. Here the fundamental objects 
are paramixed complexes and mixed complexes show up mainly in calculations. \\
In abstract terms our notion of a paracomplex can be defined most naturally
using the concept of a para-additive category.  
\begin{definition} A para-additive category is an additive category $ \mathcal{C} $ 
together with a natural isomorphism $ T $ of the identity functor $ \id: \mathcal{C} \rightarrow \mathcal{C} $.   
\end{definition}
In other words, we are given invertible morphisms $ T(M): M \rightarrow M $ for all 
objects $ M \in \mathcal{C} $ such that $ \phi  T(M) = T(N) \phi $ for all morphisms $ \phi: M \rightarrow N $. In the sequel we will 
simply write $ T $ instead of $ T(M) $. \\
Clearly any additive category is para-additive by setting $ T = \id $. More interestingly, 
it follows from proposition \ref{covparaadd} that the category $ G\LSMod $ of covariant modules for a locally compact group $ G $ 
is a para-additive category in a natural way. Remark that in this case the operator $ \id - T: M \rightarrow M  $ is usually far from being 
zero. 
\begin{definition} Let $ \mathcal{C} $ be a para-additive category. 
A paracomplex $ C = C_0 \oplus C_1 $ in $ \mathcal{C} $ is a given by objects $ C_0 $ and $ C_1 $ 
in $ \mathcal{C} $ together with  morphisms $ \partial_0: C_0 \rightarrow C_1 $ and $ \partial_1: C_1 \rightarrow C_0 $ such that
$$ 
\partial^2 = \id - T 
$$
where the differential $ \partial: C \rightarrow C_1 \oplus C_0 \cong C $ is the composition of $ \partial_0 \oplus \partial_1 $
with the canonical flip map. A chain map $ \phi: C \rightarrow D $ between two paracomplexes is a morphism from $ C $ to $ D $ that commutes with
the differentials.
\end{definition}
Remark that we consider only $ \mathbb{Z}_2 $-graded objects. The morphism $ \partial $ in a paracomplex is called a differential although 
this contradicts the classical definition of a differential. \\
In general it does not make sense to speak about the homology of a
paracomplex. Given a paracomplex $ C $ with differential $ \partial $, for instance in a category of modules over some ring, one could force 
it to become a complex by dividing out the subspace $ \partial^2(C) $ and then take homology. However, it turns out that this procedure is not 
appropriate in our context. \\
Although there is no reasonable definition of homology we can give meaning to the statement 
that two paracomplexes are homotopy equivalent: Let $ \phi, \psi: C \rightarrow D $ be two chain maps between paracomplexes. A
chain homotopy connecting $ \phi $ and $ \psi $ is a map $ \sigma: C
\rightarrow D $ of degree $ 1 $ satisfying the usual relation $ \partial \sigma
+ \sigma \partial  = \phi - \psi $. Note that the map $ \partial \sigma + \sigma \partial $ is a chain map for any 
morphism $ \sigma: C \rightarrow D $ of odd degree since $ \partial^2 $ commutes with all morphisms in $ \mathcal{C} $. 
Two paracomplexes $ C $ and $ D $ are called homotopy equivalent if there exist chain maps $ \phi: C \rightarrow D $ and 
$ \psi: D \rightarrow C $ which are inverse to each other up to chain homotopy. \\
The paracomplexes we have in mind arise from paramixed complexes that we are going to define now.
\begin{definition} Let $ \mathcal{C} $ be a para-additive category.
A paramixed complex $ M $ in $ \mathcal{C} $ is a sequence of
objects $ M_n $ together with differentials $ b $ of degree $ -1 $ and $ B $ of degree $ + 1 $ satisfying
$ b^2 = 0 $, $ B^2 = 0 $ and
\begin{equation*}
[b,B] = bB + Bb = \id - T.
\end{equation*}
\end{definition}
If $ \mathcal{C} $ is additive, that is $ T = \id $, we reobtain the notion of a mixed complex. 
In general one can define and study Hochschild homology of a paramixed complex in the usual way since the 
Hochschild operator $ b $ satisfies $ b^2 = 0 $. On the other hand we shall not try to define the cyclic homology of 
an arbitrary paramixed complex. We will see below how bivariant periodic cyclic homology can still be defined in a 
natural way. 

\section{Quasifree pro-$ G $-algebras} 
\label{secproalg}

Let $ G $ be a locally compact group and let $ A $ be a pro-$ G $-algebra. 
The space 
$ \Omega^n(A) $ of noncommutative $ n $-forms over $ A $ is defined by 
$ \Omega^n(A) = A^+ \hat{\otimes} A^{\hat{\otimes} n} $ for $ n \geq 0 $. 
We recall that $ A^+ $ denotes the unitarization of $ A $. 
From its definition as a tensor product it is clear that $ \Omega^n(A) $ becomes a 
pro-$ G $-module in a natural way. 
The differential $ d: \Omega^n(A)\rightarrow \Omega^{n + 1}(A) $ and the 
multiplication of forms $ \Omega^n(A) \hat{\otimes} \Omega^m(A) \rightarrow 
\Omega^{n + m}(A) $ are defined as usual ~\cite{CQ4} and it is clear that both are equivariant linear 
maps. Multiplication of forms yields in particular an 
$ A $-bimodule structure on $ \Omega^n(A) $ for all $ n $. 
Apart from the ordinary product of differential forms we have the Fedosov product  given by
\begin{equation*}
\omega \circ \eta = \omega \eta - (-1)^{|\omega|} d\omega d\eta
\end{equation*}
for homogenous forms $ \omega $ and $ \eta $. Consider 
the pro-$ G $-module $ \Omega^{\leq n}(A) = 
A \oplus \Omega^1(A) \oplus \cdots \oplus \Omega^n(A) $ equipped with the Fedosov product where forms above degree $ n $ are ignored. It is 
easy to check that this multiplication is associative and turns 
$ \Omega^{\leq n}(A) $ into a pro-$ G $-algebra. 
Moreover we have the usual $ \mathbb{Z}_2 $-grading on $ \Omega^{\leq n}(A) $ into 
even and odd forms. 
The natural projection $ \Omega^{\leq m}(A)
\rightarrow \Omega^{\leq n}(A) $ for $ m \geq n $ is an equivariant homomorphism
and compatible with the grading. 
Hence we get a projective system $ (\Omega^{\leq n}(A))_{n \in \mathbb{N}} $ 
of pro-$ G $-algebras. By definition the periodic differential envelope 
$ \theta\Omega(A) $ of $ A $ is the pro-$ G $-algebra obtained as the 
projective limit of this system. We define the periodic tensor algebra 
$ \mathcal{T}A $ of $ A $ to be the even part of 
$ \theta \Omega(A) $. 
If we set $ \mathcal{T}A/(\mathcal{J}A)^n := 
A \oplus \Omega^2(A) \oplus \cdots \oplus \Omega^{2n - 2}(A) $ 
we can describe $ \mathcal{T}A $ 
as the projective limit of the projective system 
$ (\mathcal{T}A/(\mathcal{J}A)^n)_{n \in \mathbb{N}} $. 
The natural projection $ \theta \Omega(A) \rightarrow A $ 
restricts to an equivariant homomorphism $ \tau_A: \mathcal{T}A \rightarrow A $. 
Since the natural inclusions 
$ A \rightarrow A \oplus \Omega^2(A) \oplus \cdots \oplus \Omega^{2n - 2}(A)$ assemble 
to give an equivariant linear section $ \sigma_A $ for $ \tau_A $ we obtain an extension 
\begin{equation*} 
 \xymatrix{
     \mathcal{J}A \;\; \ar@{>->}[r] &
         \mathcal{T}A \ar@{->>}[r]^{\tau_A} &
           A   
     }
\end{equation*}
of pro-$ G $-algebras where $ \mathcal{J}A $ is by definition the 
projective limit of the pro-$ G $-algebras $ \mathcal{J}A/(\mathcal{J}A)^n := 
\Omega^2(A) \oplus \cdots \oplus \Omega^{2n - 2}(A) $. \\
This section is devoted to the study of the pro-$ G $-algebras 
$ \mathcal{T}A $ and $ \mathcal{J}A $. Since this part of the equivariant theory is a straightforward 
extension of ordinary Cuntz-Quillen theory we have omitted some of the proofs. For 
more details we refer to ~\cite{Meyer}. \\
Let $ m^n: N^{\otimes n} \rightarrow  N $ be the iterated multiplication in an arbitrary pro-$ G $-algebra $ N $. Then 
$ N $ is called $ k $-nilpotent for $ k \in \mathbb{N} $ if the iterated multiplication
$ m^k: N^{\hat{\otimes} k} \rightarrow N $ is zero. 
It is called nilpotent if $ N $ is $ k $-nilpotent for some $ k \in \mathbb{N} $. 
We call $ N $ locally nilpotent if for every  
equivariant linear map $ f: N \rightarrow C $ with constant range $ C $ there 
exists $ n \in \mathbb{N} $ such that $ f m^n = 0 $. In particular nilpotent pro-$ G $-algebras are 
locally nilpotent. 
An extension $ 0 \rightarrow K \rightarrow E \rightarrow Q \rightarrow 0 $ of pro-$ G $-algebras is 
called locally nilpotent ($ k $-nilpotent, nilpotent) 
if $ K $ is locally nilpotent ($ k $-nilpotent, nilpotent).
\begin{lemma} \label{JALocNil} The pro-$ G $-algebra $ \mathcal{J}A $ is locally nilpotent.
\end{lemma}
\proof Let $ l: \mathcal{J}A \rightarrow C $ be an equivariant linear 
map. By the construction of projective limits it follows that there 
exists $ n \in \mathbb{N} $ such that $ l $ factors through 
$ \mathcal{J}A/(\mathcal{J}A)^n $. The pro-$ G $-algebra 
$ \mathcal{J}A/(\mathcal{J}A)^n $ is $ n $-nilpotent by the definition of the 
Fedosov product. Hence $ l m^n_{\mathcal{J}A} = 0 $ as desired. \qed  
\begin{lemma} \label{tensorlemma} Let $ N $ be a locally nilpotent pro-$ G $-algebra and let $ A $ be 
any pro-$ G $-algebra. Then the pro-$ G $-algebra $ A \hat{\otimes} N $ is locally nilpotent. 
\end{lemma} 
\proof Let $ f: A \hat{\otimes} N \rightarrow C $ be an equivariant 
linear map with constant range. By the construction of tensor 
products in $ \pro(G\LMod) $ this map can be written as $ g(f_1 \hat{\otimes} f_2) $ 
for equivariant linear maps $ f_1: A \rightarrow C_2, f_2: N \rightarrow C_2 $ 
with constant range and an equivariant bounded linear map $ g: C_1 \hat{\otimes} C_2 
\rightarrow C $. Since $ N $ is locally nilpotent there exists a natural number 
$ n $ such that $ f_2 m_N^n = 0 $. Up to a coordinate flip the $ n $-fold
multiplication in $ A \hat{\otimes} N $ is given by $ m_A^n \hat{\otimes} m_N^n $. 
This implies $ f m_{A \hat{\otimes} N}^n = 0 $ for the multiplication 
$ m_{A \hat{\otimes} N} $ in $ A \hat{\otimes} N $. Hence $ A \hat{\otimes} N $ 
is locally nilpotent. \qed \\
Next we want to study the pro-$ G $-algebra $ \mathcal{T}A $. 
In order to formulate its universal property we need another definition. 
An equivariant linear map $ l: A \rightarrow B $ 
between pro-$ G $-algebras is called a lonilcur if its curvature 
$ \omega_l: A \hat{\otimes} A \rightarrow B $ defined by 
$ \omega_l(a,b) = l(ab) - l(a)l(b) $ is 
locally nilpotent, that is, if for every equivariant linear map 
$ f: B \rightarrow C $ with 
constant range $ C $ there exists $ n \in \mathbb{N} $ such that 
$ f m^n_B \omega^{\hat{\otimes} n}_l = 0 $. The term lonilcur is an abbreviation for 
"equivariant linear map with \emph{lo}cally \emph{nil}potent \emph{cur}vature". 
It is clear that every equivariant homomorphism 
is a lonilcur because the curvature is zero in this case. 
Using the fact that $ \mathcal{J}A $ is locally nilpotent one checks 
easily that the natural map $ \sigma_A: A \rightarrow \mathcal{T}A $ 
is a lonilcur. 
\begin{prop}\label{PeriodicTensorAlg} Let $ A $ be a pro-$ G $-algebra. The 
pro-$ G $-algebra $ \mathcal{T}A $ and the equivariant linear map $ \sigma_A: A 
\rightarrow \mathcal{T}A $ satisfy the following universal property. 
If $ l: A \rightarrow B $ is a lonilcur into a pro-$ G $-algebra $ B $ 
there exists a unique equivariant homomorphism $ [[l]]: \mathcal{T}A 
\rightarrow B $ such that $ [[l]] \sigma_A = l $. 
\end{prop}
Let us now define and study quasifree pro-$ G $-algebras. 
\begin{definition} A pro-$ G $-algebra $ R $ is called $ G $-equivariantly quasifree 
if there exists an equivariant splitting homomorphism $ R \rightarrow \mathcal{T}R $ 
for the natural projection $ \tau_R $.
\end{definition}
By abuse of language we will occasionally speak of quasifree pro-$ G $-algebras instead of $ G $-equivariantly 
quasifree $ G $-algebras although the latter is the correct terminology for 
a pro-$ G $-algebra which is quasifree as a pro-algebra. \\
In the following theorem the class of quasifree pro-$ G $-algebras is characterized. 
\begin{theorem} \label{qf}
Let $ G $ be a locally compact group and let $ R $ be a pro-$ G $-algebra. Then the
following conditions are equivalent:
\begin{bnum}
\item[a)] $ R $ is $ G $-equivariantly quasifree.
\item[b)] There exists a family of equivariant homomorphisms
$ v_n: R \rightarrow \mathcal{T}R/(\mathcal{J}R)^n $ such that $ v_1 = \id $ and
$ v_{n + 1} $ is a lifting of $ v_n $.
\item[c)] For every locally nilpotent extension
$ 0 \rightarrow K \rightarrow E \rightarrow Q \rightarrow 0 $ 
of pro-$ G $-algebras and every equivariant homomorphism $ f: R \rightarrow Q $ 
there exists an equivariant lifting homomorphism $ h: R \rightarrow E $.
\item[d)] For every nilpotent extension 
$ 0 \rightarrow K \rightarrow E \rightarrow Q \rightarrow 0 $ 
of pro-$ G $-algebras and every equivariant homomorphism $ f: R \rightarrow Q $ 
there exists an equivariant lifting homomorphism $ h: R \rightarrow E $.
\item[e)] For every $ 2 $-nilpotent
extension $ 0 \rightarrow K \rightarrow E \rightarrow Q \rightarrow 0 $ 
of pro-$ G $-algebras and every equivariant homomorphism
$ f: R \rightarrow Q $ there exists an equivariant lifting homomorphism
$ h: R \rightarrow E $.
\item[f)] For every $ 2 $-nilpotent
extension $ 0 \rightarrow K \rightarrow E \rightarrow R \rightarrow 0 $
of pro-$ G $-algebras there exists an equivariant
splitting homomorphism $ R \rightarrow E $.
\item[g)] There exists an equivariant splitting homomorphism
for the natural homomorphism $ \mathcal{T}R/(\mathcal{J}R)^2 \rightarrow R $.
\item[h)] There exists an equivariant linear map
$ \phi: R \rightarrow \Omega^2(R) $ satisfying
\begin{equation*} 
\phi(xy) = \phi(x)y + x \phi(y) - dxdy
\end{equation*} 
for all $ x,y \in R $.
\item[i)] There exists an equivariant linear map
$ \nabla: \Omega^1(R) \rightarrow \Omega^2(R) $ satisfying
\begin{equation*}
\nabla(x \omega) = x \nabla(\omega), \qquad
\nabla(\omega x) = \nabla(\omega) x - \omega dx
\end{equation*}
for all $ x \in R $ and $ \omega \in \Omega^1(R) $.
\item[j)] The $ R $-bimodule $ \Omega^1(R) $ is projective in
$ \pro(G \LMod) $.
\item[k)] There exists a projective resolution $ 0 \rightarrow P_1 \rightarrow P_0 \rightarrow R^+ $ 
of the $ R $-bimodule $ R^+ $ of length $ 1 $ in $ \pro(G \LMod) $.
\end{bnum}
\end{theorem}
Let us also include the following definitions.  
\begin{definition}
A pro-$ G $-algebra $ A $ is called $ n $-dimensional (with respect to $ G $) if there exists a 
projective resolution $ 0 \rightarrow P_n \rightarrow \cdots \rightarrow P_0 \rightarrow A^+ $ of 
the $ A $-bimodule $ A^+ $ of length $ n $ in $ \pro(G \LMod) $.
\end{definition}
\begin{definition}\label{defgradconn}
Let $ A $ be a pro-$ G $-algebra and let $ n > 0 $. An equivariant graded (right) 
connection on $ \Omega^n(A) $ is an equivariant linear map $ \nabla: \Omega^n(A) \rightarrow \Omega^{n + 1}(A) $ 
such that 
$$
\nabla(x \omega) = x \nabla(\omega), \qquad \nabla(\omega x) = \nabla(\omega) x + (-1)^n \omega dx
$$
for $ x \in A $ and $ \omega \in \Omega^n(A) $. 
\end{definition}
According to theorem \ref{qf} a pro-$ G $-algebra $ A $ is $ G $-equivariantly quasifree iff 
it is $ 1 $-dimensional with respect to $ G $. As in the non-equivariant case one has the following 
characterization of $ n $-dimensional algebras.
\begin{prop}\label{propndim}
Let $ G $ be a locally compact group and let $ A $ be a pro-$ G $-algebra. Then the
following conditions are equivalent:
\begin{bnum}
\item[a)] $ A $ is $ n $-dimensional with respect to $ G $.
\item[b)] The $ A $-bimodule $ \Omega^n(A) $ is projective in $ \pro(G \LMod) $.
\item[c)] There exists an equivariant graded connection on $ \Omega^n(A) $. 
\end{bnum}
\end{prop}
A basic example of a quasifree pro-$ G $-algebra is the algebra of complex numbers 
$ \mathbb{C} $ with the trivial $ G $-action. More generally we observe the following. 
\begin{lemma} Let $ A $ be a pro-algebra equipped with the trivial
$ G $-action. If $ A $ is quasifree as a pro-algebra it is
$ G $-equivariantly quasifree.
\end{lemma}
The following result is important. 
\begin{prop}\label{TAQuasifree} Let $ A $ be any pro-$ G $-algebra. The periodic tensor algebra 
$ \mathcal{T}A $ is $ G $-equivariantly quasifree. 
\end{prop}
\proof We have to show that there exists an equivariant splitting homomorphism 
for the projection $ \tau_{\mathcal{T}A}: \mathcal{T}\mathcal{T} A \rightarrow 
\mathcal{T}A $. Let us consider the equivariant linear map 
$ \sigma^2_A = \sigma_{\mathcal{T}A} \sigma_A: A \rightarrow 
\mathcal{T}\mathcal{T} A $. We want to show that $ \sigma_A^2 $ is a lonilcur. 
First we compute the curvature $ \omega_{\sigma_A^2} $ of $ \sigma_A^2 $ as 
follows:
\begin{align*}
\omega_{\sigma_A^2}&(x,y) = \sigma^2_A(xy) - \sigma^2_A(x) \circ \sigma^2_A(y) \\
&= \sigma_{\mathcal{T}A}(\sigma_A(xy)) - 
\sigma_{\mathcal{T}A}(\sigma_A(x) \circ \sigma_A(y)) + d\sigma^2_A(x) d\sigma^2_A(y) \\
&= \sigma_{\mathcal{T}A}(\omega_{\sigma_A}(x,y)) + d\sigma^2_A(x) d\sigma^2_A(y).
\end{align*}
Consider the equivariant linear map 
$ \sigma_A = \tau_{\mathcal{T}A} \sigma^2_A $. 
Since $ \tau_{\mathcal{T}A} $ is a homomorphism we obtain 
$ \omega_{\sigma_A} = \tau_{\mathcal{T}A} \omega_{\sigma^2_A} $.
Let $ l: \mathcal{T}\mathcal{T}A \rightarrow C $ be an equivariant linear
map with constant range $ C $. Composition with 
$ \sigma_{\mathcal{T}A}: \mathcal{T}A \rightarrow \mathcal{T}\mathcal{T}A $ yields a map 
$ k = l \sigma_{\mathcal{T}A}: \mathcal{T}A \rightarrow C $ with constant range. 
Since $ \sigma_A $ is a lonilcur there exists $ n \in \mathbb{N} $ such that 
\begin{equation*}
k m^n_{\mathcal{T}A} \omega_{\sigma_A}^{\hat{\otimes} n} = 
k m^n_{\mathcal{T}A} \tau_{\mathcal{T}A}^{\hat{\otimes} n} \omega_{\sigma^2_A}^{\hat{\otimes} n} = 
k \tau_{\mathcal{T}A} m^n_{\mathcal{T}\mathcal{T}A} 
\omega_{\sigma^2_A}^{\hat{\otimes} n} = 0.
\end{equation*}
By the construction of $ \mathcal{T}\mathcal{T}A $ the map $ l $ 
factors over $ \mathcal{T}\mathcal{T}A/(\mathcal{J}(\mathcal{T}A))^m $ 
for some $ m $. Using the formula for the curvature of $ \sigma_A^2 $ 
and our previous computation we obtain $ l m_{\mathcal{T}\mathcal{T}A}^{mn} 
\omega_{\sigma^2_A}^{\hat{\otimes} mn} = 0 $. Hence $ \sigma_A^2 $ is 
a lonilcur. By the universal property of $ \mathcal{T}A $ there exists 
a homomorphism $ v = [[\sigma_A^2]]: \mathcal{T}A \rightarrow 
\mathcal{T}\mathcal{T}A $ such that $ v \sigma_A = \sigma_A^2 $. 
This implies $ (\tau_{\mathcal{T}A} v)\sigma_A = \tau_{\mathcal{T}A} 
\sigma_{\mathcal{T}A} \sigma_A = \sigma_A $.
From the uniqueness assertion of proposition \ref{PeriodicTensorAlg} 
we deduce $ \tau_{\mathcal{T}A} v = \id $. This means that 
$ \mathcal{T}A $ is quasifree. \qed \\
In connection with unital algebras the following result is useful. 
\begin{prop}\label{Unitalqf} Let $ A $ be a pro-$ G $-algebra. Then $ A $ is $ G $-equivariantly quasifree if and only if $ A^+ $ is 
$ G $-equivariantly quasifree.
\end{prop}
We will now define universal locally nilpotent extensions of pro-$ G $-algebras. 
\begin{definition}
Let $ A $ be a pro-$ G $-algebra. A universal locally nilpotent extension of $ A $ 
is an extension of pro-$ G $-algebras 
$ 0 \rightarrow N \rightarrow R \rightarrow A \rightarrow 0 $ where $ N $ is locally nilpotent and $ R $ 
is $ G $-equivariantly quasifree.
\end{definition}
We equip the Fr\'echet algebra $ C^\infty[0,1] $ of smooth functions on 
the interval $ [0,1] $ with the bounded bornology and view it as a $ G $-algebra with the 
trivial $ G $-action. 
An equivariant homotopy is an equivariant homomorphism of 
pro-$ G $-algebras 
$ h: A \rightarrow B \hat{\otimes} C^\infty[0,1] $ where $ C^\infty[0,1] $ is 
viewed as a constant pro-$ G $-algebra. 
For each $ t \in [0,1] $ evalutation at 
$ t $ defines an equivariant homomorphism $ h_t: A \rightarrow B $. Two equivariant homomorphisms are equivariantly 
homotopic if they can be connected by an equivariant homotopy. 
We will also write $ B[0,1] $ for the pro-$ G $-algebra $ B \hat{\otimes} C^\infty[0,1] $. 
\begin{prop}\label{UnivExt}
Let $ (\iota,\pi): 0 \rightarrow N \rightarrow R \rightarrow A
\rightarrow 0 $ be a universal locally nilpotent extension of $ A $.
If $ (i,p): 0 \rightarrow K \rightarrow E \rightarrow Q \rightarrow 0 $ is any
other locally nilpotent extension and $ \phi: A \rightarrow Q $ 
an equivariant homomorphism there exists a commutative diagram 
of pro-$ G $-algebras
\begin{equation*} 
 \xymatrix{
     N \;\; \ar@{>->}[r]^{\iota} \ar@{->}[d]^\xi &
         R \ar@{->>}[r]^{\pi} \ar@{->}[d]^\psi&
           A \ar@{->}[d]^\phi \\
 K \;\; \ar@{>->}[r]_i &
         E \ar@{->>}[r]_{p} &
           Q 
     }
\end{equation*}
Moreover the equivariant homomorphisms $ \xi $ and $ \psi $ are unique up to smooth homotopy. \\
More generally let $ (\xi_t,\psi_t,\phi_t) $ for $ t = 0,1 $ be equivariant
homomorphisms of extensions and let $ \Phi: A \rightarrow Q[0,1] $ be 
an equivariant homotopy connecting $ \phi_0 $ and $ \phi_1 $. Then
$ \Phi $ can be lifted to an equivariant homotopy $ (\Xi,\Psi,\Phi) $ 
between $ (\xi_0,\psi_0,\phi_0) $ and $ (\xi_1,\psi_1,\phi_1) $. 
\end{prop}
\proof Let $ v: R \rightarrow \mathcal{T}R $ be a splitting homomorphism 
for the projection $ \tau_R: \mathcal{T}R \rightarrow R $ and let 
$ s: Q \rightarrow E $ be an equivariant linear section for 
the projection $ p: E \rightarrow Q $. Since $ p(s \phi \pi) = \phi \pi $ is an 
equivariant homomorphism the curvature of $ s \phi \pi: R \rightarrow E $ 
has values in $ K $. Since by assumption $ K $ is locally nilpotent 
it follows that $ s \phi \pi $ is a lonilcur. From the universal property 
of $ \mathcal{T}R $ we obtain an equivariant homomorphism 
$ k = [[s \phi \pi]]: \mathcal{T}R \rightarrow E $ such that 
$ k \sigma_R = s \phi \pi $. Define $ \psi = k v: R \rightarrow E $. We have 
\begin{equation*}
(p k) \sigma_R = p s \phi \pi = \phi \pi = (\phi \pi \tau_R) \sigma_R
\end{equation*}
and by the uniqueness assertion in proposition \ref{PeriodicTensorAlg} we 
get $ p k = \phi \pi \tau_R $. Hence $ p \psi = p k v = \phi \pi \tau_R v = \phi \pi $ 
as desired. Moreover $ \psi $ maps $ N $ into $ K $ and restricts consequently 
to an equivariant homomorphism $ \xi: N \rightarrow K $ making the diagram commutative. \\
The assertion that $ \psi $ and $ \xi $ are uniquely defined up to smooth homotopy
follows from the more general statement about the lifting of homotopies. 
Hence let $ (\xi_t,\psi_t,\phi_t) $ for $ t = 0,1 $ and 
$ \Phi: A \rightarrow Q[0,1] $ be given as above. 
Tensoring with $ C^\infty[0,1] $ yields an extension 
$ (i[0,1], p[0,1]): 0 \rightarrow K[0,1] \rightarrow 
E[0,1] \rightarrow Q[0,1] \rightarrow 0 $ of pro-$ G $-algebras. An equivariant linear splitting $ s[0,1] $ for 
this extension is obtained by tensoring $ s $ with the identity on 
$ C^\infty[0,1] $. Since $ \Phi_t \pi = p \psi_t $ for $ t = 0,1 $ the equivariant linear map $ l: R \rightarrow E[0,1] $ defined by 
\begin{equation*}
l = s[0,1] \Phi \pi + (\psi_0 - s \phi_0 \pi) \otimes (1 - t) 
+ (\psi_1 - s \phi_1 \pi) \otimes t
\end{equation*}
satisfies $ \ev_t l = \psi_t $ for $ t = 0,1 $ and $ p[0,1] l = \Phi \pi $. 
The map $ p[0,1] l = \Phi \pi $ is a homomorphism and hence the curvature 
of $ l $ has values in $ K[0,1] $. Due to lemma \ref{tensorlemma} 
the pro-$ G $-algebra $ K[0,1] = K \hat{\otimes} C^\infty[0,1] $ is locally nilpotent. Consequently we get an equivariant homomorphism 
$ [[l]]: \mathcal{T}R \rightarrow E[0,1] $ such that $ [[l]] \sigma_R = l $. 
We define $ \Psi = [[l]] v $ and in the same way as above we obtain 
$ p[0,1] \Psi = \Phi \pi $. An easy computation shows 
$ \Psi_t = \ev_t \Psi = \psi_t $ for $ t = 0,1 $. 
Clearly $ \Psi $ restricts to an equivariant homomorphism $ \Xi: N 
\rightarrow K[0,1] $ such that $ (\Xi,\Psi,\Phi) $ becomes an equivariant  homomorphism of extensions. \qed 
\begin{prop}\label{univext2}
Let $ A $ be a pro-$ G $-algebra. The extension
$ 0 \rightarrow \mathcal{J}A \rightarrow\mathcal{T}A \rightarrow A
\rightarrow 0 $ is a universal locally nilpotent extension of $ A $.
If $ 0 \rightarrow N \rightarrow R \rightarrow A \rightarrow 0 $ is any
other universal locally nilpotent extension of $ A $ it is equivariantly
homotopy equivalent over $ A  $ to $ 0 \rightarrow \mathcal{J}A
\rightarrow\mathcal{T}A \rightarrow A \rightarrow 0 $. In particular 
$ R $ is equivariantly homotopy equivalent to $ \mathcal{T}A $ 
and $ N $ is equivariantly homotopy equivalent to $ \mathcal{J}A $.
\end{prop}
\proof The pro-$ G $-algebra $ \mathcal{J}A $ is locally nilpotent by 
lemma \ref{JALocNil}. Moreover $ \mathcal{T}A $ is quasifree 
by proposition \ref{TAQuasifree}. Hence the assertion follows from proposition \ref{UnivExt}. \qed 

\section{Equivariant differential forms}
\label{secdiffform}

In the previous section we have seen that the space of noncommutative $ n $-forms 
$ \Omega^n(A) $ for a pro-$ G $-algebra $ A $ is a pro-$ G $-module in a 
natural way. 
Let now $ A $ be any pro-$ G $-algebra and consider the covariant pro- module $ \Omega^n_G(A) = 
\mathcal{O}_G \hat{\otimes} \Omega^n(A) $. The $ G $-action on this space is 
defined by
\begin{equation*}
t \cdot (f(s) \otimes \omega) = f(t^{-1} s t) \otimes t \cdot \omega
\end{equation*}
for all $ f \in \mathcal{O}_G $ and $ \omega \in \Omega^n(A) $ and  
the $ \mathcal{O}_G $-module structure is given by multiplication. 
\begin{definition}
Let $ A $ be a pro-$ G $-algebra. The covariant pro-module $ \Omega^n_G(A) $ is 
called the space of equivariant $ n $-forms over $ A $. 
\end{definition}
Let us define operators $ d $ and 
$ b_G $ on equivariant differential forms by
\begin{equation*}
d(f(s) \otimes \omega) = f(s) \otimes d\omega
\end{equation*}
and
\begin{equation*}
b_G(f(s) \otimes \omega dx) = (-1)^n(f(s) \otimes (\omega
x - (s^{-1} \cdot x) \omega))
\end{equation*}
for $ \omega \in \Omega^n(A) $ and $ x \in A $. We remark that the definition of the 
operator $ b_G $ goes back at least to the work of Brylinski ~\cite{Brylinski1}.
Moreover in order to clarify our notation we point out that one may view elements in 
$ \Omega^n_G(A) $ as functions from $ G $ to $ \Omega^n(A) $. 
In particular the precise meaning of the last formula is that evaluation of 
$ b_G(f \otimes \omega dx) \in \Omega^n_G(A) $ at the group element $ s \in G $ 
yields $ (-1)^n(f(s) (\omega x - (s^{-1} \cdot x) \omega)) \in \Omega^n(A) $. \\
Having this in mind we want to study the properties of the operators 
$ d $ and $ b_G $. As in the non-equivariant case we clearly have $ d^2 = 0 $. 
The  operator $ b_G $ should be thought of as a twisted version of the ordinary
Hochschild boundary. We compute for $ \omega \in \Omega^n(A) $ and 
$ x,y \in A $ 
\begin{align*}
b_G^2(f(s) &\otimes \omega dx dy) = 
b_G((-1)^{n + 1}(f(s) \otimes \omega dx y - 
f(s) \otimes (s^{-1} \cdot y) \omega dx)) \\
&= b_G((-1)^{n + 1}(f(s) \otimes \omega d(xy) - f(s) \otimes \omega x dy 
- f(s) \otimes (s^{-1} \cdot y) \omega dx)) \\
&= (-1)^n(-1)^{n + 1}(f(s) \otimes \omega x y - f(s) \otimes s^{-1} \cdot (xy)\omega \\
&\qquad - (f(s) \otimes \omega x y - f(s) \otimes (s^{-1} \cdot y) \omega x) \\
&\qquad - (f(s) \otimes (s^{-1} \cdot y) \omega x - 
f(s) \otimes (s^{-1} \cdot x)(s^{-1} \cdot y)\omega)) = 0. 
\end{align*}
This shows $ b_G^2 = 0 $ and hence $ b_G $ is an ordinary differential. We 
will call $ b_G $ the equivariant Hochschild operator. \\
Similar to the non-equivariant case we construct an equivariant
Karoubi operator $ \kappa_G $ and an equivariant Connes operator 
$ B_G $ out of $ d $ and $ b_G $. We define
\begin{equation*}
\kappa_G = \id - (b_G d + d b_G)
\end{equation*}
and on $ \Omega^n_G(A) $ we set
\begin{equation*}
B_G = \sum_{j = 0}^n \kappa_G^j d.
\end{equation*}
Using that $ \kappa_G $ commutes with $ d $ and $ d^2 = 0 $ we obtain 
$ B_G^2 = 0 $. 
Let us record the following explicit formulas on $ \Omega^n_G(A) $.
For $ n \geq 1 $ we have
\begin{equation*}
\kappa_G(f(s) \otimes \omega dx) = (-1)^{n - 1} f(s) \otimes
(s^{-1} \cdot dx) \omega
\end{equation*}
and we obtain $ \kappa_G(f(s) \otimes x) = f(s) \otimes s^{-1} \cdot x $ 
for $ f(s) \otimes x \in \Omega^0_G(A) $. For the Connes operator we compute
\begin{equation*}
B_G(f(s) \otimes x_0dx_1 \cdots dx_n) = \sum_{i = 0}^n (-1)^{ni}
f(s) \otimes s^{-1} \cdot(dx_{n + 1 - i} \cdots dx_n)dx_0 \cdots
dx_{n - i}.
\end{equation*}
In addition we have the symmetry operator $ T $ which is defined on any covariant pro-module and takes 
the form
\begin{equation*}
T(f(s) \otimes \omega) = f(s) \otimes s^{-1} \cdot \omega 
\end{equation*}
on $ \Omega^n_G(A) $. It is easy to check that all operators constructed so far are 
covariant. \\
In order to keep the formulas readable we will frequently write $ b $ instead of 
$ b_G $ in the sequel and use similar simplifications for the other operators. \\
We need the following lemma concerning relations between the
operators constructed above. 
See ~\cite{CQ2} for the corresponding 
assertion in the non-equivariant context.  
\begin{lemma}\label{diffformformel}
On $ \Omega^n_G(A) $ the following relations hold:
\begin{bnum}
\item[a)] $ \kappa^{n + 1} d = T d $
\item[b)] $ \kappa^n = T + b \kappa^n d $
\item[c)] $ b \kappa^n = b T $
\item[d)] $ \kappa^{n + 1} = (\id - db) T $
\item[e)] $ (\kappa^{n + 1} - T)(\kappa^n - T) = 0 $
\item[f)] $ Bb + bB = \id - T $
\end{bnum}
\end{lemma}
\proof a) follows directly from the explicit formula for $ \kappa $ from above. b) Using again the formula for $ \kappa $ we compute 
\begin{align*}
\kappa^n(f(s) & \otimes x_0dx_1 \cdots dx_n) = 
f(s) \otimes s^{-1} \cdot(dx_1 \cdots dx_n)x_0 \\
&= f(s) \otimes s^{-1} \cdot (x_0dx_1 \cdots dx_n) + 
(-1)^n b(f(s) \otimes s^{-1} \cdot (dx_1 \cdots dx_n) dx_0) \\
&= T(f(s) \otimes x_0dx_1 \cdots dx_n) + 
b \kappa^n d(f(s) \otimes x_0dx_1 \cdots dx_n).
\end{align*}
c) follows by applying the Hochschild boundary $ b $ to both
sides of b). d) Apply $ \kappa $ to b) and use a).
e) is a consequence of b) and d). f) We compute
\begin{align*}
Bb +& bB = \sum_{j = 0}^{n - 1} \kappa^j db + \sum_{j = 0}^n b\kappa^j d
= \sum_{j = 0}^{n - 1} \kappa^j(db + bd) + \kappa^n bd \\
&= \id - \kappa^n(\id - bd) 
= \id - \kappa^n(\kappa + db) =
\id - T + dbT - Tdb = \id - T
\end{align*}
where we use d) and b) and the fact that $ T $ commutes with 
covariant maps due to proposition \ref{covparaadd}. \qed \\
Let us summarize this discussion as follows.
\begin{prop} Let $ A $ be a pro-$ G $-algebra. The space $ \Omega_G(A) $
of equivariant differential forms is a paramixed complex in the
category $ \pro(G\LSMod) $ of covariant pro-modules and all the operators constructed above are covariant.
\end{prop}
As for ordinary differential forms we define 
$ \Omega^{\leq n}_G(A) = \Omega^0_G(A) \oplus \Omega^1_G(A) \oplus \cdots \oplus \Omega^n_G(A) $ for all $ n \geq 0 $. We have the usual $ \mathbb{Z}_2 $-grading on 
$ \Omega^{\leq n}_G(A) $ into even and odd forms. 
The natural projection $ \Omega^{\leq m}_G(A)
\rightarrow \Omega^{\leq n}_G(A) $ for $ m \geq n $ is a covariant homomorphism
and compatible with the grading. 
Hence we obtain a projective system $ (\Omega^{\leq n}_G(A))_{n \in \mathbb{N}} $ 
and we let $ \theta\Omega_G(A) $ be the corresponding projective limit. \\
Using lemma \ref{relproj} we easily obtain the following fact.  
\begin{lemma} \label{localprj}
For any pro-$ G $-algebra $ B $ the covariant pro-module 
$ \theta\Omega_G(B \hat{\otimes} \mathcal{K}_G) $ is relatively projective. 
\end{lemma}

\section{The equivariant $ X $-complex}
\label{secX}

In this section we define and study the equivariant $ X $-complex. 
Apart from the periodic tensor algebra introduced in section 3.1 this object is the main ingredient in the definition of equivariant 
periodic cyclic homology. \\
Consider the paramixed complex $ \Omega_G(A) $ of equivariant differential forms 
for a pro-$ G $-algebra $ A $ which was defined in the previous section. 
Following Cuntz and Quillen ~\cite{CQ2} we define the $ n $-th level 
of the Hodge tower associated to $ \Omega_G(A) $ by 
\begin{equation*}
\theta^n \Omega_G(A) = \bigoplus_{j = 0}^{n - 1} \Omega^j_G(A) \oplus \Omega_G^n(A)/b(\Omega^{n + 1}_G(A)).
\end{equation*}
It is easy to see that the operators $ d $ and $ b $ descend to 
$ \theta^n \Omega_G(A) $. Consequently the same holds true for $ \kappa $ 
and $ B $. Using the natural grading into even and odd forms we see that 
$ \theta^n \Omega_G(A) $ together with the boundary operator $ B + b $ becomes a paracomplex. 
For $ m \geq n $ there exists a natural covariant chain map $ \theta^m \Omega_G(A) \rightarrow 
\theta^n\Omega_G(A) $. By definition the the Hodge tower $ \theta \Omega_G(A) $ of $ A $ is the projective limit 
of the projective system $ (\theta^n \Omega_G(A))_{n \in \mathbb{N}} $ obtained in this way. \\
We emphasize that $ \theta^n \Omega_G(A) $ for an arbitrary pro-$ G $-algebra $ A $ is a projective systems of not necessarily separated 
covariant modules. However, we will only have to work with these objects
in the case they are in fact projective systems of separated spaces. \\
We define the Hodge filtration on $ \theta^n\Omega_G(A) $ by 
\begin{equation*}
F^k \theta^n\Omega_G(A) = b(\Omega^{k + 1}_G(A)) \oplus 
\bigoplus_{j = k + 1}^{n - 1} \Omega^j_G(A) \oplus 
\Omega_G^n(A)/b(\Omega^{n + 1}_G(A)). 
\end{equation*}
Clearly $ F^k \theta^n \Omega_G(A) $ is closed under $ b $ and $ B $. The Hodge filtration 
on $ \theta^n\Omega_G(A) $ is a finite decreasing filtration such that 
$ F^{-1} \theta^n\Omega_G(A) =  \theta^n\Omega_G(A) $ and 
$ F^n \theta^n\Omega_G(A) =  0 $. Remark that these definitions can be extended to arbitrary 
paramixed complexes of covariant modules in a straightforward way. 
\begin{definition} Let $ A $ be a pro-$ G $-algebra. The equivariant $ X $-complex 
$ X_G(A) $ of $ A $ is the paracomplex $ \theta^1 \Omega_G(A) $. Explicitly, we have  
\begin{equation*}
    \xymatrix{
      {X_G(A) \colon \ }
        {\Omega^0_G(A)\;} \ar@<1ex>@{->}[r]^-{d} &
          {\;\Omega^1_G(A)/ b(\Omega^2_G(A)).} 
            \ar@<1ex>@{->}[l]^-{b} 
               } 
\end{equation*}
\end{definition}
Let us point out that, despite of our terminology, $ X_G(A) $ is usually only a paracomplex and not a complex. 
Moreover we remark that we will only be interested in the equivariant $ X $-complex $ X_G(A) $ in the case that $ A $ is 
quasifree. Recall from theorem \ref{qf} that the $ A $-bimodule $ \Omega^1(A) $ is a projective object
in $ \pro(G\LMod) $ if $ A $ is a quasifree pro-$ G $-algebra. 
It follows easily that $ \Omega^1_G(A)/b(\Omega^2_G(A)) $ is a projective system of separated spaces in this 
case. \\
The following lemma shows how the equivariant $ X $-complex behaves with respect to unitarizations. This will be 
useful later on.
\begin{lemma} \label{XA}
For every pro-$ G $-algebra $ A $ the natural homomorphisms $ A \rightarrow A^+ $ and $ \mathbb{C} \rightarrow A^+ $ induce an 
isomorphism of paracomplexes
$$
X_G(A) \oplus \mathcal{O}_G[0] \cong X_G(A^+).
$$
\end{lemma}
\proof We have an evident isomorphism $ q_0: X^0_G(A) \oplus \mathcal{O}_G \cong X^0_G(A^+) $ in degree zero given by the identification
$$
X^0_G(A) \oplus \mathcal{O}_G = \mathcal{O}_G \cotimes A \oplus \mathcal{O}_G = \mathcal{O}_G \cotimes A^+ = X^0_G(A^+).
$$
Let $ q_1: X^1_G(A) \rightarrow X^1_G(A^+) $ be the map induced by the inclusion homomorphism. In order to construct 
an inverse of $ q_1 $ consider the map $ p_1: \mathcal{O}_G \cotimes \Omega^1(A^+) \rightarrow \mathcal{O}_G \cotimes \Omega^1(A) $ given by 
$$
p_1(f \otimes (a_0,\alpha_0)d(a_1, \alpha_1)) = f \otimes a_0 da_1 + f \otimes \alpha_0 da_1, 
\qquad p_1(f \otimes d(a_1,\alpha_1)) = f \otimes da_1.
$$
It is straightforward to verify that $ p_1 $ descends to a covariant map $ X^1_G(A^+) \rightarrow X^1_G(A) $. 
Moreover one checks easily $ p_1 q_1 = \id $. To prove $ q_1 p_1 = \id $ observe first that in $ X^1_G(A^+) $ we have 
$$
f \otimes (0,1) d(0,1) = f \otimes (0,1) d((0,1)(0,1)) = 2 f \otimes (0,1) d(0,1) 
$$
and hence $ f \otimes (0,1) d(0,1) = 0 $. This implies 
$$
f \otimes (a_0,\alpha_0) d(0,1) = f \otimes (a_0,\alpha_0)d((0,1)(0,1)) = 2 f \otimes (a_0,\alpha_0)d(0,1) = 0.
$$
Now we compute 
\begin{align*}
q_1 p_1((f \otimes &(a_0,\alpha_0)d(a_1, \alpha_1)) = f \otimes (a_0,0) d(a_1,0) + f \otimes \alpha_0 d(a_1, 0) \\
&= f \otimes (a_0,0)d(a_1,0) + f \otimes (0,\alpha_0) d(a_1,0) \\
&= f \otimes (a_0,\alpha_0)d(a_1,0) = (a_0, \alpha_0)d(a_1,\alpha_1)
\end{align*}
and 
$$
q_1 p_1((f \otimes d(a_1, \alpha_1)) = f \otimes d(a_1,0) =  f \otimes d(a_1,\alpha_1).
$$
Finally one checks easily that the map $ q $ is compatible with the differentials. This finishes the proof. \qed \\
If we set $ A = 0 $ in lemma \ref{XA} we obtain a simple description of the equivariant $ X $-complex of the complex numbers. 
\begin{lemma}\label{XC}
The equivariant $ X $-complex $ X_G(\mathbb{C}) $ of the complex numbers 
$ \mathbb{C} $ can be identified with the trivial supercomplex $ \mathcal{O}_G[0] $.
\end{lemma}
We are interested in the equivariant $ X $-complex of the periodic tensor algebra 
$ \mathcal{T}A $ of a pro-$ G $-algebra $ A $. The first goal is to
relate the covariant pro-module $ X_G(\mathcal{T}A) $ to equivariant differential 
forms over $ A $. If we denote the even part of $ \theta \Omega_G(A) $ by 
$ \theta \Omega^{ev}_G(A) $ we obtain a covariant isomorphism 
\begin{equation*}
X^0_G(\mathcal{T}A) = \mathcal{O}_G \hat{\otimes} \mathcal{T}A \cong 
\theta\Omega^{ev}_G(A)
\end{equation*}
according to the definition of $ \mathcal{T}A $. \\
Before we consider $ X^1_G(\mathcal{T}A) $ we have to make a
convention. We use the letter $ D $ for the equivariant linear map
$ \mathcal{T}A \rightarrow \Omega^1(\mathcal{T}A) $ usually denoted by 
$ d $. This will help us not to confuse this map with the differential 
$ d $ in $ \mathcal{T}A = \theta \Omega^{ev}(A) $. \\
As in ~\cite{Meyer} we obtain the following assertion. 
\begin{prop} \label{Omegaiso} Let $ A $ be any pro-$ G $-algebra. The following maps 
are equivariant linear isomorphisms. 
\begin{alignat*}{2}
&\mu_1: (\mathcal{T}A)^+ \hat{\otimes} A \hat{\otimes} (\mathcal{T}A)^+ \rightarrow 
\Omega^1(\mathcal{T}A), & \qquad & \mu_1(x \otimes a \otimes y) = x D\sigma_A(a) y \\ 
& \mu_2: (\mathcal{T}A)^+ \hat{\otimes} A \rightarrow \mathcal{T}A, 
& \qquad &\mu_2(x \otimes a) = x \circ \sigma_A(a) \\
&\mu_3: A \hat{\otimes}(\mathcal{T}A)^+ \rightarrow \mathcal{T}A, 
& \qquad &\mu_3(a \otimes x) = \sigma_A(a) \circ x. 
\end{alignat*}
Hence $ \Omega^1(\mathcal{T}A) $ is a free $ \mathcal{T}A $-bimodule and 
$ \mathcal{T}A $ is free as a left and right $ \mathcal{T}A $-module. 
\end{prop} 
Using proposition \ref{Omegaiso} we see that the map $ \mu_1: (\mathcal{T}A)^+ \hat{\otimes} A \hat{\otimes} (\mathcal{T}A)^+ \rightarrow \Omega^1(\mathcal{T}A) $ 
induces a covariant isomorphism 
$ \mathcal{O}_G \hat{\otimes} (\mathcal{T}A)^+ \hat{\otimes} A \hat{\otimes} (\mathcal{T}A)^+ \cong \Omega^1_G(\mathcal{T}A) $. Identifying 
equivariant commutators under this isomorphism yields 
a covariant isomorphism
\begin{equation*}
\Omega^1_G(\mathcal{T}A)/b_G(\Omega^2_G(\mathcal{T}A)) \cong 
\mathcal{O}_G \hat{\otimes} (\mathcal{T}A)^+ 
\hat{\otimes} A.
\end{equation*}
Using again $ \mathcal{T}A = \theta\Omega^{ev}(A) $ we obtain a covariant isomorphism
\begin{equation*}
X^1_G(\mathcal{T}A) \cong \theta\Omega^{odd}_G(A)
\end{equation*}
where $ \theta\Omega^{odd}_G(A) $ is the odd part of $ \theta\Omega_G(A) $. \\
Having identified $ X_G(\mathcal{T}A) $ and $ \theta\Omega_G(A) $ as 
covariant pro-modules we want to compare the differentials on both sides. To this 
end let $ f(s) \otimes xda $ be an element of $ \theta\Omega_G^{odd}(A) $ where 
$ x \in \mathcal{T}A \cong \theta\Omega_G^{ev}(A) $ and $ a \in A $. The differential 
$ X^1_G(\mathcal{T}A) \rightarrow X^0_G(\mathcal{T}A) $ in the equivariant 
$ X $-complex corresponds to 
\begin{align*}
\partial_1(f(s)&\otimes xda) = f(s)\otimes (x\circ a - (s^{-1} \cdot a)\circ x) \\
&= f(s)\otimes (x a - (s^{-1} \cdot a) x - dx da + (s^{-1} \cdot da) dx) \\
&= b(f(s)\otimes xda) - (\id + \kappa) d(f(s)\otimes xda).
\end{align*}
To compute the other differential we map $ \Omega^1_G(\mathcal{T}A) $ to 
$ \mathcal{O}_G \cotimes (\mathcal{T}A)^+ \hat{\otimes} A \hat{\otimes} (\mathcal{T}A)^+ $ 
using the inverse of the isomorphism $ \mu_1 $ in proposition \ref{Omegaiso} 
and compose with the covariant map 
$ \mathcal{O}_G \cotimes (\mathcal{T}A)^+ \hat{\otimes} A \hat{\otimes} (\mathcal{T}A)^+ \rightarrow 
\theta\Omega^{odd}_G(A) $ sending 
$ f(s) \otimes x_0 \otimes a \otimes x_1 $ to 
$ f(s) \otimes (s^{-1} \cdot x_1) \circ x_0 da $. The derivation rule for $ D $ yields the explicit formula 
\begin{align*}
\partial_0(&f(s) \otimes x_0 dx_1 \cdots dx_{2n}) = 
f(s) \otimes D(x_0 dx_1 \cdots dx_{2n}) \\
&= f(s) \otimes s^{-1} \cdot (dx_1 \cdots dx_{2n}) Dx_0 \\
&\qquad+ \sum_{j = 1}^n f(s) \otimes s^{-1} \cdot (dx_{2j + 1} \cdots dx_{2n}) \circ x_0 dx_1 \cdots dx_{2j - 2} D(x_{2j - 1}x_{2j}) \\
&\qquad- \sum_{j = 1}^n f(s) \otimes s^{-1} \cdot (dx_{2j + 1} \cdots dx_{2n}) \circ x_0 dx_1 \cdots dx_{2j - 2} 
\circ x_{2j - 1} Dx_{2j} \allowdisplaybreaks[3]\\
&\qquad- \sum_{j = 1}^n f(s) \otimes s^{-1} \cdot (x_{2j} dx_{2j + 1} \cdots dx_{2n}) \circ x_0 dx_1 \cdots 
dx_{2j - 2} Dx_{2j - 1} \allowdisplaybreaks[2] \\
&= \sum_{j = 0}^{2n} f(s) \otimes s^{-1} \cdot (dx_j \cdots dx_{2n}) dx_0 dx_1 \cdots dx_{j - 1} \\
&\qquad- \sum_{j = 1}^n b(f(s) \otimes s^{-1} \cdot (dx_{2j + 1} \cdots dx_{2n}) x_0 dx_1 \cdots dx_{2j - 1} dx_{2j} \\
&= B(f(s) \otimes x_0 dx_1 \cdots dx_{2n}) - \sum_{j = 0}^{n - 1}
\kappa^{2j} b(f(s) \otimes x_0 dx_1 \cdots dx_{2n})
\end{align*} 
for the operator corresponding to the differential $ X^0_G(\mathcal{T}A) \rightarrow X^1_G(\mathcal{T}A) $. 
This can be summarized as follows. 
\begin{prop}\label{Xdiff}
Under the identification $ X_G(\mathcal{T}A) \cong \theta\Omega_G(A) $ as above
the differentials of the equivariant $ X $-complex correspond to
\begin{alignat*}{2}
\partial_1 &= b - (\id + \kappa)d \qquad &&\text{on}\;\; \theta\Omega_G^{odd}(A)\\
\partial_0 &= - \sum_{j = 0}^{n - 1}
\kappa^{2j} b + B \qquad && \text{on}\;\; \Omega_G^{2n}(A).
\end{alignat*}
\end{prop}
We would like to show that the paracomplexes $ X_G(\mathcal{T}A) $ and 
$ \theta\Omega_G(A) $ are covariantly homotopy equivalent. However, at this point 
we cannot proceed as in the nonequivariant case. \\
Let us recall the situation for the ordinary $ X $-complex. The proof of
the homotopy equivalence between $ X(\mathcal{T}A) $ and $ \theta\Omega(A) $ given 
by Cuntz and Quillen ~\cite{CQ2}, ~\cite{CQ4} is based on the spectral decomposition of the Karoubi operator $ \kappa $. 
This decomposition is obtained from the polynomial relation
\begin{equation*}
(\kappa^{n + 1} - \id)(\kappa^n - \id) = 0
\end{equation*}
which holds on $ \Omega^n(A) $. Remark that this formula is related to the
fact that the cyclic permutation operator is of finite order on
$ \Omega^n(A) $. \\
In the equivariant theory the situation is different. The
equivariant cyclic permutation operator is in general of infinite
order, due to lemma \ref{diffformformel} e) the relevant relation for $ \kappa $
is 
\begin{equation*}
(\kappa^{n + 1} - T)(\kappa^n - T) = 0
\end{equation*}
on $ \Omega^n_G(A) $. Hence the proof from ~\cite{CQ2} cannot be carried over directly. \\
However, some additional work will in fact yield the following theorem. 
\begin{theorem}\label{homotopyeq} For any pro-$ G $-algebra $ A $ the paracomplexes $ X_G(\mathcal{T}A) $ and 
$ \theta\Omega_G(A) $ are covariantly homotopy equivalent. 
\end{theorem}
Due to proposition \ref{Xdiff} it suffices to prove that the 
paracomplexes $ (\theta\Omega_G(A), \partial) $ and 
$ (\theta\Omega_G(A), B + b) $ are 
covariantly homotopy equivalent. 
We define $ c_{2n} = c_{2n + 1} = (-1)^n n! $ for all $ n $. Consider the 
isomorphism 
$ c: \theta\Omega_G(A) \rightarrow \theta\Omega_G(A) $ given by 
$ c(\omega) = c_n\, \omega $ for 
$ \omega \in \Omega^n_G(A) $ and let 
$ \delta = c^{-1} (B + b) c $ be the boundary corresponding to 
$ B + b $ under this isomorphism. It is easy to check that 
\begin{equation*}
\delta = B - n b \qquad \text{on}\;\; \Omega_G^{2n}(A)
\end{equation*}
and
\begin{equation*}
\delta = - \frac{1}{n + 1}B + b \qquad \text{on}\;\; \Omega_G^{2n + 1}(A).
\end{equation*}
Hence in order to prove theorem \ref{homotopyeq} it is enough to show that 
$ (\theta\Omega_G(A),\partial) $ and $ (\theta\Omega_G(A),\delta) $ are covariantly 
homotopy equivalent. \\
In ordinary Cuntz-Quillen theory one proceeds by considering certain operators associated to the spectral decomposition of the operator $ \kappa^2 $. 
These operators are polynomials in $ \kappa^2 $ and explicit formulas can be found in ~\cite{Meyer}. 
Since we do not have a spectral decomposition of $ \kappa^2 $ in the equivariant situation we will work directly with these polynomials. \\
We begin with the operator $ N_n $ which is given by 
\begin{equation*} 
N_n = N_n(\kappa^2) = \frac{1}{n} \sum_{j = 0}^{n - 1} \kappa^{2j}
\end{equation*} 
for $ n \geq 1 $ and by $ N_0 = \id $. \\
Due to lemma \ref{diffformformel} a) we have $ \kappa^{2n + 1}d = T d $ on $ \Omega_G^{2n}(A) $. 
Hence we get  
\begin{equation} \label{cqheq1}
(\id - \kappa^2) N_{2n + 1} B = \frac{1}{2n + 1} (\id - \kappa^{2(2n + 1)}) B = \frac{1}{2n + 1}(\id - T^2) B
\end{equation}
on $ \Omega_G^{2n}(A) $. Similarly we have 
\begin{equation} \label{cqheq2}
(\id - \kappa^2) N_{2n + 1} b = \frac{1}{2n + 1} (\id - \kappa^{2(2n + 1)}) b = \frac{1}{2n + 1}(\id - T^2) b
\end{equation}
on $ \Omega_G^{2n + 1}(A) $ since $ \kappa^{2n + 1}b = Tb $ on $ \Omega_G^{2n + 1}(A) $ by lemma \ref{diffformformel} c). 
Next we define the polynomials $ f_n $ and $ g_n $ by 
\begin{equation*} 
f_n = f_n(\kappa^2) = N_n(\kappa^2) N_{n + 1}(\kappa^2) (\id + (n - \frac{1}{2})(\id - \kappa^2))
\end{equation*}
and 
\begin{equation*} 
g_n = g_n(\kappa^2) = - \biggl(n - \frac{1}{2}\biggr) N_n N_{n + 1} + N_n \frac{N_{n + 1} - \id}{\kappa^2 - \id} + \frac{N_n - \id}{\kappa^2 - \id}
\end{equation*}
for all $ n \geq 0 $. In addition we set $ f_j = \id $ and $ g_j = 0 $ for all negative integers $ j $. It is easy to check that each $ g_n $ is in 
fact a polynomial in $ \kappa^2 $ and that we have  
\begin{equation} \label{cqeq2}
g_n (\id - \kappa^2) = \id - f_n
\end{equation}
for all $ n $. We define covariant maps $ F_j $ by  
\begin{equation*} 
F_{2n - 1} = F_{2n} = f_{2n - 2} f_{2n - 1} f_{2n} 
\end{equation*} 
for all $ n $ and let 
$ F: \theta\Omega_G(A) \rightarrow \theta\Omega_G(A) $ be the operator which is given on $ j $-forms by $ F_j $. \\
We have to investigate the compatibility of the operator $ F $ with the differentials $ \partial $ and $ \delta $. 
Let us first determine the failure of $ F $ to define a chain map from $ (\theta\Omega_G(A),\partial) $ to $ (\theta\Omega_G(A),\partial) $. 
Using equations (\ref{cqeq2}) and (\ref{cqheq1}) we get on $ \Omega^{2n}_G(A) $ 
\begin{align*}\label{cqeq6}
\partial_0 F &- F \partial_0= B F_{2n} - \sum_{j = 0}^{n - 1} \kappa^{2j} b F_{2n} - F_{2n + 1} B + F_{2n - 1} \sum_{j = 0}^{n - 1} \kappa^{2j} b \\
&= (F_{2n} - F_{2n + 1}) B \\
&= f_{2n}(f_{2n - 2} f_{2n - 1} - f_{2n + 1} f_{2n + 2}) B \\
&= -f_{2n}((\id - f_{2n - 2}) f_{2n - 1} + (\id - f_{2n - 1}) - (\id - f_{2n + 2}) f_{2n + 1} - (\id - f_{2n + 1})) B \\
&= -f_{2n}(g_{2n - 2} f_{2n - 1} + g_{2n - 1} - g_{2n + 2} f_{2n + 1} - g_{2n + 1})(\id - \kappa)^2 B \\
&= (\id - T)Q_{2n}
\end{align*}
where 
\begin{align*}
Q_{2n} &= -\frac{1}{2n + 1} N_{2n} (\id + (2n - \frac{1}{2})(\id - \kappa^2)) \times \\
&\qquad \times  (g_{2n - 2} f_{2n - 1} + g_{2n - 1} - g_{2n + 2} f_{2n + 1} - g_{2n + 1}) (\id + T) B.
\end{align*}
Similarly, using equation (\ref{cqheq2}) we have on $ \Omega^{2n + 1}_G(A) $ 
\begin{align*}
\partial_1 F - F \partial_1 &= b F_{2n + 1} - (\id + \kappa) d F_{2n + 1} - F_{2n} b + F_{2n + 2} (\id + \kappa) d \\
&= (F_{2n + 1} - F_{2n}) b \\
&= f_{2n} (g_{2n - 2} f_{2n - 1} + g_{2n - 1} - g_{2n + 2} f_{2n + 1} - g_{2n + 1})(\id - \kappa)^2 b \\
&= (\id - T) Q_{2n + 1} 
\end{align*}
where 
\begin{align*}
Q_{2n + 1} &=  \frac{1}{2n + 1} N_{2n} (\id + (2n - \frac{1}{2})(\id - \kappa^2)) \times \\
&\qquad \times (g_{2n - 2} f_{2n - 1} + g_{2n - 1} - g_{2n + 2} f_{2n + 1} - g_{2n + 1}) (\id + T) b.
\end{align*}
An analogous computation is needed to determine the deviation of $ F $ to define a chain 
map from $ (\theta\Omega_G(A),\delta) $ to $ (\theta\Omega_G(A),\delta) $. We get on $ \Omega^{2n}_G(A) $ 
\begin{align*}
\delta_0 F - F \delta_0 &= B F_{2n} - nb F_{2n} - F_{2n + 1} B + n F_{2n - 1} b \\
&= (F_{2n} - F_{2n + 1}) B = (\id - T) Q_{2n}
\end{align*}
and 
\begin{align*}
\delta_1 F - F \delta_1 &= b F_{2n + 1} - \frac{1}{n + 1} B F_{2n + 1} - F_{2n} b + \frac{1}{n + 1} F_{2n + 2} B \\
&= (F_{2n + 1} - F_{2n}) b = (\id - T) Q_{2n + 1}
\end{align*}
on $ \Omega^{2n + 1}_G(A) $. 
Let $ Q: \theta\Omega_G(A) \rightarrow \theta\Omega_G(A) $ be the operator which is given on $ n $-forms by 
$ Q_n $. Then the previous computation yields  
\begin{equation}\label{eqcq}
\partial F - F \partial = (\id - T) Q, \qquad \delta F - F \delta = (\id - T) Q.
\end{equation}
The operator $ Q $ satisfies the following identities. We have on $ \Omega^{2n}_G(A) $
\begin{align*}
\partial_1 Q = \delta_1 Q = b Q_{2n}, \qquad Q \partial_0 = Q \delta_0 = Q_{2n + 1} B
\end{align*}
and similarly on $ \Omega_G^{2n + 1}(A) $ 
\begin{align*}
\partial_0 Q = \delta_0 Q = B Q_{2n + 1}, \qquad Q \partial_1 = Q \delta_1 = Q_{2n} b. 
\end{align*}
Since $ b Q_{2n} + Q_{2n + 1} B = 0 $ and $ B Q_{2n + 1} + Q_{2n} b = 0 $ we deduce 
\begin{lemma} \label{cqlemma4} The operator $ Q $ satisfies the relations
\begin{equation*}
\partial Q = \delta Q, \qquad Q \partial = Q \delta. 
\end{equation*}
Moreover
\begin{align*}
\partial Q + Q \partial = 0, \qquad \delta Q + Q \delta = 0, 
\end{align*}
that is, $ Q $ is a chain map of odd degree for both boundary operators. 
\end{lemma}
Using lemma \ref{cqlemma4} we define the operator $ P: \theta \Omega_G(A) \rightarrow \theta \Omega_G(A) $ by 
\begin{equation}
P = F + \frac{1}{2} Q \partial = F - \frac{1}{2} \partial Q = F + \frac{1}{2} Q \delta = F - \frac{1}{2} \delta Q
\end{equation} 
and calculate using equation (\ref{eqcq}) 
\begin{align*}
\partial P - P \partial &= \partial(F - \frac{1}{2}\,\partial Q) - 
(F + \frac{1}{2} Q \partial)\partial 
= \partial F - F \partial -\frac{1}{2}\, \partial^2 Q - \frac{1}{2}\, Q \partial^2 \\
&= (\id - T)Q - (\id - T)Q = 0.  
\end{align*}
In the same way we get 
$$
\delta P - P \delta = 0
$$
which shows that $  P $ defines a chain map from $ (\theta \Omega_G(A), \partial) $ to itself 
and also a chain map from $ (\theta \Omega_G(A), \delta) $ to itself. \\
Next we shall prove that these chain maps are homotopic to the identity. 
First observe that   
$$
\id - F_{2n - 1} = \id - F_{2n} = (\id - f_{2n - 2}) + (\id - f_{2n - 1}) f_{2n - 2} + (\id - f_{2n}) f_{2n - 1} f_{2n - 2}. 
$$
Hence if we set  
\begin{equation}\label{Seq1}
S_{2n - 1} = S_{2n} = g_{2n - 2} + g_{2n - 1} f_{2n - 2} + g_{2n} f_{2n - 1} f_{2n - 2}
\end{equation}
and let $ S: \theta \Omega_G(A) \rightarrow \theta \Omega_G(A) $ be the operator given on $ n $-forms by $ S_n $ we get 
$$
\id - F = (\id - \kappa^2) S. 
$$
Observe that we also have 
$$
\id - F_{2n - 1} = \id - F_{2n} = (\id - f_{2n}) + (\id - f_{2n - 1}) f_{2n} + (\id - f_{2n - 2}) f_{2n - 1} f_{2n} 
$$
which implies 
\begin{equation}\label{Seq2}
S_{2n - 1} = S_{2n} = g_{2n} + g_{2n - 1} f_{2n} + g_{2n - 2} f_{2n - 1} f_{2n}.
\end{equation}
Combining equations (\ref{Seq1}) and (\ref{Seq2}) we get 
\begin{equation}\label{Seq3}
S_{2n} - S_{2n + 2} = f_{2n}(g_{2n - 1} - g_{2n + 1} + g_{2n - 2} f_{2n - 1} - g_{2n + 2} f_{2n + 1}).
\end{equation}
Let us consider the chain map $ P: (\theta \Omega_G(A), \partial) \rightarrow (\theta \Omega_G(A), \partial) $.  
We define 
\begin{equation*}
h_{2n} = (\id + \kappa) d - b, \qquad h_{2n + 1} = 0 
\end{equation*}
and calculate
\begin{equation*}
\partial h + h \partial = -(b - (\id + \kappa)d)^2 = (\id + \kappa)(bd + db) = (\id + \kappa)(\id - \kappa) 
= \id - \kappa^2. 
\end{equation*}
It follows that $ \id - \kappa^2 $ is homotopic to zero with respect to the boundary $ \partial $. 
Now we set 
$$ 
H_{2n} = h_{2n} S_{2n} +\frac{1}{2} Q_{2n}, \qquad H_{2n + 1} = 0
$$ 
and compute on $ \Omega^{2n}_G(A) $
\begin{align*}
\partial H + H \partial &= \partial h_{2n} S_{2n} + \frac{1}{2}\,\partial Q_{2n} = \id - F_{2n} + \frac{1}{2}\, \partial Q_{2n} = \id - P_{2n}. 
\end{align*}
Observe that by lemma \ref{diffformformel} a) we have on $ \Omega^{2n}_G(A) $ 
\begin{align} \label{Seq4}
N_{2n + 1} &(\id + \kappa) d = \frac{1}{2n + 1} \sum_{j = 0}^{2n} \kappa^{2j} (\id + \kappa) d \\
&= \frac{1}{2n + 1} \sum_{j = 0}^{2n} \kappa^j (\id + \kappa^{2n + 1}) d = \frac{1}{2n + 1} (\id + T) B. \notag
\end{align}
Hence using equation (\ref{Seq3}) and (\ref{Seq4}) we get on $ \Omega^{2n + 1}_G(A) $ 
\begin{align*}
h_{2n + 2}&(S_{2n} - S_{2n + 2}) b = N_{2n} (\id + (2n - \frac{1}{2})(\id - \kappa^2)) \times \\
&\qquad \times (g_{2n - 1} - g_{2n + 1} + g_{2n - 2} f_{2n - 1} - g_{2n + 2} f_{2n + 1}) N_{2n + 1} (\id + \kappa) db = - Q_{2n} b
\end{align*}
and compute on $ \Omega^{2n + 1}_G(A) $ 
\begin{align*}
\partial H + H \partial &= - h_{2n + 2} S_{2n + 2} (\id + \kappa) d + h_{2n} S_{2n} b + \frac{1}{2} Q_{2n} b \\
&= \id - F_{2n + 1} + h_{2n + 2}(S_{2n} - S_{2n + 2}) b + \frac{1}{2} Q_{2n} b \\
&= \id - F_{2n + 1} - Q_{2n} b  + \frac{1}{2} Q_{2n} b = \id - P_{2n + 1}. 
\end{align*}
We now consider the chain map $ P:(\theta \Omega_G(A), \delta) \rightarrow (\theta \Omega_G(A), \delta) $. Let us define 
\begin{equation*}
l_{2n} = (\id + \kappa)d, \qquad l_{2n + 1} = -\frac{1}{n + 1}(\id + \kappa)d
\end{equation*}
for all $ n $. Then the equation 
\begin{equation*}
[\delta, c^{-1}d c] = c^{-1}[B + b,d]c = c^{-1}(bd + db) c = c^{-1}(\id - \kappa) c 
= \id - \kappa
\end{equation*}
implies
\begin{equation*}
\delta l + l \delta = (\id + \kappa) (\id - \kappa) = \id - \kappa^2. 
\end{equation*}
It follows that $ \id - \kappa^2 $ is homotopic to zero with respect to the boundary $ \delta $. 
Now we set 
$$
L_{2n} = l_{2n} S_{2n} + \frac{1}{2} Q_{2n},\qquad L_{2n + 1} = l_{2n + 1} S_{2n + 1}
$$
and compute on $ \Omega^{2n}_G(A) $
\begin{align*}
\delta L + L \delta &= S_{2n} \delta l + S_{2n} l \delta + \frac{1}{2}\,\delta Q_{2n} = \id - P_{2n}. 
\end{align*}
On $ \Omega^{2n + 1}_G(A) $ we get 
\begin{align*}
\delta L + L \delta &= \delta l_{2n + 1} S_{2n + 1} + l_{2n} S_{2n} b + \frac{1}{2} Q_{2n} b \\
&= \id - F_{2n + 1} + l_{2n}(S_{2n} - S_{2n + 2}) b + \frac{1}{2} Q_{2n} b \\
&= \id - F_{2n + 1} + h_{2n}(S_{2n} - S_{2n + 2}) b + \frac{1}{2} Q_{2n} b \\
&= \id - F_{2n + 1} - \frac{1}{2} Q \delta = \id - P_{2n + 1}. 
\end{align*}
We summarize this discussion as follows. 
\begin{prop} \label{Pprop}
We have 
$$
\id - P = \partial H + H \partial, \qquad \id - P = \delta L + L \delta, 
$$
that is, the chain map $ \id - P $ is homotopic to zero with respect to both boundary operators. 
\end{prop}
Let us now determine the failure of $ F $ to define a chain map from $ (\theta\Omega_G(A),\delta) $ to $ (\theta\Omega_G(A),\partial) $. 
Using the relation $ \kappa^{2n} b = Tb $ on $ \Omega_G^{2n}(A) $ we compute 
$$
(\id - \kappa^2) N_{2n} b = \frac{1}{2n} (\id - \kappa^{2(2n)}) b = \frac{1}{2n}(\id - T^2) b
$$
on $ \Omega_G^{2n}(A) $ for $ n > 0 $. Hence we have on $ \Omega_G^{2n}(A) $ for $ n > 0 $
\begin{align*}\label{cqeq10}
\delta_0 F - F \partial_0 &= B F_{2n} - nb F_{2n} + F_{2n - 1} \sum_{j = 0}^{n - 1} \kappa^{2j} b - F_{2n + 1} B \\
&= (F_{2n} - F_{2n + 1})B -(\id - \kappa^2)  \sum_{j = 0}^{n - 2} (n - j - 1) \kappa^{2j} F_{2n} b \\
&= (\id - T) Q_{2n} - (\id - T) \sum_{j = 0}^{n - 2} (n - j - 1) \kappa^{2j} K_n b   
\end{align*}
where 
$$
K_n = \frac{1}{2n} f_{2n - 2} f_{2n - 1} N_{2n + 1} (\id +(2n - \frac{1}{2})(\id - \kappa^2)) (\id + T). 
$$
Similarly, on $ \Omega^{2n - 1}_G(A) $  we have $ \kappa^{2n} d = Td $ and hence 
$$
(\id - \kappa^2) N_{2n} d = \frac{1}{2n} (\id - \kappa^{2(2n)}) d = \frac{1}{2n}(\id - T^2) d.
$$
Using this we compute on $ \Omega_G^{2n - 1}(A) $  
\begin{align*}
\delta_1 F - F \partial_1 &= b F_{2n - 1} - \frac{1}{n} B  F_{2n - 1} - F_{2n - 2} b + F_{2n}(\id + \kappa)d \\
&= b(F_{2n - 1} - F_{2n - 2}) -\biggl(\frac{1}{n}\sum_{j = 0}^{2n - 1} \kappa^j F_{2n - 1} - (\id + \kappa) F_{2n}\biggl) d \\
&= (\id - T) Q_{2n - 1} + \frac{1}{n} (\id - \kappa^2) (\id + \kappa) \sum_{j = 0}^{n - 2} (n - j - 1) \kappa^{2j} F_{2n} d \\
&= (\id - T) Q_{2n - 1} + \frac{1}{n} (\id - T) (\id + \kappa) \sum_{j = 0}^{n - 2} (n - j - 1) \kappa^{2j} K_n d. 
\end{align*}
Hence if we set 
\begin{equation*}
R_{2n} = - \sum_{j = 0}^{n - 2} (n - j - 1) \kappa^{2j} K_n b, \qquad 
R_{2n - 1} = \frac{1}{n} (\id + \kappa) \sum_{j = 0}^{n - 2} (n - j - 1) \kappa^{2j} K_n d
\end{equation*}
for $ n >0 $ and $ R_0 = 0 $ we get  
\begin{equation*}
\delta F - F \partial = (\id - T)(Q + R), 
\end{equation*}
where, as before, $ R $ is given by $ R_n $ in degree $ n $. Similarly, we obtain on $ \Omega^{2n}_G(A) $ 
\begin{align*}
\partial_0 F - F \delta_0 &= BF_{2n} - \sum_{j = 0}^{n - 1} \kappa^{2j}b F_{2n} - F_{2n + 1} B + n F_{2n - 1} b \\
&= (F_{2n} - F_{2n + 1}) B + (\id - \kappa^2) \sum_{j = 0}^{n - 2} (n - j - 1)\kappa^{2j} F_{2n} b \\
&= (\id - T)Q_{2n} + (\id - T) \sum_{j = 0}^{n - 2} (n - j - 1) \kappa^{2j} K_n b   
\end{align*}
and on $ \Omega^{2n - 1}_G(A) $  
\begin{align*}\label{cqeq9}
\partial_1 F - F \delta_1 &= bF_{2n - 1} - (\id + \kappa) d F_{2n-1} - F_{2n - 2} b + \frac{1}{n} F_{2n} B \\
&= (F_{2n - 1} - F_{2n - 2})b - \frac{1}{n} (\id - \kappa^2) (\id + \kappa) \sum_{j = 0}^{n - 2} (n - j - 1) \kappa^{2j} F_{2n} d \\
&= (\id - T) Q_{2n - 1}  - \frac{1}{n} (\id - T) (\id + \kappa) \sum_{j = 0}^{n - 2} (n - j - 1) \kappa^{2j} K_n d. 
\end{align*}
Hence we have  
\begin{equation*}
\partial F - F \delta = (\id - T)(Q - R). 
\end{equation*} 
The operator $ R $ satisfies the identities 
$$
\delta R + R \partial = - \frac{1}{n} B R_{2n} - R_{2n - 1} \sum_{j = 0}^{n - 1} \kappa^{2j} b, 
\qquad \partial R + R \delta = -(\id + \kappa)d R_{2n} - R_{2n - 1} n b  
$$
on $ \Omega^{2n}_G(A) $ and 
$$
\delta R + R \partial = -nb R_{2n - 1} - R_{2n} (\id + \kappa)d, 
\qquad \partial R + R \delta = - \sum_{j = 0}^{n - 1} \kappa^{2j} b R_{2n - 1} - \frac{1}{n} R_{2n} B  
$$
on $ \Omega^{2n - 1}_G(A) $. Moreover we have on $ \Omega^{2n}_G(A) $ 
$$ 
FR - RF = F_{2n - 1} R_{2n} - R_{2n} F_{2n}= 0
$$
and similarly 
$$ 
FR - RF = F_{2n} R_{2n - 1} - R_{2n - 1} F_{2n - 1}= 0
$$
on $ \Omega^{2n - 1}_G(A) $. Finally one easily checks $ RQ = QR = 0 $. We summarize this as follows.  
\begin{lemma}\label{cqlemma5} 
We have the relations 
$$
\delta F - F \partial = (\id - T)(Q + R), \qquad \partial F - F \delta = (\id - T)(Q - R)
$$
as well as 
$$
\delta R + R \partial = 0, \qquad  \partial R + R \delta = 0
$$
and 
$$
[F,R] = FR - RF = 0, \qquad RQ = QR = 0. 
$$
\end{lemma}
Let us define $ \phi: (\theta \Omega_G(A),\partial) \rightarrow (\theta \Omega_G(A),\delta) $ and 
$ \psi: (\theta \Omega_G(A),\delta) \rightarrow (\theta \Omega_G(A),\partial) $ by 
\begin{equation*}
\phi = P + \frac{1}{2}\, R \partial = P - \frac{1}{2} \delta R, \qquad 
\psi = P + \frac{1}{2}\, \partial R = P - \frac{1}{2} \, R \delta.
\end{equation*} 
Using lemma \ref{cqlemma5} and lemma \ref{cqlemma4} one verifies that $ \phi $ and $ \psi $ are chain maps. Let us 
prove that $ \phi \psi $ is homotopic to the identity. According to lemma \ref{cqlemma5} one has
\begin{align*}
\phi\psi &= (P + \frac{1}{2} R \partial)(P + \frac{1}{2} \partial R) \\
&= P^2 + \frac{1}{2}(R \partial P + P \partial R) + \frac{1}{4}\, R \partial^2 R \\
&= P^2 - \frac{1}{2}(\delta R (F + \frac{1}{2} Q \delta) + (F - \frac{1}{2} \delta Q) R \delta) + \frac{1}{4} R^2 (\id - T) \\
&= P^2 - \frac{1}{2}(\delta R F + R F \delta) + \frac{1}{4} R^2 (\id - T).
\end{align*} 
Consider the first term in the last expression. By proposition \ref{Pprop} the map $ P $ is homotopic to the identity 
with respect to the boundary $ \delta $. Hence the same holds true for the 
chain map $ P^2 $. The second term is obviously homotopic to zero. The last term is homotopic to zero 
since $ R^2 $ is a chain map with respect to the boundary $ \delta $ according to lemma \ref{cqlemma5}. We conclude that $ \phi \psi $ 
is homotopic to the identity. In the same way one shows that $ \psi \phi $ is homotopic to the identity. \\
This finishes the proof of theorem \ref{homotopyeq}. \qed 

\section{Equivariant periodic cyclic homology}
\label{secHPdef}

In this section we define bivariant equivariant periodic cyclic homology for pro-$ G $-algebras. 
\begin{definition}
Let $ G $ be a locally compact group and let $ A $ and $ B $ be pro-$ G $-algebras. The bivariant equivariant periodic
cyclic homology of $ A $ and $ B $ is
\begin{equation*}
HP^G_*(A,B) =
H_*(\SHom_G(X_G(\mathcal{T}(A \hat{\otimes} \mathcal{K}_G)),X_G(\mathcal{T}(B \hat{\otimes} \mathcal{K}_G)))).
\end{equation*}
\end{definition}
There are some explanations in order. On the right hand side of this definition we take homology with respect to the 
usual boundary in a $ \Hom $-complex given by 
\begin{equation*}
\partial(\phi) = \phi \partial_A - (-1)^{|\phi|} \partial_B \phi
\end{equation*}
for a homogeneous element $ \phi\in \SHom_G(X_G(\mathcal{T}(A \hat{\otimes} \mathcal{K}_G)),X_G(\mathcal{T}(B \hat{\otimes} \mathcal{K}_G))) $ 
where $ \partial_A $ and $ \partial_B $ denote the boundary operators of 
$ X_G(\mathcal{T}(A \hat{\otimes} \mathcal{K}_G)) $ and 
$ X_G(\mathcal{T}(B \hat{\otimes} \mathcal{K}_G)) $, respectively. 
However, in order to take homology we have to check that we indeed obtain a supercomplex in this way 
since the equivariant $ X $-complexes are only paracomplexes. \\
From the definition of the equivariant $ X $-complex we know $ \partial_A^2 = \id - T $ and 
$ \partial_B^2 = \id - T $. Using these relations we compute 
\begin{equation*}
\partial^2(\phi) =
\phi\, \partial_A^2 + (-1)^{|\phi|} (-1)^{|\phi| - 1} \partial_B^2\,\phi 
= \phi(\id - T) - (\id - T)\phi = T \phi - \phi T 
\end{equation*}
and hence $ \partial^2(\phi) = 0 $ follows from proposition \ref{covparaadd}. Thus the failure of the individual differentials 
to satisfy $ \partial^2 = 0 $ is cancelled out in the $ \Hom $-complex. This shows that our definition of 
$ HP^G_* $ makes sense. \\
Let us discuss some basic properties of the equivariant homology groups defined above. 
Clearly $ HP^G_* $ is a bifunctor, contravariant in the first variable and covariant in
the second variable. As usual we define $ HP^G_*(A) =
HP^G_*(\mathbb{C}, A) $ to be the equivariant periodic cyclic
homology of $ A $ and $ HP_G^*(A) = HP^G_*(A, \mathbb{C}) $ to be
equivariant periodic cyclic cohomology. There is a natural product
\begin{equation*}
HP^G_i(A,B) \times HP^G_j(B,C) \rightarrow HP^G_{i + j}(A,C), \qquad
(x,y) \mapsto x\cdot y
\end{equation*}
induced by the composition of maps. This product is obviously associative.
Every equivariant homomorphism $ f: A \rightarrow B $ defines an
element in $ HP^G_0(A,B) $ denoted by $ [f] $. The element $ [\id] \in
HP^G_0(A,A) $ is simply denoted by $ 1 $ or $ 1_A $. An element
$ x \in HP^G_*(A,B) $ is called invertible if there
exists an element $ y \in HP^G_*(B,A) $ such that
$ x \cdot y = 1_A $ and $ y \cdot x = 1_B $. An invertible element
of degree zero will also be called an $ HP^G $-equivalence.
Such an element induces isomorphisms
$ HP^G_*(A,D) \cong HP^G_*(B,D) $ and $ HP^G_*(D,A) \cong HP^G_*(D,B) $ for
all $ G $-algebras $ D $.
An $ HP^G $-equivalence exists if and only if the paracomplexes
$ X_G(\mathcal{T}(A \hat{\otimes} \mathcal{K}_G)) $ and
$ X_G(\mathcal{T}(B \hat{\otimes} \mathcal{K}_G)) $ are covariantly homotopy
equivalent. 

\section{Homotopy invariance}

In this section we show that $ HP^G_* $ is invariant under smooth equivariant homotopies in both variables. \\
Let $ B $ be a pro-$ G $-algebra and consider the Fr\'echet algebra $ C^\infty[0,1] $ 
of smooth functions on the interval $ [0,1] $. We denote by $ B[0,1] $ the 
pro-$ G $-algebra $ B \hat{\otimes} C^\infty[0,1] $ where the action on 
$ C^\infty[0,1] $ is trivial. By definition a (smooth) equivariant homotopy is an equivariant 
homomorphism $ \Phi: A \rightarrow B[0,1] $ of pro-$ G $-algebras. 
Evaluation at a point $ t \in [0,1] $ yields an equivariant homomorphism 
$ \Phi_t: A \rightarrow B $. Two equivariant homomorphisms from $ A $ to $ B $ are called 
equivariantly homotopic if they can be connected by an equivariant homotopy. \\
A homology theory $ h_* $ for algebras is called homotopy invariant if 
the induced maps $ h_*(\phi_0) $ and $ h_*(\phi_1) $ are equal whenever $ \phi_0 $ and $ \phi_1 $ are homotopic homomorphisms. 
In our situation we will prove the following assertion. 
\begin{theorem}[Homotopy invariance] \label{homotopyinv} Let $ A $ and $ B $ be 
pro-$ G $-algebras
and let $ \Phi: A \rightarrow B[0,1] $ be a smooth equivariant homotopy.
Then the elements $ [\Phi_0] $ and $ [\Phi_1] $ in $ HP^G_0(A,B) $ are
equal. Hence the functor $ HP^G_* $ is homotopy invariant in both variables
with respect to smooth equivariant homotopies. \\
More generally the elements $ [\Phi_0] $ and $ [\Phi_1] $ in 
$ H_0(\SHom_G(X_G(A), X_G(B))) $ are equal provided $ A $ is quasifree. 
\end{theorem}
We recall that $ \theta^2 \Omega_G(A) $ is the 
paracomplex $ \Omega^0_G(A) \oplus \Omega^1_G(A) \oplus \Omega^2_G(A)/b(\Omega^3_G(A)) $ with the usual differential $ B + b $ 
and the grading into even and odd forms for any pro-$ G $-algebra $ A $. 
Clearly there is a natural map of paracomplexes $ \xi^2: \theta^2 \Omega_G(A) \rightarrow X_G(A) $. 
The first step in the proof of theorem \ref{homotopyinv} is to 
show that $ \xi^2 $ is a covariant homotopy equivalence provided 
$ A $ is equivariantly quasifree. \\
Let us consider the following more general situation. Assume that $ A $ is a pro-$ G $-algebra and 
let $ \nabla: \Omega^n(A) \rightarrow \Omega^{n + 1}(A) $ be an equivariant graded connection. 
Recall from definition \ref{defgradconn} that $ \nabla $ satisfies  
\begin{equation*}
\nabla(x \omega) = x \nabla(\omega), \qquad
\nabla(\omega x) = \nabla(\omega) x + (-1)^n \omega dx
\end{equation*}
for all $ x \in A $ and $ \omega \in \Omega^n(A) $. 
We extend $ \nabla $ to forms of higher degree by setting 
$ \nabla(a_0 da_1 \cdots da_m) = \nabla(a_0 da_1 \cdots da_n)da_{n + 1} \cdots da_m $. 
Moreover we put $ \nabla(\omega) = 0 $ if the degree of $ \omega $ is smaller than $ n $. 
Then we have 
\begin{equation*}
\nabla(a \omega) = a \nabla(\omega), \qquad
\nabla(\omega \eta) = \nabla(\omega) \eta  + (-1)^{|\omega|} \omega d\eta
\end{equation*}
for $ a \in A $ and differential forms $ \omega $ and $ \eta $.  
Let us define a covariant map 
$ \nabla_G: \theta \Omega_G(A) \rightarrow \theta \Omega_G(A) $ 
by the formula
\begin{equation*}
\nabla_G(f(s) \otimes \omega) = f(s) \otimes \nabla(\omega). 
\end{equation*}
\begin{prop} \label{homotopyinv1} 
Let $ A $ be a pro-$ G $-algebra and let 
$ \nabla: \Omega^n_G(A) \rightarrow \Omega^{n + 1}_G(A) $ be an equivariant graded connection. 
Then the covariant map $ [b,\nabla_G] = b \nabla_G + \nabla_G b $ 
is an idempotent operator on $ \theta \Omega_G(A) $ and defines a retraction 
for the natural map $ F^n \theta \Omega_G(A) \rightarrow \theta \Omega_G(A) $. \\
It follows that $ \theta\Omega_G(A) $ and $ \theta^n\Omega_G(A) $ 
are covariantly homotopy equivalent with respect to the Hochschild operators if $ A $ 
is $ n $-dimensional with respect to $ G $.  
\end{prop}
\proof Let us compute the commutator of $ b $ and $ \nabla_G $.
Take $ \omega \in \Omega^j(A) $ for $ j > n $. For 
$ a \in A $ we obtain
\begin{align*}
[b,\nabla_G]&(f(s) \otimes \omega da) = 
b(f(s) \otimes \nabla(\omega)da) + \nabla_G(b(f(s) \otimes \omega da)) \\
&= (-1)^{j + 1}(f(s) \otimes \nabla(\omega)a - 
f(s) \otimes (s^{-1} \cdot a) \nabla(\omega)) \\
&\qquad 
+ (-1)^j (\nabla_G(f(s) \otimes \omega a - f(s) \otimes (s^{-1} \cdot a)\omega)) \allowdisplaybreaks[2] \\
&= (-1)^j\biggl(f(s) \otimes (s^{-1} \cdot a)\nabla(\omega) - 
f(s) \otimes \nabla(\omega) a \\
&\qquad + f(s) \otimes \nabla(\omega a) - 
f(s) \otimes \nabla((s^{-1} \cdot a)\omega)\biggr) \\
&= (-1)^j\biggl(f(s) \otimes (s^{-1} \cdot a)\nabla(\omega) - 
f(s) \otimes \nabla(\omega) a + f(s) \otimes \nabla(\omega) a \\
&\qquad + (-1)^j f(s) \otimes \omega da - 
f(s) \otimes (s^{-1} \cdot a)\nabla(\omega)\biggr) \\
&= f(s) \otimes \omega da 
\end{align*}
Hence $ [b,\nabla_G] = \id $ on $ \Omega^j_G(A) $ for $ j > n $. Since 
$ [b,\nabla_G] $ commutes with $ b $ this holds also on
$ b(\Omega^{n + 1}_G(A)) $. Let us determine the behaviour 
of $ [b,\nabla_G] $ on $ \Omega^j_G(A) $ for $ j \leq n $. 
Clearly $ [b,\nabla_G] = 0 $ on $ \Omega^j_G(A) $ for $ j < n $ since 
$ \nabla_G $ vanishes on $ \Omega^j_G(A) $ and $ \Omega^{j - 1}_G(A) $ in this case.
On $ \Omega^n_G(A) $ we have $ [b,\nabla_G] = b \nabla_G $ because 
$ \nabla_G $ is zero on $ \Omega^{n - 1}_G(A) $. 
Hence 
\begin{equation*}
[b,\nabla_G][b,\nabla_G] = b\nabla_G b\nabla_G = b(\id - b \nabla_G) \nabla_G = b \nabla_G = [b,\nabla_G]\quad\text{on}\;\; \Omega^1_G(A)
\end{equation*}
and it follows that $ [b,\nabla_G] $ is idempotent. 
The range of the map $ [b,\nabla_G] = b \nabla_G $ restricted to 
$ \Omega^n_G(A) $ is contained in $ b(\Omega^{n + 1}_G(A)) $. Equality holds because 
$ [b,\nabla_G] $ is equal to the identity on $ b(\Omega^{n + 1}_G(A)) $ as we 
have seen before. It follows that $ [b,\nabla_G] $ maps $ \theta \Omega_G(A) $ to $ F^n \theta\Omega_G(A) $ and 
is a retraction of the natural map from $ F^n \theta\Omega_G(A) $ into $ \theta\Omega_G(A) $. 
Hence the map $ \id - [b, \nabla_G]: \theta^n\Omega_G(A) \rightarrow \theta\Omega_G(A) $ is inverse to the natural projection 
up to homotopy with respect to the Hochschild boundary. \qed 
\begin{prop} \label{homotopyinv1a} Let $ A $ be a $ G $-equivariantly quasifree pro-$ G $-algebra. 
Then the map $ \xi^2: \theta^2 \Omega_G(A) \rightarrow X_G(A) $ is a covariant 
homotopy equivalence. 
\end{prop}
\proof Since $ A $ is quasifree there exists by theorem \ref{qf} an equivariant 
graded connection $ \nabla: \Omega^1(A) \rightarrow \Omega^2(A) $.  
We use the covariant map $ \nabla_G $ defined above to construct an inverse of $ \xi^2 $ up to 
homotopy. In order to do this consider the commutator of 
$ \nabla_G $ with the boundary $ B + b $. 
Clearly we have $ [\nabla_G, B + b] = [\nabla_G, B] + [\nabla_G, b] $. Since 
$ [\nabla_G, B] $ has degree $ + 2 $ we see from proposition \ref{homotopyinv1} that 
$ \id - [\nabla_G, B + b] $ maps 
$ F^j \Omega_G(A) $ to $ F^{j + 1} \Omega_G(A) $ for all $ j \geq 1 $. 
This implies in particular that $ \id - [\nabla_G, B + b] $ descends to a covariant map $ \nu: X_G(A) \rightarrow \theta^2 \Omega_G(A) $. 
Using that $ \nabla_G $ is covariant we see that $ \nu $ is a chain map. Explicitly 
we have 
\begin{alignat*}{2}
\nu &= \id - \nabla_G d &&\text{on}\;\;\Omega^0_G(A) \\
\nu &= \id - [\nabla_G,b] = \id - b \nabla_G \qquad && \text{on}\;\; \Omega^1_G(A)/b(\Omega^2_G(A))
\end{alignat*}
and we deduce $ \xi^2 \nu = \id $. Moreover 
$ \nu \xi^2 = \id - [\nabla_G, B + b] $ is homotopic to the 
identity. This yields the assertion. \qed \\
Now let $ \Phi: A \rightarrow B[0,1] $ be an equivariant homotopy. 
The derivative of $ \Phi $ is an equivariant linear map 
$ \Phi': A \rightarrow B[0,1] $. 
If we view $ B[0,1] $ as a bimodule over itself the map 
$ \Phi' $ is a derivation with respect to $ \Phi $ in the sense that 
$ \Phi'(xy) = \Phi'(x) \Phi(y) + \Phi(x) \Phi'(y) $ for $ x,y \in A $. 
We define a covariant map $ \eta: \Omega^n_G(A) \rightarrow \Omega^{n - 1}_G(B) $ 
for $ n > 0 $ by 
\begin{equation*}
\eta(f(s) \otimes x_0dx_1 \dots dx_n) = 
\int_0^1 f(s) \otimes \Phi_t(x_0) \Phi'_t(x_1) d\Phi_t(x_2) \cdots d\Phi_t(x_n) dt.
\end{equation*}
Since integration is a bounded linear map we see that $ \eta $ is bounded. 
In addition we set $ \eta = 0 $ on $ \Omega^0_G(A) $. 
Using the fact that $ \Phi' $ is a derivation with respect to $ \Phi $ we 
compute 
\begin{align*}
\eta& b(f(s) \otimes x_0dx_1 \dots dx_n) = 
\int_0^1 f(s) \otimes \Phi_t(x_0 x_1) \Phi'_t(x_2) d\Phi_t(x_3) \cdots d\Phi_t(x_n) \\
&\qquad - f(s) \otimes \Phi_t(x_0) \Phi'_t(x_1 x_2) d\Phi_t(x_3) \cdots d\Phi_t(x_n) \\
&\qquad + f(s) \otimes \Phi_t(x_0) \Phi'_t(x_1) \Phi_t(x_2) d\Phi_t(x_3) 
\cdots d\Phi_t(x_n)\\
&\qquad - (-1)^n f(s) \otimes \Phi_t(x_0) \Phi'_t(x_1) (d\Phi_t(x_2) \cdots 
d\Phi_t(x_{n - 1})) \Phi_t(x_n) \\
&\qquad + (-1)^n f(s) \otimes \Phi_t((s^{-1} \cdot x_n)x_0) \Phi'_t(x_1) d\Phi_t(x_2) \cdots d\Phi_t(x_{n - 1}) dt\\
&= \int_0^1 (-1)^{n - 1} (f(s) \otimes \Phi_t(x_0) \Phi'_t(x_1) (d\Phi_t(x_2) \cdots 
d\Phi_t(x_{n - 1})) \Phi_t(x_n) \\
&\qquad - f(s) \otimes \Phi_t((s^{-1}\cdot x_n)x_0) \Phi'_t(x_1) d\Phi_t(x_2) \cdots 
d\Phi_t(x_{n - 1})) dt \\
&= - b \eta(f(s) \otimes x_0dx_1 \cdots dx_n).
\end{align*}
This implies that $ \eta $ maps $ b(\Omega^3_G(A)) $ into 
$ b(\Omega^2_G(B) $ and hence induces a covariant map 
$ \eta: \theta^2 \Omega_G(A) \rightarrow X_G(B) $. 
\begin{lemma} \label{homotopyinv2} We have $ X_G(\Phi_1) \xi^2 - X_G(\Phi_0) \xi^2 = \partial \eta + \eta \partial $. Hence the chain maps 
$ X_G(\Phi_t) \xi^2: \theta^2 \Omega_G(A) \rightarrow 
X_G(B) $ for $ t = 0,1 $ are covariantly homotopic. 
\end{lemma}
\proof We compute both sides on $ \Omega^j_G(A) $ for $ j = 0,1,2 $. For 
$ j = 0 $ we have 
\begin{equation*}
[\partial,\eta](f(s) \otimes x) = \eta(f(s) \otimes dx) 
= \int_0^1 f(s) \otimes \Phi'_t(x) dt = f(s) \otimes \Phi_1(x) - 
f(s) \otimes \Phi_0(x).
\end{equation*}
For $ j = 1 $ we get 
\begin{align*}
[\partial,\eta]&(f(s) \otimes x_0dx_1) = d\eta(f(s) \otimes x_0dx_1) 
+ \eta B(f(s) \otimes x_0dx_1) \\
&= \int_0^1 f(s) \otimes d(\Phi_t(x_0) \Phi'_t(x_1)) + 
f(s) \otimes \Phi'_t(x_0) d\Phi_t(x_1) - \\
&\qquad f(s) \otimes \Phi'_t(s^{-1} \cdot x_1) d\Phi_t(x_0) dt \\
&= \int_0^1 b(f(s) \otimes d\Phi_t(x_0) d\Phi'_t(x_1)) 
+ \frac{\partial}{\partial t}\biggl(f(s) \otimes \Phi_t(x_0) d\Phi_t(x_1)\biggr) dt \\
&= f(s) \otimes \Phi_1(x_0)d\Phi_1(x_1) - f(s) \otimes \Phi_0(x_0)d\Phi_0(x_1).
\end{align*}
Here we can forget about the term 
\begin{equation*}
\int_0^1 b(f(s) \otimes d\Phi_t(x_0) d\Phi'_t(x_1)) dt 
\end{equation*}
since the range of $ \eta $ is $ X_G(B) $. Finally, on 
$ \Omega^3_G(A)/b(\Omega^2_G(A)) $ we have $ \partial \eta + \eta \partial 
= \eta b + b \eta = 0 $ due to the computation above. \qed \\
Now we come back to the proof of theorem \ref{homotopyinv}. Let 
$ \Phi: A \rightarrow B[0,1] $ be an equivariant homotopy. 
Tensoring both sides with $ \mathcal{K}_G $ we obtain an equivariant homotopy 
$ \Phi \hat{\otimes}\mathcal{K}_G: A \hat{\otimes} \mathcal{K}_G \rightarrow (B \hat{\otimes} \mathcal{K}_G)[0,1] $. 
The map $ \Phi \hat{\otimes}\mathcal{K}_G $ 
induces an equivariant homomorphism 
$ \mathcal{T}(\Phi \hat{\otimes}\mathcal{K}_G): \mathcal{T}(A \hat{\otimes} \mathcal{K}_G) 
\rightarrow \mathcal{T}((B \hat{\otimes} \mathcal{K}_G)[0,1]) $. 
Now consider the equivariant linear map 
\begin{equation*}
l: B \hat{\otimes} \mathcal{K}_G \hat{\otimes} C^\infty[0,1] 
\rightarrow \mathcal{T}(B \hat{\otimes} \mathcal{K}_G) \hat{\otimes} C^\infty[0,1], 
\qquad l(b \otimes T \otimes f) = \sigma_{B \hat{\otimes} \mathcal{K}_G}(b \otimes T) \otimes f.
\end{equation*}
Since $ \sigma_{B \hat{\otimes} \mathcal{K}_G} $ is a lonilcur it follows that 
the same holds true for $ l $. Hence we obtain an equivariant 
homomorphism $ [[l]]: \mathcal{T}((B \hat{\otimes} \mathcal{K}_G)[0,1]) 
\rightarrow \mathcal{T}(B \hat{\otimes} \mathcal{K}_G)[0,1] $ due to 
proposition \ref{PeriodicTensorAlg}. 
Composition of $ \mathcal{T}(\Phi \hat{\otimes}\mathcal{K}_G) $ with the homomorphism 
$ [[l]] $ yields an equivariant homotopy 
$ \Psi = [[l]] \mathcal{T}(\Phi \hat{\otimes}\mathcal{K}_G): 
\mathcal{T}(A \hat{\otimes} \mathcal{K}_G) \rightarrow \mathcal{T}(B \hat{\otimes} \mathcal{K}_G)[0,1] $. 
From the definition of $ \Psi $ it follows easily that 
$ \Psi_t = \mathcal{T}(\Phi_t \hat{\otimes}\mathcal{K}_G) $ for all $ t $. 
Since $ \mathcal{T}(A \hat{\otimes} \mathcal{K}_G) $ is quasifree 
we can apply proposition \ref{homotopyinv1} and lemma \ref{homotopyinv2} 
to obtain $ [\Phi_0] = [\Phi_1] \in HP^G_0(A,B) $. The second assertion 
of theorem \ref{homotopyinv} follows directly from proposition \ref{homotopyinv1} and lemma \ref{homotopyinv2}. This finishes the proof of theorem 
\ref{homotopyinv}. \\
Let us note a formula for the chain homotopy $ h $ connecting $ X_G(\Phi_0) $ and $ X_G(\Phi_1) $ obtained 
above in the case that $ A $ is equivariantly quasifree. Since $ A $ is quasifree 
there exists according to theorem \ref{qf} an equivariant linear map $ \phi: A \rightarrow \Omega^2(A) $ 
satisfying $ \phi(xy) = \phi(x) y + x \phi(y) - dx dy $. Using the map $ \phi $ one obtains  
\begin{align*}
&h_0(f(s) \otimes x_0) = - \eta(f(s) \otimes \phi(x_0))\\
&h_1(f(s) \otimes x_0dx_1) = \eta(f(s) \otimes x_0dx_1) - \eta b(f(s) \otimes x_0 \phi(x_1))
\end{align*}
for the homotopy $ h: X_G(A) \rightarrow X_G(B) $. \\
As a first application of homotopy invariance we show that 
$ HP^G_* $ can be computed using arbitrary universal locally 
nilpotent extensions.
\begin{prop}\label{homotopyunivext}
Let $ 0 \rightarrow I \rightarrow R \rightarrow A \rightarrow 0 $ be
a universal locally nilpotent extension of the pro-$ G $-algebra
$ A $. Then $ X_G(R) $ is covariantly homotopy equivalent to 
$ X_G(\mathcal{T}A) $ in a canonical way. 
More precisely, any morphism of extensions $ (\xi,\phi,\id) $ from $ 0 \rightarrow \mathcal{J}A \rightarrow\mathcal{T}A \rightarrow A
\rightarrow 0 $ to 
$ 0 \rightarrow I \rightarrow R \rightarrow A \rightarrow 0 $ induces a covariant 
homotopy equivalence $ X_G(\phi): X_G(\mathcal{T}A) \rightarrow X_G(R) $. The class 
of this homotopy equivalence in $ H_*(\SHom_G(X_G(\mathcal{T}A), X_G(R))) $ is independent of the choice of $ \phi $.
\end{prop}
\proof From propositions \ref{UnivExt} and \ref{univext2} it follows that 
$ \phi: \mathcal{T}A \rightarrow R $ is an equivariant homotopy equivalence 
of algebras. Hence $ X_G(\phi): X_G(\mathcal{T}A) \rightarrow X_G(R) $ is a covariant 
homotopy equivalence due to theorem \ref{homotopyinv}. Since $ \phi $ is unique 
up to equivariant homotopy it follows that the class of this homotopy 
equivalence does not depend on the particular choice of $ \phi $. \qed \\
In particular there is a natural covariant homotopy equivalence between 
$ X_G(\mathcal{T}A) $ and $ X_G(A) $ if $ A $ itself is quasifree.

\section{Stability}

In this section we want to investigate stability properties of 
$ HP^G_* $. We will show that $ HP^G_* $ is stable with respect to
tensoring with the algebras $ l(b) $ associated to an equivariant bounded bilinear pairing $ b: W \times V \rightarrow \mathbb{C} $
that were introduced in section \ref{secborn}. \\
First we consider a special class of pairings. 
\begin{definition}
Let $ V $ and $ W $ be $ G $-modules. An equivariant bilinear pairing $ b: W \times V \rightarrow \mathbb{C} $ 
is called  admissible if there are subspaces $ N_W \subset W $ and $ N_V \subset V $ 
where the $ G $-action is trivial such that the restriction of $ b $ to $ N_W \times N_V $ is nonzero. 
\end{definition}
Now let $ A $ be a $ G $-algebra and let $ b: W \times V \rightarrow \mathbb{C} $ be an
admissible pairing. Let $ N_W \subset W $ and $ N_V \subset V $ be the corresponding subspaces. 
By assumption we may choose $ w \in N_W $ and $ v \in N_V $ such that $ b(w,v) = 1 $. Then $ p = v \otimes w $ 
is an element of $ l(b) $ and clearly $ p $ is $ G $-invariant. Consider the equivariant homomorphism 
$ \iota_A: A \rightarrow A \hat{\otimes} l(b), \iota_A(a) = a \otimes p $.  
\begin{theorem}\label{StabLemma} Let $ A $ be a pro-$ G $-algebra and
let $ b: W \times V \rightarrow \mathbb{C} $ be an admissible pairing. Then
the class $ [\iota_A] \in H_0(\SHom_G(X_G(\mathcal{T}A),
X_G(\mathcal{T}(A \hat{\otimes} l(b))))) $ is invertible.
\end{theorem}
\proof We have to find an inverse for $ [\iota_A] $. Our argument
is a generalization of a well-known proof of stability in the nonequivariant 
case. \\
First observe that the canonical equivariant linear map 
$ A \hat{\otimes} l(b) \rightarrow 
\mathcal{T}A \hat{\otimes} l(b) $ is a lonilcur and induces consequently 
an equivariant homomorphism $ \lambda_A: \mathcal{T}(A \hat{\otimes}
l(b)) \rightarrow \mathcal{T}A \hat{\otimes} l(b) $.
Define the map $ tr_A: X_G(\mathcal{T}A \hat{\otimes} l(b))
\rightarrow X_G(\mathcal{T}A) $ by
\begin{equation*}
tr_A(f(s) \otimes x \otimes T) = tr_s(T) f(s) \otimes x
\end{equation*}
and
\begin{equation*}
tr_A(f(s) \otimes x_0 \otimes T_0\, d(x_1 \otimes T_1)) = 
tr_s(T_0 T_1) f(s) \otimes x_0 dx_1.
\end{equation*}
Here we use the twisted trace $ tr_s: l(b) \rightarrow \mathbb{C} $ defined by 
\begin{equation*}
tr_s(v \otimes w) = b(w, s\cdot v) = b(s^{-1} \cdot w,v)
\end{equation*}
for $ v \otimes w \in V \cotimes W $ and $ s \in G $. \\
Now it is easily verified that
\begin{equation*}
tr_s(T_0 T_1) = tr_s((s^{-1}\cdot T_1)T_0)
\end{equation*}
for all $ T_0, T_1 \in l(b) $. \\
One checks that $ tr_A $ is a covariant map of
paracomplexes. We define $ \tau_A = tr_A \circ
X_G(\lambda_A) $ and claim that $ [\tau_A] $ is an inverse for 
$ [\iota_A] $. Using the relation $ p\, U_s = p $ one computes 
$ [\iota_A] \cdot [\tau_A] = 1 $. We have to show that $ [\tau_A] \cdot
[\iota_A] = 1 $. Consider the equivariant homomorphisms $ i_j: A
\hat{\otimes} l(b) \rightarrow A \hat{\otimes} l(b) \hat{\otimes} l(b) $ for $ j = 1,2 $ given by
\begin{align*}
&i_1(a \otimes T) = a \otimes T \otimes p \\
&i_2(a \otimes T) = a \otimes p \otimes T
\end{align*}
As before we see $ [i_1] \cdot [\tau_{A\hat{\otimes} l(b)}] = 1 $
and we determine $ [i_2] \cdot [\tau_{A \hat{\otimes} l(b)}] =
[\tau_A] \cdot [\iota_A] $. Let us show that the maps $ i_1 $ and 
$ i_2 $ are equivariantly homotopic. We shall define an invertible  multiplier $ \sigma $ of $ l(b) \cotimes l(b) $ such that 
conjugation with $ \sigma $ yields the natural coordinate flip of $ l(b) \cotimes l(b) $ sending $ k_1 \cotimes k_2 $ to 
$ k_2 \cotimes k_1 $ as follows. Using that $ l(b) \cotimes l(b) \cong V \cotimes W \cotimes V \cotimes W $ as $ G $-modules we set 
$$
\sigma \cdot (v_1 \otimes w_1 \otimes v_2 \otimes w_2) = v_2 \otimes w_1 \otimes v_1 \otimes w_2 
$$
and 
$$
(v_1 \otimes w_1 \otimes v_2 \otimes w_2) \cdot \sigma = v_2 \otimes w_1 \otimes v_1 \otimes w_2. 
$$
It is clear that these formulas define equivariant bounded linear maps $ l(b) \cotimes l(b) \rightarrow l(b) \cotimes l(b) $. 
Moreover we have $ \sigma \cdot (k l) = (\sigma \cdot k)l $, $ (kl) \cdot \sigma = k (l\cdot \sigma) $ and 
$ (k \cdot \sigma) l = k(\sigma \cdot l) $ for all $ k,l \in l(b) \cotimes l(b) $ which means by definition that $ \sigma $ is a multiplier of $ l(b) \cotimes l(b) $. 
We have  $ \sigma \cdot (\sigma \cdot k) = k = (k \cdot \sigma) \cdot \sigma $ and 
$ \Ad(\sigma)(k_1 \otimes k_2) = \sigma \cdot (k_1 \otimes k_2) \cdot \sigma 
= k_2 \otimes k_1 $. Consider for $ t \in [0,1] $ the invertible multiplier  $ \sigma_t = \cos(\pi t/2) \id + \sin(\pi t/2) \sigma $ with inverse 
given by  $ \sigma_t^{-1} = \cos(\pi t/2) \id - \sin(\pi t/2) \sigma $. The family 
$ \sigma_t $ depends smoothly on $ t $ and we have 
$ \sigma_0 = \id $ and $ \sigma_1 = \sigma $. Now the formula 
$ \Ad(\sigma_t)(k) = \sigma_t \cdot  k \cdot \sigma_t^{-1} $ defines equivariant 
homomorphisms $ \Ad(\sigma_t): l(b) \cotimes l(b) \rightarrow l(b) \cotimes l(b) $. We use $ \Ad(\sigma_t) $ to define an 
equivariant homomorphism $ h_t: A \hat{\otimes} l(b) \rightarrow 
A \hat{\otimes} l(b) \cotimes l(b) $ by $ h_t(a \otimes k) = a \otimes \Ad(\sigma_t)(k \otimes p) $. One computes 
$ h_0 = i_1 $ and $ h_1 = i_2 $ and the family $ h_t $ again depends 
smoothly on $ t $. Hence we have indeed defined a smooth homotopy
between $ i_1 $ and $ i_2 $. Theorem \ref{homotopyinv}  
yields $ [i_1] = [i_2] $ and hence $ [\tau_A] \cdot [\iota_A] = 1 $. \qed \\
Now we can prove the following stability theorem.
\begin{theorem}[Stability] Let $ A $ be a pro-$ G $-algebra and let
$ b: W \times V $ be any nonzero equivariant bilinear pairing. Then there exists
an invertible element in $ HP^G_0(A, A \hat{\otimes} l(b)) $.
Hence there are natural isomorphisms
\begin{equation*}
HP^G_*(A \hat{\otimes} l(b), B) \cong HP^G_*(A,B), \qquad
HP^G_*(A,B) \cong HP^G_*(A, B \hat{\otimes} l(b))
\end{equation*}
for all pro-$ G $-algebras $ A $ and $ B $.
\end{theorem}
\proof Consider the natural pairing $ \D(G) \times \D(G) \rightarrow \mathbb{C} $ used in the definition of $ \mathcal{K}_G $.  
The tensor product $ l(b) \cotimes \mathcal{K}_G $ is isomorphic to the algebra $ l(\D(G) \cotimes V, \D(G) \cotimes W) $ associated to 
the tensor product pairing. We have a natural equivariant isomorphism 
$ \alpha: \D(G) \cotimes V \cong \D(G) \cotimes V_\tau $ given by $ \alpha(f)(t) = t^{-1} \cdot f(t) $ 
where $ V_\tau $ is the space $ V $ equipped with the trivial $ G $-action. In the same way we obtain an equivariant 
isomorphism $ \D(G) \cotimes W \cong \D(G) \cotimes W_\tau $ which we will also denote by $ \alpha $. 
These isomorphisms are compatible with the pairings in the sense that  
$$
b(\alpha(g), \alpha(f)) = \int_G b(\alpha(g)(t), \alpha(f)(t)) dt = \int_G b(t^{-1} \cdot g(t), t^{-1} \cdot f(t)) dt = b(g, f) 
$$
for $ g \in \D(G) \cotimes W, f \in \D(G) \cotimes V $ where we use the fact that the pairing $ W \times V \rightarrow \mathbb{C} $ is equivariant. 
It follows that we obtain an equivariant isomorphism $ l(\D(G) \cotimes V, \D(G) \cotimes W) \rightarrow l(\D(G) \cotimes V_\tau, \D(G) \cotimes W_\tau) $ given 
by $ \alpha(f \cotimes g) = \alpha(f) \cotimes \alpha(g) $. 
In other words, we have an equivariant isomorphism of  $ G $-algebras
\begin{equation*}
\mathcal{K}_G \cotimes l(b) \cong \mathcal{K}_G \cotimes l(b_\tau)
\end{equation*}
where $ b_\tau = b: W_\tau \times V_\tau \rightarrow \mathbb{C} $. 
Now we can apply theorem \ref{StabLemma} with $ A $ replaced by 
$ A \hat{\otimes} \mathcal{K}_G $ and $ l(b) $ replaced by
$ l(b_\tau) $ to obtain the assertion. \qed \\
As an application of theorem \ref{StabLemma} we obtain a simpler description of $ HP^G_* $ if $ G $ is a compact group.
\begin{prop} \label{defcomp} Let $ G $ be a compact group. Then we have 
a natural isomorphism
\begin{equation*}
HP^G_*(A,B) \cong H_*(\SHom_G(X_G(\mathcal{T}A),X_G(\mathcal{T}B)))
\end{equation*}
for all pro-$ G $-algebras $ A $ and $ B $.
\end{prop}
\proof If $ G $ is compact the trivial one-dimensional representation is contained in 
$ \D(G) $. Hence the pairing used in the definition of $ \mathcal{K}_G $ is admissible in this case.\qed 

\section{Excision}

The goal of this section is the proof of excision in equivariant periodic cyclic homology. We 
consider an extension 
\begin{equation}\label{eqext}
   \xymatrix{
     K\;\; \ar@{>->}[r]^{\iota} & E \ar@{->>}[r]^{\pi} & Q \ar@/^/@{.>}[l]^{\sigma}
     }
\end{equation}
of pro-$ G $-algebras where $ \sigma: Q \rightarrow E $ is an equivariant linear 
splitting for the quotient map $ \pi: E \rightarrow Q $.\\
Let $ X_G(\mathcal{T}E:\mathcal{T}Q) $ be the kernel of the map $ X_G(\mathcal{T}\pi):
X_G(\mathcal{T}E) \rightarrow X_G(\mathcal{T}Q)) $ induced by $ \pi $. The splitting $ \sigma $ 
yields a direct sum decomposition 
$ X_G(\mathcal{T}E) = X_G(\mathcal{T}E:\mathcal{T}Q) \oplus X_G(\mathcal{T}Q) $ 
of covariant pro-modules. The resulting extension 
$$ 
   \xymatrix{
      X_G(\mathcal{T}E:\mathcal{T}Q)\;\; \ar@{>->}[r] & X_G(\mathcal{T}E) \ar@{->>}[r] & X_G(\mathcal{T}Q) 
     }
$$ 
of paracomplexes induces long exact sequences in homology in both variables. 
Moreover there is a natural covariant map 
$ \rho: X_G(\mathcal{T}K) \rightarrow X_G(\mathcal{T}E:\mathcal{T}Q) $ of paracomplexes.
Our main result is the following generalized excision theorem. 
\begin{theorem}\label{Excision2} The map $ \rho: X_G(\mathcal{T}K) \rightarrow
X_G(\mathcal{T}E:\mathcal{T}Q) $ is a covariant homotopy equivalence.
\end{theorem}
As a consequence we get excision in equivariant periodic cyclic homology.
\begin{theorem}[Excision]\label{Excision} Let $ A $ be a pro-$ G $-algebra and let 
$ (\iota, \pi): 0 \rightarrow K \rightarrow E \rightarrow Q \rightarrow 0 $ be
an extension of pro-$ G $-algebras with a linear splitting. Then there are two natural exact sequences
\begin{equation*}
\xymatrix{
 {HP^G_0(A,K)\;} \ar@{->}[r] \ar@{<-}[d] &
      HP^G_0(A,E) \ar@{->}[r] &
        HP^G_0(A,Q) \ar@{->}[d] \\
   {HP^G_1(A,Q)\;} \ar@{<-}[r] &
    {HP^G_1(A,E)}  \ar@{<-}[r] &
     {HP^G_1(A,K)} \\
}
\end{equation*}
and
\begin{equation*}
\xymatrix{
    {HP^G_0(Q,A)\;} \ar@{->}[r] \ar@{<-}[d] &
       HP^G_0(E,A) \ar@{->}[r] &
          HP^G_0(K,A) \ar@{->}[d] \\
    {HP^G_1(K,A)\;} \ar@{<-}[r] &
      {HP^G_1(E,A)}  \ar@{<-}[r] &
        {HP^G_1(Q,A)} \\
}
\end{equation*}
The horizontal maps in these diagrams are induced by the maps
in the extension.
\end{theorem} 
We point out that in theorem \ref{Excision} we only require a pro-linear splitting for the quotient 
homomorphism $ \pi: E \rightarrow Q $. Let us first show how theorem \ref{Excision2} implies theorem \ref{Excision}. 
Tensoring the extension given in theorem \ref{Excision} with $ \mathcal{K}_G $ yields 
an extension 
\begin{equation}\label{eqext2}
   \xymatrix{
     K \cotimes \mathcal{K}_G\;\; \ar@{>->}[r] & E \cotimes \mathcal{K}_G \ar@{->>}[r] & Q \cotimes \mathcal{K}_G 
     }
\end{equation}
of pro-$ G $-algebras with a linear splitting. Due to lemma 
\ref{relproj} the pro-$ G $-module $ Q \cotimes \mathcal{K}_G $ is relatively projective. It follows that we obtain in fact an 
equivariant linear splitting for extension (\ref{eqext2}). Now we can apply theorem \ref{Excision2} to this extension and obtain 
the claim by considering long exact sequences in homology. \\
Our proof of theorem \ref{Excision2} is an adaption of the method used in ~\cite{Meyer}
to prove excision in cyclic homology theories. Consider the extension (\ref{eqext}) and let $ \mathfrak{L} \subset \mathcal{T}E $ be the left ideal generated by
$ K \subset \mathcal{T}E $. Using proposition \ref{Omegaiso} we see that 
\begin{equation}\label{eq2}
(\mathcal{T}E)^+ \hat{\otimes} K \rightarrow \mathfrak{L}, \qquad
x \otimes k \mapsto x \circ k
\end{equation}
is an equivariant linear isomorphism. Moreover we obtain from this description an equivariant linear retraction for the inclusion
 $ \mathfrak{L} \rightarrow  \mathcal{T}E $. Clearly $ \mathfrak{L} $ is a pro-$ G $-algebra since the ideal
$ K \subset E $ is $ G $-invariant. The natural projection 
$ \tau_E: \mathcal{T}E \rightarrow E $ induces an equivariant 
homomorphism $ \tau: \mathfrak{L} \rightarrow K $ and $ \sigma_E $ restricted 
to $ K $ is an equivariant linear splitting for $ \tau $. 
Hence we obtain an extension 
\begin{equation*}
  \xymatrix{
     N\;\; \ar@{>->}[r] & \mathfrak{L} \ar@{->>}[r]^\tau & K 
     }
\end{equation*}
of pro-$ G $-algebras. The inclusion $ \mathfrak{L} 
\rightarrow \mathcal{T}E $ induces a morphism of extensions
from $ 0 \rightarrow N \rightarrow \mathfrak{L} \rightarrow K \rightarrow 0 $ 
to $ 0 \rightarrow \mathcal{J}E \rightarrow \mathcal{T}E \rightarrow E 
\rightarrow 0 $. In particular we have a natural equivariant homomorphism 
$ i: N \rightarrow \mathcal{J}E $ and it is easy to see that there exists 
an equivariant linear map $ r: \mathcal{J}E \rightarrow N $ 
such that $ r i = \id $. Using this retraction we want to show that $ N $ 
is locally nilpotent. If $ l: N \rightarrow C $ is an equivariant linear 
map with constant range $ C $ we compute $ l m_N^n = l p i m_N^n = 
l p m_{\mathcal{J}E}^n i^{\hat{\otimes} n} $ where $ m_N $ and $ m_{\mathcal{J}E} $ 
denote the multiplication maps in $ N $ and $ \mathcal{J}E $, respectively. Since 
$ l p: \mathcal{J}E \rightarrow C $ is an equivariant linear map with constant range the claim 
follows from the fact that $ \mathcal{J}E $ is locally nilpotent. \\
We will establish theorem \ref{Excision2} by showing
\begin{theorem}\label{Excision3} With the notations as above we have
\begin{bnum}
\item[a)] The pro-$ G $-algebra $ \mathfrak{L} $ is quasifree.
\item[b)] The inclusion map $ \mathfrak{L} \rightarrow \mathcal{T}E $ induces a
covariant homotopy equivalence $ \psi: X_G(\mathfrak{L}) \rightarrow
X_G(\mathcal{T}E:\mathcal{T}Q) $.
\end{bnum}
\end{theorem}
Let us indicate how theorem \ref{Excision3} implies theorem \ref{Excision2}.
The map $ \rho $ is the composition of the natural maps
$ X_G(\mathcal{T}K) \rightarrow X_G(\mathfrak{L}) $ and
$ X_G(\mathfrak{L}) \rightarrow X_G(\mathcal{T}E:\mathcal{T}Q) $. Since
$ \mathfrak{L} $ is quasifree by part a) it follows that $ 0 \rightarrow 
N \rightarrow \mathfrak{L} \rightarrow K \rightarrow 0 $ is a universal locally nilpotent 
extension of $ K $. Hence the first map is a covariant homotopy 
equivalence due to proposition \ref{homotopyunivext}. The second map 
is a covariant homotopy equivalence by part b). It follows that $ \rho $
itself is a covariant homotopy equivalence. \\
We need some notation. The equivariant linear section
$ \sigma: Q \rightarrow E $ induces an equivariant linear map
$ \sigma_L: \Omega^n(Q) \rightarrow \Omega^n(E) $ definded by
\begin{equation*}
\sigma_L(q_0 dq_1 \cdots dq_n) =
\sigma(q_0) d\sigma(q_1) \dots d\sigma(q_n).
\end{equation*}
Here $ \sigma $ is extended to an equivariant linear 
map $ Q^+ \rightarrow E^+ $ in the obvious way by requiring $ \sigma(1) = 1 $. \\
We also need a right-handed version of the map $ \sigma_L $. 
In order to explain this correctly consider first an arbitrary pro-$ G $-algebra 
$ A $. There is a natural equivariant isomorphism 
$ \Omega^1(A) \cong A \hat{\otimes} A^+ $ 
of right $ A $-modules. This follows easily from the description of $ \Omega^1(A) $ 
as the kernel of the multiplication map $ A^+ \hat{\otimes} A^+ \rightarrow A^+ $.
More generally we obtain equivariant linear isomorphisms $ \Omega^n(A) \cong 
A^{\hat{\otimes} n} \hat{\otimes} A^+ $ for all $ n $. 
Using these identifications we define the equivariant linear map 
$ \sigma_R: \Omega(Q) \rightarrow \Omega(E) $ by
\begin{equation*}
\sigma_R(dq_1 \cdots dq_n q_{n + 1}) =
d\sigma(q_1) \dots d\sigma(q_n) \sigma(q_{n + 1})
\end{equation*}
which is the desired right-handed version of $ \sigma_L $. As in ~\cite{Meyer} we obtain the following assertion. 
\begin{lemma}\label{exlem1} The following maps are equivariant linear 
isomorphisms:
\begin{align*}
\mu_L: &(\mathcal{T}Q)^+ \oplus (\mathcal{T}E)^+ \hat{\otimes} K \hat{\otimes} (\mathcal{T}Q)^+
\rightarrow (\mathcal{T}E)^+ \\
&q_1 \oplus (x \otimes k \otimes q_2) \mapsto
\sigma_L(q_1) + x \circ k \circ \sigma_L(q_2) \\
&\qquad \\
\mu_R: &(\mathcal{T}Q)^+ \oplus (\mathcal{T}Q)^+ \hat{\otimes} K \hat{\otimes} (\mathcal{T}E)^+ \rightarrow (\mathcal{T}E)^+ \\
&q_1 \oplus (q_2 \otimes k \otimes x) \mapsto
\sigma_R(q_1) + \sigma_R(q_2) \circ k \circ x
\end{align*}
\end{lemma}
Equation (\ref{eq2}) and lemma \ref{exlem1} yield an equivariant linear 
isomorphism
\begin{equation}\label{eq3}
\mathfrak{L}^+ \hat{\otimes} (\mathcal{T}Q)^+ \cong (\mathcal{T}E)^+, \qquad
l \otimes q \mapsto l \circ \sigma_L(q).
\end{equation}
This isomorphism is obviously left $ \mathfrak{L} $-linear and it follows that 
$ (\mathcal{T}E)^+ $ is a free left $ \mathfrak{L} $-module. Furthermore we
get from lemma \ref{exlem1}
\begin{equation*}
(\mathcal{T}Q)^+ \hat{\otimes} K \hat{\otimes} \mathfrak{L}^+ \cong 
(\mathcal{T}Q)^+ \hat{\otimes} K \oplus
(\mathcal{T}Q)^+ \hat{\otimes} K \hat{\otimes} (\mathcal{T}E)^+ \hat{\otimes} K 
\cong (\mathcal{T}E)^+ \hat{\otimes} K \cong \mathfrak{L}.
\end{equation*}
It follows that the equivariant linear map 
\begin{equation}\label{eq4}
(\mathcal{T}Q)^+ \hat{\otimes} K \hat{\otimes} \mathfrak{L}^+ \rightarrow \mathfrak{L}, \qquad
q \otimes k \otimes l \mapsto \sigma_R(q) \circ k \circ l
\end{equation}
is an isomorphism. This map is right $ \mathfrak{L} $-linear and we see that
$ \mathfrak{L} $ is a free right $ \mathfrak{L} $-module. \\
Denote by $ \mathfrak{J} $ the kernel of the map
$ \mathcal{T}\pi: \mathcal{T}E \rightarrow \mathcal{T}Q $. Using again 
lemma \ref{exlem1} we see that
\begin{equation}\label{eq5}
(\mathcal{T}Q)^+ \hat{\otimes} K \hat{\otimes} (\mathcal{T}E)^+ \cong \mathfrak{J}, \qquad
q \otimes k \otimes x \mapsto \sigma_R(q) \circ k \circ x
\end{equation}
is a right $ \mathcal{T}E $-linear isomorphism. In a similar way we 
have a left $ \mathcal{T}E $-linear isomorphism
\begin{equation*}
(\mathcal{T}E)^+ \hat{\otimes} K \hat{\otimes} (\mathcal{T}Q)^+ \cong \mathfrak{J}, \qquad
x \otimes k \otimes q \mapsto x \circ k \circ \sigma_L(q). 
\end{equation*}
Together with equation (\ref{eq2}) this yields
\begin{equation}\label{eq6}
\mathfrak{L} \hat{\otimes} (\mathcal{T}Q)^+ \cong \mathfrak{J}, \qquad
l \otimes q \mapsto l \circ \sigma_L(q),
\end{equation}
and using equation (\ref{eq3}) we get
\begin{equation}\label{eq7}
\mathfrak{L} \hat{\otimes}_{\mathfrak{L}^+} (\mathcal{T}E)^+ \cong \mathfrak{J},
\qquad l \otimes x \mapsto l \circ x.
\end{equation}
Now one constructs a free resolution of the $ \mathfrak{L} $-bimodule
$ \mathfrak{L}^+ $. First let $ A $ be any pro-$ G $-algebra and consider
the extension of $ A $-bimodules in $ \pro(G\LMod) $ given by
\begin{equation*}
\xymatrix{
      {B^A_\bullet \colon \ }
        {\Omega^1(A)\;\;} \ar@<1ex>@{>->}[r]^-{\alpha_1} &
          {\;A^+ \cotimes A^+} \ar@<1ex>@{->>}[r]^-{\alpha_0}
            \ar@<1ex>@{.>>}[l]^-{h_1} &
            {\;\;A^+} \ar@<1ex>@{>.>}[l]^-{h_0}
      } 
\end{equation*}
Here the maps are defined as follows:
\begin{alignat*}{2}
&\alpha_1(x Dy z) = x y \otimes z - x \otimes y z, 
\qquad &\alpha_0(x \otimes y) = xy \\
&h_1(x \otimes y) = Dx y, \qquad \qquad &h_0(x) = 1 \otimes x
\end{alignat*}
It is easy to check that $ \alpha h + h \alpha = \id $. 
The complex $ B_\bullet^A $ is a projective 
resolution of the $ A $-bimodule $ A^+ $ in $ \pro(G\LMod) $ iff $ A $ is quasifree. 
Define a subcomplex $ P_\bullet $ of $  B_\bullet^{\mathcal{T}E} $ as follows:
\begin{align*}
&P_0 = (\mathcal{T}E)^+ \hat{\otimes} \mathfrak{L} +
\mathfrak{L}^+ \hat{\otimes} \mathfrak{L}^+ \subset
(\mathcal{T}E)^+ \hat{\otimes} (\mathcal{T}E)^+ \\
&P_1 = (\mathcal{T}E)^+ D\mathfrak{L} \subset \Omega^1(\mathcal{T}E).
\end{align*}
There exists an equivariant linear retraction $ B_\bullet^{\mathcal{T}E}
\rightarrow P_\bullet $ for the inclusion $ P_\bullet \rightarrow  B_\bullet^{\mathcal{T}E} $. 
Since $ \mathfrak{L} $ is a left ideal in $ \mathcal{T}E $ we see that
the boundary operators in $ B_\bullet^{\mathcal{T}E} $ restrict 
to $ P_\bullet $ and turn $ P_1 \rightarrow P_0 \rightarrow \mathfrak{L}^+ $ into a complex. It is clear that $ P_0 $ and $ P_1 $ inherit a natural 
$ \mathfrak{L} $-bimodule structure from $ B_0^{\mathcal{T}E} $ and 
$ B_1^{\mathcal{T}E} $, respectively. 
Moreover the homotopy $ h $ restricts to a contracting homotopy for the 
complex $ P_1 \rightarrow P_0 \rightarrow \mathfrak{L}^+ $. Hence 
$ P_\bullet $ is a resolution of $ \mathfrak{L}^+ $ by 
$ \mathfrak{L} $-bimodules in $ \pro(G\LMod) $. 
Next we show that the $ \mathfrak{L} $-bimodules $ P_0 $ and $ P_1 $ are free. 
Using equation (\ref{eq3}) we obtain the isomorphism
\begin{align}\label{eq8}
\mathfrak{L}^+ \hat{\otimes}& \mathfrak{L}^+ \oplus 
\mathfrak{L}^+ \hat{\otimes} \mathcal{T}Q \hat{\otimes} \mathfrak{L} \cong P_0, \\
& (l_1 \otimes l_2) \oplus (l_3 \otimes q \otimes l_4) \mapsto 
l_1 \otimes l_2 + (l_3 \circ \sigma_L(q)) \otimes l_4 \notag.
\end{align}
Since $ \mathfrak{L} $ is a free right $ \mathfrak{L} $-module by 
(\ref{eq4}) we see that $ P_0 $ is a free $ \mathfrak{L} $-bimodule. 
Now consider $ P_1 $. 
We claim that 
\begin{equation*}
P_1 = \Omega^1(\mathcal{T}E) \circ K + (\mathcal{T}E)^+ DK.
\end{equation*}
The inclusion $ (\mathcal{T}E)^+ DK \subset P_1 $ is clear and it is easy to see that 
$ \Omega^1(\mathcal{T}E) \circ K \subset P_1 $. Conversely, for 
$ x_0 D(x_1 \circ k) \in P_1 $ with $ x_0,x_1 \in (\mathcal{T}E)^+ $ we compute
\begin{equation*}
x_0 D(x_1 \circ k) = x_0 (D x_1) \circ k + x_0 \circ x_1 Dk  
\end{equation*}
which is contained in $ \Omega^1(\mathcal{T}E) \circ K + (\mathcal{T}E)^+ DK $. 
This yields the claim. Under the isomorphism $ \Omega^1(\mathcal{T}E) 
\cong (\mathcal{T}E)^+ \hat{\otimes} E \hat{\otimes} (\mathcal{T}E)^+ $ from 
proposition \ref{Omegaiso} the space $ \Omega^1(\mathcal{T}E) \circ K $ corresponds to 
$ (\mathcal{T}E)^+ \hat{\otimes} E \hat{\otimes} (\mathcal{T}E)^+ \circ K 
= (\mathcal{T}E)^+ \hat{\otimes} E \hat{\otimes} \mathfrak{L} $ and 
$ (\mathcal{T}E)^+ DK $ corresponds to $ (\mathcal{T}E)^+ \hat{\otimes} K \hat{\otimes} 1 $. 
Hence 
\begin{align} \label{eq9}
((\mathcal{T}E)^+ \hat{\otimes} K &\hat{\otimes} \mathfrak{L}^+) \oplus 
((\mathcal{T}E)^+ \hat{\otimes} Q \hat{\otimes} \mathfrak{L}) \rightarrow P_1, \\
&(x_1 \otimes k \otimes l_1) \oplus (x_2 \otimes q \otimes l_2)
\mapsto x_1 Dk l_1 + x_2 D\sigma(q) l_2 \notag 
\end{align}
is an equivariant linear isomorphism. Since $ (\mathcal{T}E)^+ $ 
is a free left $ \mathfrak{L} $-module by equation (\ref{eq3}) and 
$ \mathfrak{L} $ is a free right $ \mathfrak{L} $-module by equation (\ref{eq4}) 
we deduce that $ P_1 $ is a free $ \mathfrak{L} $-bimodule. 
Consequently we have established that 
$ P_\bullet $ is a
free $ \mathfrak{L} $-bimodule resolution of $ \mathfrak{L}^+ $ in the
category $ \pro(G\LMod) $. According to theorem \ref{qf} this finishes the proof 
of part a) of theorem \ref{Excision3}. \\
We need some more notation. Let $ X^\beta_G(\mathcal{T}E) $ be the complex 
obtained from 
$ X_G(\mathcal{T}E) $ by replacing the differential 
$ \partial_1: X^1_G(\mathcal{T}E) \rightarrow X^0_G(\mathcal{T}E) $ by zero. In the same way we proceed for $ X_G(\mathcal{T}E:\mathcal{T}Q) $. 
Moreover let $ M $ be an $ \mathfrak{L} $-bimodule 
in $ \pro(G\LMod) $. We define the covariant module 
$ (\mathcal{O}_G \hat{\otimes} M)/[\;,\,]_G $ as the quotient of 
$ \mathcal{O}_G \hat{\otimes} M $ by twisted commutators 
$ f(s) \otimes ml - f(s) \otimes (s^{-1} \cdot l)m $ where $ l \in \mathfrak{L} $ 
and $ m \in M $. \\
Now we continue the proof of theorem \ref{Excision3}. 
The inclusion $ P_\bullet \rightarrow B_\bullet^{\mathcal{T}E} $ is an
$ \mathfrak{L} $-bimodule homomorphism and induces a chain map
\begin{equation*}
\phi: (\mathcal{O}_G \hat{\otimes} P_\bullet)/[\;,\,]_G
\rightarrow (\mathcal{O}_G \hat{\otimes} B_\bullet^{\mathcal{T}E})/[\;,\,]_G \cong
X^\beta_G(\mathcal{T}E) \oplus \mathcal{O}_G[0].
\end{equation*}
Let us determine the image of $ \phi $. We use equations (\ref{eq8}) and 
(\ref{eq6}) to obtain
\begin{align*}
(\mathcal{O}_G \hat{\otimes} P_0)&/[\;,\,]_G \cong
\mathcal{O}_G \hat{\otimes} (\mathfrak{L}^+ \oplus
\mathfrak{L} \hat{\otimes}  \mathcal{T}Q) \\
&\cong \mathcal{O}_G \oplus
(\mathcal{O}_G \hat{\otimes} \mathfrak{L} \hat{\otimes} (\mathcal{T}Q)^+)
\cong \mathcal{O}_G \oplus (\mathcal{O}_G \hat{\otimes} \mathfrak{J})
\subset \mathcal{O}_G \hat{\otimes} (\mathcal{T}E)^+. 
\end{align*}
Using equations (\ref{eq9}) and (\ref{eq7}) we get 
\begin{align*}
(\mathcal{O}_G \hat{\otimes} P_1)&/[\;,\,]_G\cong
\mathcal{O}_G \hat{\otimes} ((\mathcal{T}E)^+ \hat{\otimes} K) \oplus
\mathcal{O}_G \hat{\otimes} (\mathfrak{L}\hat{\otimes}_{\mathfrak{L}^+} (\mathcal{T}E)^+ \hat{\otimes} Q) \\
&\cong \mathcal{O}_G \hat{\otimes} ((\mathcal{T}E)^+ \hat{\otimes} K) \oplus
\mathcal{O}_G \hat{\otimes} \mathfrak{J} \hat{\otimes} Q \\
&\cong \mathcal{O}_G \hat{\otimes} ((\mathcal{T}E)^+ DK +
\mathfrak{J} D \sigma(Q)) \subset \Omega^1_G(\mathcal{T}E).
\end{align*}
This implies that $ \phi $ induces a covariant isomorphism of
chain complexes 
\begin{equation*}
(\mathcal{O}_G \hat{\otimes} P_\bullet)/[\;,\,]_G \cong
X^\beta_G(\mathcal{T}E:\mathcal{T}Q) \oplus \mathcal{O}_G[0].
\end{equation*}
With these preparations we can prove part b) of theorem 
\ref{Excision3}. 
\begin{prop}
The natural map $ \psi: X_G(\mathfrak{L}) \rightarrow
X_G(\mathcal{T}E:\mathcal{T}Q) $ is split injective and we have
\begin{equation*}
X_G(\mathcal{T}E:\mathcal{T}Q) = X_G(\mathfrak{L}) \oplus C_\bullet
\end{equation*}
with a covariantly contractible paracomplex $ C_\bullet $.
Hence $ X_G(\mathcal{T}E:\mathcal{T}Q) $ and $ X_G(\mathfrak{L}) $
are covariantly homotopy equivalent.
\end{prop}
\proof The standard resolution $ B_\bullet^\mathfrak{L} $ of
$ \mathfrak{L}^+ $ is a subcomplex
of $ P_\bullet $. Since $ P_\bullet $ itself is a free $ \mathfrak{L} $-bimodule resolution of $ \mathfrak{L}^+ $ the inclusion map 
$ f_\bullet: B_\bullet^\mathfrak{L} \rightarrow P_\bullet $ is a homotopy equivalence. Explicitly 
set $ M_0 = \mathfrak{L}^+ \hat{\otimes} \mathcal{T}Q \hat{\otimes} \mathfrak{L} $ and define $ g: M_0 \rightarrow P_0 $ by
\begin{equation*}
g(l_1 \otimes q \otimes l_2) = l_1 \circ \sigma_L(q) \otimes l_2 -
l_1 \otimes \sigma_L(q) \circ l_2.
\end{equation*}
Using equation (\ref{eq8}) it is easy to check that
$ f_0 \oplus g: \mathfrak{L}^+ \hat{\otimes} \mathfrak{L}^+ \oplus M_0
\rightarrow P_0 $ is an isomorphism. Furthermore we have
$ \alpha_0 g = 0 $. Since the complex $ P_\bullet $ is exact this
implies $ P_1 = \ker \alpha_0 \cong \Omega^1(\mathfrak{L}) \oplus M_0 $.
Set $ M_1 = M_0 $ and define the boundary $ M_1 \rightarrow M_0 $ to be
the identity map. The complex $ M_\bullet $ of $ \mathfrak{L} $-bimodules
is obviously contractible and
$ P_\bullet \cong B_\bullet^\mathfrak{L} \oplus M_\bullet $.
Applying the functor $ (\mathcal{O}_G \hat{\otimes} -)/[\;,\,]_G $ we obtain covariant
isomorphisms
\begin{align*}
X^\beta_G(\mathcal{T}E:\mathcal{T}Q) &\oplus \mathcal{O}_G[0] \cong
(\mathcal{O}_G \hat{\otimes} P_\bullet)/[\;,\,]_G
\cong (\mathcal{O}_G \hat{\otimes} B_\bullet^\mathfrak{L})/[\;,\,]_G \oplus
(\mathcal{O}_G \hat{\otimes} M_\bullet)/[\;,\,]_G \\
&\cong X^\beta_G(\mathfrak{L}) \oplus \mathcal{O}_G[0] \oplus
(\mathcal{O}_G \hat{\otimes} M_\bullet)/[\;,\,]_G.
\end{align*}
One checks that the two copies of $ \mathcal{O}_G $ are identified under this isomorphism. Moreover the map $ X^\beta_G(\mathfrak{L})
\rightarrow X^\beta_G(\mathcal{T}E:\mathcal{T}Q) $ arising from these
identifications is equal to $ \psi $. Hence $ \psi $ is split injective. Let 
$ C_\bullet $ be the image
of $ (\mathcal{O}_G \hat{\otimes} M_\bullet)/[\;,\,]_G $ in
$ X^\beta_G(\mathcal{T}E:\mathcal{T}Q) $. One checks that 
$ C_0 $ is the range of the map 
\begin{equation*}
\mathcal{O}_G \hat{\otimes} \mathfrak{L} \hat{\otimes} 
\mathcal{T}Q \rightarrow X^0_G(\mathcal{T}E), \qquad 
f(s) \otimes l \otimes q \mapsto f(s) \otimes l \circ s_L(q) - 
f(s) \otimes (s^{-1} \cdot s_L(q)) \circ l 
\end{equation*}
and that $ C_1 $ is the range of the map 
\begin{equation*}
\mathcal{O}_G \hat{\otimes} \mathfrak{L} \hat{\otimes} 
\mathcal{T}Q \rightarrow X^1_G(\mathcal{T}E), \qquad 
f \otimes l \otimes q \mapsto f \otimes l Ds_L(q).
\end{equation*}
The boundary $ C_1 \rightarrow C_0 $ is the boundary induced from 
$ X_G(\mathcal{T}E:\mathcal{T}Q) $. 
On the other hand the boundary $ \partial_0: X^0_G(\mathcal{T}E:\mathcal{T}Q)
\rightarrow X^1_G(\mathcal{T}E:\mathcal{T}Q) $ does not vanishes on $ C_0 $. However, we have $ \partial^2 = \id - T $ and 
this implies that $ C_\bullet $ is a 
sub-paracomplex of $ X_G(\mathcal{T}E:\mathcal{T}Q) $. 
Since $ \psi $ is compatible 
with $ \partial_0 $ we obtain the desired direct sum decomposition
\begin{equation*}
X_G(\mathcal{T}E:\mathcal{T}Q) \cong X_G(\mathfrak{L}) \oplus C_\bullet.
\end{equation*}
It is clear that the paracomplex $ C_\bullet $ is covariantly
contractible. \qed \\
This completes the proof of theorem \ref{Excision2}. 

\section{The exterior product}
\label{secext}

In this section we construct the exterior product for equivariant periodic cyclic homology. The exterior 
product is a generalization of the obvious composition product $ HP^G_*(A,B) \times HP^G_*(B,C) \rightarrow HP^G_*(A,C) $ 
discussed in section \ref{secHPdef} and an analogue of the exterior product 
in $ KK $-theory. Our discussion follows essentially the construction in the non-equivariant case 
given by Cuntz and Quillen \cite{CQ4}. \\
We need some preparations. First we define the tensor product of paracomplexes of covariant modules. 
Let $ C $ and $ D $ be paracomplexes of covariant modules and assume that $ C $ or $ D $ is a 
projective $ \mathcal{O}_G $-module. Then the tensor product $ C \boxtimes D $ of $ C $ and $ D $ is the paracomplex 
defined as follows. The space underlying $ C \boxtimes D $ is given by 
$$ (C \boxtimes D)_0 = C_0 \cotimes_{\mathcal{O}_G} D_0 \oplus C_1 \cotimes_{\mathcal{O}_G} D_1, \qquad 
(C \boxtimes D)_1 = C_1 \cotimes_{\mathcal{O}_G} D_0 \oplus C_0 \cotimes_{\mathcal{O}_G} D_1. 
$$
Observe that $ C \boxtimes D $ is complete due to our projectivity assumption. 
The group $ G $ acts diagonally and $ \mathcal{O}_G $ acts by multiplication. Using the fact that 
$ \mathcal{O}_G $ is commutative we see that $ C \boxtimes D $ becomes a covariant module in this way. \\
It remains to define the boundary operator in $ C \boxtimes D $. The usual formula for the differential 
in a tensor product of complexes is not appropriate since this formula does not yield a paracomplex. 
Instead we define the differential $ \partial $ in $ C \boxtimes D $ by 
$$
\partial_0 = 
\begin{pmatrix}
\partial \otimes \id & - \id \otimes \partial \\
\id \otimes \partial & \partial \otimes T
\end{pmatrix} 
\qquad 
\partial_1 = 
\begin{pmatrix}
\partial \otimes T & \id \otimes \partial \\
-\id \otimes \partial & \partial \otimes \id
\end{pmatrix}.
$$
It is straightforward to check that $ \partial^2 = \id - T $ in $ C \boxtimes D $. Hence the tensor product $ C \boxtimes D $ 
is again a paracomplex. \\
Now let $ I $ be a $ G $-invariant ideal in a pro-$ G $-algebra $ R $ and define a paracomplex 
$ \mathcal{H}^2_G(R,I) $ by 
$$ 
\mathcal{H}^2_G(R,I)^0 = \mathcal{O}_G \cotimes R/(\mathcal{O}_G \cotimes I^2 + b(\mathcal{O}_G \cotimes IdR))
$$ 
in degree zero and by  
$$
\mathcal{H}^2_G(R,I)^1 = \mathcal{O}_G \cotimes \Omega^1(R)/(b(\Omega^2_G(R)) + \mathcal{O}_G \cotimes I \Omega^1(R))
$$
in degree one where the boundary operators are induced from $ X_G(R) $. This paracomplex is the equivariant analogue of the 
corresponding quotient of the ordinary $ X $-complex considered in \cite{CQ2}. \\
Let $ A $ and $ B $ be pro-$ G $-algebras. In the same way as explained in \cite{CQ1} we see that the unital free product $ A^+ * B^+ $ of 
$ A^+ $ and $ B^+ $ can be written as  
$$
A^+ * B^+ = A^+ \cotimes B^+ \oplus \bigoplus_{j > 0} \Omega^j(A) \cotimes \Omega^j(B)
$$
where the multiplication is given by the Fedosov product
$$
(x_1 \otimes y_1) \circ (x_2 \otimes y_2) = x_1 x_2 \otimes y_1 y_2 - (-1)^{|x_1|} x_1 dx_2 \otimes dy_1 y_2. 
$$
An element $ a_0da_1 \cdots da_n \otimes b_0 db_1 \cdots db_n $ corresponds 
to $ a_0 b_0 [a_1, b_1]\cdots [a_n,b_n] $ in the free product under this identification 
where $ [x,y] = xy - yx $ denotes the ordinary commutator. \\
Consider the extension 
\begin{equation*}
  \xymatrix{
     I\;\; \ar@{>->}[r] & A^+ * B^+ \ar@{->>}[r]^\pi & A^+ \cotimes B^+ 
     }
\end{equation*}
of pro-$ G $-algebras where $ I $ is the kernel of the canonical homomorphism $ \pi: A^+ * B^+ \rightarrow A^+ \cotimes B^+ $. 
Using the description of the free product explained above we have 
$$
I^k = \bigoplus_{j \geq k} \Omega^j(A) \cotimes \Omega^j(B) 
$$
for the powers of the ideal $ I $. \\
Let us abbreviate $ R = A^+ * B^+ $ and define a covariant map $ \phi: X_G(A^+) \boxtimes X_G(B^+) \rightarrow \mathcal{H}^2_G(R,I) $ by 
\begin{align*}
&\phi(f(t) \otimes x \otimes y) = f(t) \otimes xy \\
&\phi(f(t) \otimes x_0dx_1 \otimes y_0dy_1) = f(t) \otimes x_0 (t^{-1} \cdot y_0)[x_1, t^{-1} \cdot y_1]  \\
&\phi(f(t) \otimes x \otimes y_0 dy_1) = f(t) \otimes xy_0 dy_1 \\
&\phi(f(t) \otimes x_0dx_1 \otimes y) = f(t) \otimes x_0dx_1y 
\end{align*}
where again $ [x,y] = xy - yx $ denotes the commutator.
\begin{prop} \label{tensorprop}
The map $ \phi: X_G(A^+) \boxtimes X_G(B^+) \rightarrow \mathcal{H}^2_G(R,I) $ defined above is a covariant 
isomorphism of paracomplexes for all pro-$ G $-algebras $ A $ and $ B $. 
\end{prop}
\proof According to the description of the free product using noncommutative differential forms we have 
an equivariant isomorphism 
$$
A^+ \cotimes B^+ \oplus \Omega^1(A) \cotimes \Omega^1(B) \cong R/I^2
$$
of $ A^+ \cotimes B^+ $-bimodules
This induces an isomorphism 
$$
X_G^0(A^+) \boxtimes X_G^1(B^+) \oplus X^1_G(A) \boxtimes X^1_G(B) \cong \mathcal{H}^2_G(R,I)^0
$$ 
and using lemma \ref{XA} we deduce
$$
X_G^0(A^+) \boxtimes X_G^1(B^+) \oplus X^1_G(A^+) \boxtimes X^1_G(B^+) \cong \mathcal{H}^2_G(R,I)^0. 
$$
After applying the covariant isomorphism $ T $ to $ X^1_G(B^+) $ this isomorphism can be 
identified with the map $ \phi $ in degree zero. \\
The inclusion maps $ A^+ \rightarrow R $ and $ B^+ \rightarrow R $ induce an equivariant $ R $-bimodule 
homomorphism
$$
R \cotimes_A \Omega^1(A) \cotimes_A R \oplus R \cotimes_B \Omega^1(B) \cotimes_B R \rightarrow \Omega^1(R).
$$
Tensoring with $ A^+ \cotimes B^+ $ over $ R $ on both sides we obtain a map 
$$
B^+ \cotimes \Omega^1(A) \cotimes B^+ \oplus A^+ \cotimes \Omega^1(B) \cotimes A^+ \rightarrow \Omega^1(R)/(I \Omega^1(R) + \Omega^1(R)I).
$$
Using the fact that $ R $ is unital we see as in \cite{CQ1} that this map determines an isomorphism 
$$
X^1_G(A) \boxtimes X^0_G(B^+) \oplus X^0_G(A^+) \boxtimes X^1_G(B) \cong \mathcal{H}^2_G(R,I)^1
$$
and by lemma \ref{XA} we obtain an isomorphism 
$$
X^1_G(A^+) \boxtimes X^0_G(B^+) \oplus X^0_G(A^+) \boxtimes X^1_G(B^+) \cong \mathcal{H}^2_G(R,I)^1
$$
which can be identified with the map $ \phi $ in degree one. \\
It remains to show that $ \phi $ is a chain map. To illustrate the occurence of the operator $ T $ we compute  
\begin{align*}
&\partial \phi(f(t) \otimes x_0dx_1 \otimes y_0dy_1) = f(t) \otimes d(x_0 (t^{-1} \cdot y_0)[x_1, t^{-1} \cdot y_1]) \\
&= f(t) \otimes x_0 (t^{-1} \cdot y_0) [dx_1, t^{-1} \cdot y_1] + 
f(t) \otimes x_0 (t^{-1} \cdot y_0) [x_1, d(t^{-1} \cdot y_1)] \\
&= f(t) \otimes x_0 dx_1 t^{-1} \cdot(y_0y_1) - f(t) \otimes x_0 dx_1 y_0 y_1 \\
&\qquad + f(t) \otimes x_0 x_1 t^{-1} \cdot (y_0 dy_1) - f(t) \otimes (t^{-1} \cdot x_1) x_0 t^{-1} \cdot (y_0 dy_1) \\
&= \phi \partial(f(t) \otimes x_0dx_1 \otimes y_0dy_1).
\end{align*}
The other cases are treated in a similar way. \qed 
\begin{lemma}\label{lemmafreeqf}
Let $ A $ and $ B $ be equivariantly quasifree pro-$ G $-algebras. Then the free product $ A^+ * B^+ $ 
is equivariantly quasifree.
\end{lemma}
\proof Let $ 0 \rightarrow K \rightarrow E \rightarrow Q \rightarrow 0 $ be a locally nilpotent extension 
of pro-$ G $-algebras and let $ f: A^+ * B^+ \rightarrow Q $ be an equivariant homomorphism. 
Since $ A^+ * B^+ $ is unital and $ \mathbb{C} $ is equivariantly quasifree we can lift the 
homomorphism $ \mathbb{C} \rightarrow Q $ induced by $ f $ to an equivariant homomorphism $ \mathbb{C} \rightarrow E $. 
We denote by $ e $ be the idempotent in $ E $ that corresponds to this lifting as well as 
its image in $ Q $. Then $ 0 \rightarrow e K e \rightarrow e E e \rightarrow Q \rightarrow 0 $ is again a locally nilpotent 
extension and the pro-$ G $-algebra $ e E e $ is unital. Since $ A $ and $ B $ are assumed to be quasifree there exist 
equivariant homomorphisms $ h_A: A \rightarrow e E e $ and $ h_B: B \rightarrow e E e $ lifting 
the maps $ A \rightarrow e Q e $ and $ B \rightarrow e Q e $ determined by $ f $. Extending $ h_A $ and $ h_B $ to the 
unitarizations and using the universal property 
of the free product we obtain a lifting $ h: A^+ * B^+ \rightarrow e E e $ for $ f $. Composing 
$ h $ with the evident map $ e E e \rightarrow E $ yields the claim. \qed \\
Next we discuss an analogue of the perturbation lemma \cite{Kassel}. 
Let $ C $ and $ D $ be paracomplexes. We shall assume that $ C $ and $ D $ are equipped with boundary operators $ b $ and $ B $ satisfying
$$ 
b^2 = 0 = B^2, \qquad Bb + bB = \id - T 
$$
such that $ \partial = B + b $. Consider the diagram 
\begin{equation*}
\xymatrix{
   D \ar@{->}[r]^i & C \ar@{->}[r]^p & D \\
 }
\end{equation*}
where $ i $ and $ p $ are chain maps with respect to the Hochschild operator $ b $ 
and assume that $ h: C \rightarrow C $ is an operator such that 
\begin{equation*}
pi = \id, \qquad ip = \id + (bh + hb).
\end{equation*}
Moreover we assume that $ p $ is a chain map with respect to $ B $. We will 
call such data a deformation retraction of $ C $ onto $ D $. A deformation 
retraction is called special if in addition the relations 
\begin{equation*}
hi = 0, \qquad ph = 0,\qquad h^2 = 0
\end{equation*}
hold. It is easy to see that any deformation retraction can be turned into a special deformation retraction. More 
precisely, if we define a new operator $ k $ by 
\begin{equation*}
k = (bh + hb)h(bh + hb) 
\end{equation*}
we get $ ki = pk = 0 $ since $ bh + hb = ip - \id $ and $ pi = \id $. Moreover one calculates
\begin{align*}
bk + kb &= (ip - \id)bh(ip - \id) + (ip - \id)hb(ip - \id) \\
&= (ip - \id)(ip - \id)(ip - \id) = ip - \id. 
\end{align*}
Hence $ k $ is again a deformation retraction. We define a map $ l $ 
by 
$$ 
l = -kbk
$$ 
and clearly get again $ li = pl = 0 $. Moreover we calculate 
\begin{align*}
bl + lb &= -bkbk - kbkb = -(ip - \id - kb)bk - kb(ip - \id - bk) = bk + kb = ip - \id 
\end{align*}
using that $ i $ and $ p $ are chain maps with respect to $ b $. The relation $ ip = \id + (bk + kb) $ 
implies $ k + bk^2 + k b k = 0 $ and $ k + k^2 b + kbk = 0 $. Combining these equations we obtain
$ bk^2 - k^2b = 0 $ and compute  
\begin{equation*}
l^2 = k b k^2 b k = k b^2 k^3 = 0. 
\end{equation*}
Hence we have constructed a special deformation retration.   
\begin{lemma}\label{perturblemma} Let $ C $ and $ D $ be paracomplexes and assume that $ l $ is a special deformation 
retraction of $ C $ onto $ D $. Then we have  
\begin{equation*}
[(lB)^ji, b] = -[(lB)^{j - 1}i, B]
\end{equation*}
and 
\begin{equation*}
[(lB)^j,b]l = B(lB)^{j - 1}l.
\end{equation*}
for all $ j > 0 $. 
\end{lemma}
\proof We use induction on $ j $. Consider the first expression. For $ j = 1 $ we 
have 
\begin{align*}
lBib - blBi &= lBbi + (lb + \id - ip)Bi \\
&= lBbi + lbBi + Bi - ipBi \\
&= Bi - iBpi = Bi - iB 
\end{align*}
since $ l(\id - T)i = li(\id - T) = 0 $. 
Assume that the claim is proved for $ j $ and compute 
\begin{align*}
[(lB)(lB)^j i,b] &= (lB)[(lB)^j i,b] + [lB,b](lB)^j i\\
&= -(lB)[(lB)^{j - 1}i,B] + (lBb - blB)(lB)^j i \\
&= -(lB)(lB)^{j - 1}i B - (lb + bl)B(lB)^j i \\
&= -(lB)^jiB + (\id - ip)B(lB)^j i \\
&= -(lB)^jiB + B(lB)^j i \\
&= -[(lB)^j i, B]
\end{align*}
using $ l^2 = 0 $. In order to prove the second formula we proceed in the same way. For 
$ j = 1 $ we have 
\begin{equation*}
lBbl - blBl = -lbBl - blBl = (\id - ip)Bl = Bl. 
\end{equation*}
Assume that the claim is proved for $ j $. Then we get 
\begin{align*}
[(lB)(lB)^j,b]l &= (lB)[(lB)^j,b]l + [lB,b](lB)^j l \\
&= lB(B(lB)^{j - 1}l) + (\id - ip) B(lB)^j l \\
&= B(lB)^j l. 
\end{align*}
This finishes the proof. \qed \\
Now let $ (\mathcal{T}A)^+ $ be the unitarized periodic tensor algebra of a pro-$ G $-algebra $ A $. 
According to theorem \ref{qf} there exists a resolution of 
$ (\mathcal{T}A)^+ $ by projective $ (\mathcal{T}A)^+ $-bimodules of length $ 1 $. 
If $ B $ is a second pro-$ G $-algebra we obtain a projective 
resolution of length $ 2 $ of the pro-$ G $-algebra $ C = (\mathcal{T}A)^+ \cotimes (\mathcal{T}B)^+ $ 
by tensoring the resolutions of $ (\mathcal{T}A)^+ $ and $ (\mathcal{T}B)^+ $. 
Using proposition \ref{propndim} we obtain an equivariant graded connection 
$ \nabla: \Omega^2(C) \rightarrow \Omega^3(C) $ for $ C $. According to proposition \ref{homotopyinv1} this yields a 
covariant homotopy equivalence between the Hochschild complexes of $ \theta \Omega_G(C) $ and $ \theta^2\Omega_G(C) $. \\
Let $ p: \theta \Omega_G(C) \rightarrow \theta^2 \Omega_G(C) $ be the natural projection,
$ i: \theta^2\Omega_G(C) \rightarrow \theta \Omega_G(C) $ be given by 
$ i = \id - [b,\nabla_G] $ and $ h = - \nabla_G: \theta \Omega_G(C) \rightarrow \theta \Omega_G(C) $. 
This defines a deformation retraction of $ \theta \Omega_G(C) $ onto $ \theta^2\Omega_G(C) $. 
Let $ l: \theta \Omega_G(C) \rightarrow \theta \Omega_G(C) $ be the special 
deformation retraction associated to $ h $ in the way described above. Since $ l $ increases the degree of a differential 
form by $ 1 $ the formula 
$$ 
K = \sum_{j = 0}^\infty (lB)^j 
$$ 
yields a well-defined operator $ K: \theta \Omega_G(C) \rightarrow \theta\Omega_G(C) $. 
We define in addition $ I = Ki $, $ H = Kl $ and $ P = p $. 
Then one has 
\begin{equation*}
PI = pKi = p \sum_{j = 0}^\infty (lB)^ji = pi = \id. 
\end{equation*}
The first relation of lemma \ref{perturblemma} yields $ [Ki,b] = -[Ki,B] $ 
and hence $ [I,B + b] = 0 $. Consequently $ I: \theta^2\Omega_G(C) \rightarrow \theta \Omega_G(C) $ is a chain 
map with respect to the total boundary $ B + b $. The second relation of lemma 
\ref{perturblemma} implies $ [K,b]l = BKl $ and from the definition of $ K $ 
we see $ K = \id + KlB $. This implies 
\begin{align*}
IP &= Kip = K + Kbl + Klb = K + BKl + bKl + Klb \\
&= \id + KlB + BKl + bKl + Klb = \id + [H,B + b].  
\end{align*}
Hence we have proved the following result.  
\begin{prop}\label{extprodAB}
Let $ A $ and $ B $ be pro-$ G $-algebras. Then the natural projection 
$ \theta \Omega_G((\mathcal{T}A)^+ \cotimes (\mathcal{T}B)^+) \rightarrow 
\theta^2 \Omega_G((\mathcal{T}A)^+ \cotimes (\mathcal{T}B)^+ ) $ is a covariant homotopy equivalence. 
\end{prop}
Now we are ready to prove the main theorem needed for the construction of the exterior product. 
\begin{theorem}\label{extprodnu}
Let $ A $ and $ B $ be pro-$ G $-algebras. Then there exists a natural covariant homotopy equivalence 
$$
X_G((\mathcal{T}A)^+ ) \boxtimes X_G((\mathcal{T}B)^+) \simeq X_G(\mathcal{T}(A^+ \cotimes B^+))
$$
of paracomplexes. 
\end{theorem}
\proof Let us write $ Q = (\mathcal{T}A)^+ \cotimes (\mathcal{T}B)^+ $ and consider the extension 
\begin{equation*}
  \xymatrix{
     I\;\; \ar@{>->}[r] & R \ar@{->>}[r]^\pi & Q
     }
\end{equation*}
where $ R = (\mathcal{T}A)^+ * (\mathcal{T}B)^+ $ is the unital free product of $ (\mathcal{T}A)^+ $ and $ (\mathcal{T}B)^+ $ and 
$ I $ is the kernel of the canonical homomorphism $ \pi: R \rightarrow Q $. By proposition \ref{tensorprop} we have a natural 
isomorphism 
$$
X_G((\mathcal{T}A)^+ ) \boxtimes X_G((\mathcal{T}B)^+) \cong \mathcal{H}^2_G(R,I)
$$
of paracomplexes. Define pro-$ G $-algebra $ \mathcal{R} $ and $ \mathcal{I} $ by taking 
the projective limit of the pro-$ G $-algebras $ R/I^n $ and $ I/I^n $, respectively. Then $ \mathcal{I} $ 
is locally nilpotent and we obtain an extension  
\begin{equation*}
  \xymatrix{
     \mathcal{I}\;\; \ar@{>->}[r] & \mathcal{R} \ar@{->>}[r]^\pi & Q
     }
\end{equation*}
of pro-$ G $-algebras. Since $ (\mathcal{T}A)^+ $ and $ (\mathcal{T}B)^+ $ are equivariantly quasifree the same holds true 
for $ R $ according to lemma \ref{lemmafreeqf}. It follows easily that $ \mathcal{R} $ is equivariantly quasifree as well. Hence we have in 
fact constructed a universal locally nilpotent extension of $ Q $. 
Due to proposition \ref{univext2} we deduce that $ \mathcal{T}Q $ and $ \mathcal{R} $ are equivariantly homotopy equivalent 
relative to $ Q $ and according to proposition \ref{homotopyinv} there exists a natural covariant homotopy equivalence 
$ X_G(\mathcal{R}) \simeq X_G(\mathcal{T}Q) $. It is easy to see that the chain maps between $ X_G(\mathcal{R}) $ and $  X_G(\mathcal{T}Q) $ 
implementing this homotopy equivalence induce chain maps between the quotients $ \mathcal{H}^2_G(\mathcal{R}, \mathcal{I}) $ 
and $ \mathcal{H}^2_G(\mathcal{T}Q, \mathcal{J}Q) $. Using the explicit formula written down after theorem \ref{homotopyinv}
we see that the corresponding chain homotopies also descend to operators on 
$ \mathcal{H}^2_G(\mathcal{R}, \mathcal{I}) $ and $ \mathcal{H}^2_G(\mathcal{T}Q, \mathcal{J}Q) $, respectively. 
Hence we obtain a natural covariant homotopy equivalence 
$$
\mathcal{H}^2_G(\mathcal{R}, \mathcal{I}) \simeq \mathcal{H}^2_G(\mathcal{T}Q, \mathcal{J}Q). 
$$
Next observe that there exists an obvious map $ X_G(R/I^n) \rightarrow \mathcal{H}^2_G(R,I) $ 
for $ n > 1 $. This implies that the projection $ R \rightarrow R/I^n $ induces an isomorphism 
$ \mathcal{H}^2_G(R,I) \rightarrow \mathcal{H}_G(R/I^n,I/I^n) $ for all $ n > 1 $. Hence we obtain a natural 
isomorphism
$$
\mathcal{H}^2_G(R,I) \cong \mathcal{H}^2_G(\mathcal{R}, \mathcal{I}).
$$
The definition of $ \mathcal{H}^2_G $ is made in such a way that the covariant homotopy equivalence  
$ X_G(\mathcal{T}Q) \simeq \theta \Omega_G(Q) $ obtained in theorem \ref{homotopyeq} induces a homotopy equivalence 
$ \mathcal{H}^2_G(\mathcal{T}Q, \mathcal{J}Q) \simeq \theta^2 \Omega_G(Q) $. 
We apply proposition \ref{extprodAB} to obtain 
$$
\theta^2 \Omega_G(Q) \simeq \theta \Omega_G(Q).
$$
Again by theorem \ref{homotopyeq} we have a natural homotopy equivalence $ \theta \Omega_G(Q) \simeq X_G(\mathcal{T}Q) $. 
Finally recall that tensor products of the form $ \mathcal{J}C \cotimes D $ with arbitrary 
pro-$ G $-algebras $ C $ and $ D $ are locally nilpotent by lemma \ref{tensorlemma}. 
Using this fact we obtain a natural covariant homotopy equivalence 
$$
X_G(\mathcal{T}Q) \simeq X_G(\mathcal{T}(A^+ \cotimes B^+))
$$
by applying the excision theorem \ref{Excision2} to the tensor products of the extensions 
$ 0 \rightarrow \mathcal{J}A \rightarrow (\mathcal{T}A)^+ 
\rightarrow A^+ \rightarrow 0 $ and $ 0 \rightarrow \mathcal{J}B \rightarrow (\mathcal{T}B)^+ \rightarrow B^+ \rightarrow 0 $. \\
Assembling these isomorphisms and homotopy equivalences yields the assertion. \qed 
\begin{cor} \label{exteriorprod} 
Let $ A $ and $ B $ be arbitrary pro-$ G $-algebras. Then there 
exists a natural covariant homotopy equivalence 
$$
X_G(\mathcal{T}A) \boxtimes X_G(\mathcal{T}B) \simeq X_G(\mathcal{T}(A \cotimes B))
$$
of paracomplexes. 
\end{cor}
\proof For every pro-$ G $-algebra $ D $ there exists a natural commutative diagram  
\begin{equation*} 
    \xymatrix{
     X_G(\mathcal{T}D)\, \ar@{->}[r] \ar@{->}[d]^{\id} & 
       X_G(\mathcal{T}(D^+))  \ar@{->}[r] \ar@{->}[d]^{\simeq}\ar@{<.}[l]&
           X_G(\mathcal{T}\mathbb{C}) \ar@{->}[d]^{\simeq} \\
  X_G(\mathcal{T}D)\, \ar@{->}[r] & 
         X_G((\mathcal{T}D)^+) \ar@{->}[r] \ar@{<.}[l]&
           X_G(\mathbb{C})   
     }
\end{equation*}
Using this we obtain the assertion from theorem \ref{extprodnu} by applying the excision 
theorem \ref{Excision2} to all possible tensor products of the extensions $ 0 \rightarrow A \rightarrow A^+ \rightarrow \mathbb{C} 
\rightarrow 0 $ and $ 0 \rightarrow B \rightarrow B^+ \rightarrow \mathbb{C} \rightarrow 0 $. \qed \\
Let $ A, B $ and $ D $ be pro-$ G $-algebras and define a map 
$$ 
\tau_D: HP^G_*(A,B) \rightarrow HP^G_*(A \cotimes D, B \cotimes D) 
$$
as follows. On the level of complexes we send a map 
$ \phi: X_G(\mathcal{T}(A \cotimes \mathcal{K}_G)) \rightarrow X_G(\mathcal{T}(B \cotimes \mathcal{K}_G)) $ 
to the map 
\begin{align*} 
\tau_D(\phi): X_G(\mathcal{T}(A &\cotimes D \cotimes \mathcal{K}_G)) \simeq 
X_G(\mathcal{T}(A \cotimes \mathcal{K}_G)) \boxtimes X_G(\mathcal{T}D) \\
&\xymatrix{\ar@{->}[r]^{\!\!\!\!\!\!\!\!\!\!\!\!\!\!\!\!\!\!\!\!\!\!\!\!\!\!\!\!\!\!\!\!\!\!\!\!\!\!\!\!\!\! \phi \cotimes \id} &  
X_G(\mathcal{T}(B \cotimes \mathcal{K}_G)) \boxtimes X_G(\mathcal{T}D)    
     } 
\simeq X_G(\mathcal{T}(B \cotimes D \cotimes \mathcal{K}_G)) 
\end{align*}
and consider the map induced in homology. Here we have used theorem \ref{exteriorprod} and 
suppressed the canonical isomorphisms corresponding to rearrangements of tensor products. \\
We can now proceed to define the exterior product. Let $ A_1, A_2, D, B_1, B_2 $ be pro-$ G $-algebras and 
let $ \phi \in HP^G_*(A_1, B_1 \cotimes D) $ and $ \psi \in HP^G_*(D \cotimes A_2, B_2) $ be two elements. 
After reordering the tensor factors we can thus use the ordinary composition product 
to compose $ \tau_{A_2}(\phi) \in HP^G_*(A_1 \cotimes A_2 \cotimes D, B_1 \cotimes A_2) $ 
and $ \tau_{B_1}(\psi) \in HP^G_*(B_1 \cotimes A_2, B_1 \cotimes B_2) $ and obtain 
$$ 
\phi \cotimes_D \psi = \tau_{A_2}(\phi) \cdot \tau_{B_1}(\psi)
$$ 
in $ HP^G_*(A_1 \cotimes A_2, B_1 \cotimes B_2) $. 
The following theorem summarizes some properties of the exterior product and 
is easily proved by inspecting the constructions.  
\begin{theorem} Let $ A_1, B_1, D, A_2, B_2 $ be pro-$ G $-algebras. The exterior product 
$$
HP^G_*(A_1, B_1 \cotimes D) \times HP^G_*(D \cotimes A_2, B_2) \rightarrow 
HP^G_*(A_1 \cotimes A_2, B_1 \cotimes B_2) 
$$
is bilinear, contravariantly functorial in $ A_1 $ and $ A_2 $ and covariantly functorial in 
$ B_1 $ and $ B_2 $. \\
The exterior product $ HP^G_*(A_1, \mathbb{C} \cotimes D) \times HP^G_*(D \cotimes \mathbb{C}, B_2) \rightarrow 
HP^G_*(A_1, B_2) $ can be identified with the composition product 
$ HP^G_*(A_1, D) \times HP^G_*(D, B_2) \rightarrow HP^G_*(A_1, B_2) $.
\end{theorem} 

\section{Compact Lie groups and the Cartan model}
\label{seclie}

After having studied the general homological properties of $ HP^G_* $ we shall 
now consider a more concrete situation. We will also show that our definition of 
equivariant cyclic homology generalizes previous constructions in the literature. \\
Let $ G $ be a compact group. Using proposition \ref{defcomp}, the fact that 
the trivial $ G $-algebra $ \mathbb{C} $ is quasifree, lemma \ref {XC} and theorem \ref{homotopyeq} we see that our definition of equivariant 
cyclic homology of a $ G $-algebra $ A $ reduces to 
$$
HP^G_*(A) = HP^G_*(\mathbb{C},A) = H_*(\SHom_G(\mathcal{O}_G[0], \theta \Omega_G(A)) = H_*(\varprojlim_n \theta^n \Omega_G(A)^G)
$$
in this case. Here $ \Omega_G(A)^G $ denotes the space of $ G $-invariant elements in $ \Omega_G(A) $. It is easy to check that 
$ T = \id $ on $ \Omega_G(A)^G $ which implies immediately that the invariant forms $ \Omega_G(A)^G $ are a mixed complex in a natural way. 
Moreover, $ HP^G_*(A) $ is just the cyclic homology of this mixed complex in the usual sense ~\cite{Loday}. 
Hence there are $ SBI $-sequences and other standard tools in order to compute these groups. In particular there 
is also a natural definition of equivariant Hochschild homology $ HH^G_*(A) $ and equivariant cyclic homology $ HC^G_* $ in this case. \\
Moreover we essentially reobtain the definition of equivariant cyclic homology for compact Lie groups 
as it has been introduced in the work of Brylinski ~\cite{Brylinski1}, ~\cite{Brylinski2}. The only difference is 
that Brylinski works with topological vector spaces whereas we use bornological vector spaces. \\
Let us now consider the important special case of a compact Lie group acting smoothly on a compact manifold $ M $. 
We remark that in this case there is no difference between the topological and the bornological approach. 
It turns out that the equivariant periodic cyclic homology of $ C^\infty(M) $ is closely related 
to the equivariant $ K $-theory of $ M $. The following theorem was obtained by Brylinski ~\cite{Brylinski1} and independently by 
Block ~\cite{Block}. 
\begin{theorem}\label{Bryl1}
Let $ G $ be a compact Lie group acting smoothly on a smooth compact manifold 
$ M $. There exists an equivariant Chern character 
\begin{equation*}
ch_G: K^*_G(M) \rightarrow HP^G_*(C^\infty(M)) 
\end{equation*}
which induces an isomorphism 
\begin{equation*}
HP_*^G(C^\infty(M)) \cong \mathcal{R}(G) \otimes_{R(G)} K^*_G(M) 
\end{equation*}
where $ R(G) $ is the representation ring of $ G $ 
and $ \mathcal{R}(G) = C^\infty(G)^G $ is the algebra of smooth conjugation 
invariant functions on $ G $. 
\end{theorem}
Here of course $ \mathcal{R}(G) $ is viewed as an $ R(G) $-module using the character map. \\
Block and Getzler have obtained a description of $ HP^G_*(C^\infty(M)) $ in terms of equivariant differential 
forms ~\cite{BG}. More precisely, there exists a $ G $-equivariant sheaf $ \Omega(M,G) $ 
over the group $ G $ itself viewed as a $ G $-space with the adjoint action. 
The stalk $ \Omega(M,G)_s $ at a group element $ s \in G $ 
is given by germs of $ G_s $-equivariant smooth maps from $ \mathfrak{g}^s $ to 
$ \mathcal{A}(M^s) $. Here $ M^s = \{x \in M|\, s \cdot x = x \} $ is the fixed point 
set of $ s $, $ G^s $ is the centralizer of $ s $ in $ G $ and $ \mathfrak{g}^s $ 
is the Lie algebra of $ G_s $. In particular the stalk $ \Omega(M,G)_e $
at the identity element $ e $ is given by 
\begin{equation*}
\Omega(M,G)_e = C^\infty_0(\mathfrak{g}, \mathcal{A}(M))^G
\end{equation*}
where $ C^\infty_0 $ is the notation for smooth germs at $ 0 $. 
Hence $ \Omega(M,G)_e $ can be viewed as a certain completion of the classical Cartan model $ \mathcal{A}_G(M) $. 
The global sections $ \Gamma(G, \Omega(M,G)) $ of the sheaf $ \Omega(M,G) $
are called global equivariant differential forms and will be denoted 
by $ \mathcal{A}(M,G) $. 
There exists a natural differential on $ \mathcal{A}(M,G) $ extending the Cartan differential. 
Block and Getzler establish an equivariant version of the Hochschild-Kostant-Rosenberg theorem and deduce the following result. 
\begin{theorem} \label{Bryl2} Let $ G $ be a compact Lie group acting smoothly on 
a smooth compact manifold $ M $. Then there is a natural isomorphism
\begin{equation*}
HP^G_*(C^\infty(M)) \cong H^*(\mathcal{A}(M,G)). 
\end{equation*}
\end{theorem}
This theorem shows that equivariant cyclic homology can be viewed as a "delocalized" noncommutative version of the Cartan model.  
Theorem \ref{Bryl2} also shows that the language of equivariant sheaves is necessary to describe equivariant cyclic homology 
appropriately. Combining theorem \ref{Bryl1} and theorem \ref{Bryl2} one obtains the following result.
\begin{theorem}
Let $ G $ be a compact Lie group acting smoothly on a smooth compact manifold $ M $. 
Then there exists a natural isomorphism 
\begin{equation*}
\mathcal{R}(G) \otimes_{R(G)} K^*_G(M) \cong H^*(\mathcal{A}(M,G)). 
\end{equation*}
\end{theorem}
Hence, up to an ``extension of scalars'', the equivariant $ K $-theory of manifolds 
can be described using global equivariant differential forms. \\
We emphasize that we do \emph{not} define $ HH_*^G $ 
and $ HC^G_* $ for non-compact groups. It seems to be unclear how a reasonable 
definition of such theories should look like. Clearly one would like to 
have $ SBI $-sequences and a relation to equivariant periodic cyclic 
homology $ HP^G_* $ similar to the one for compact groups. \\ 
Finally, we mention that for finite groups our definition of equivariant periodic cyclic cohomology is 
compatible with the constructions in ~\cite{KKL1}. 

\section{The Green-Julg theorem}
\label{secGJ}

The Green-Julg theorem ~\cite{Green}, \cite{Julg} asserts that 
for a compact group $ G $ the equivariant $ K $-theory $ K^G_*(A) $ of a 
$ G $-$ C^* $-algebra $ A $ is naturally isomorphic to the ordinary 
$ K $-theory $ K_*(A \rtimes G) $ of the crossed product $ C^* $-algebra 
$ A \rtimes G $. \\
In this section we prove an analogue of the Green-Julg theorem 
in cyclic homology. In its original form this result is due to Brylinski 
~\cite{Brylinski1}, ~\cite{Brylinski2} who studied smooth actions of compact 
Lie groups. Independently this version of the Green-Julg theorem was obtained by 
Block ~\cite{Block}. We follow the work of Bues ~\cite{Bues1}, ~\cite{Bues2} and prove a 
variant of this theorem for pro-algebras and arbitrary compact groups. 
Some ingredients in the proof show up in a similar way in the computation 
of the cyclic cohomology of crossed products in general \cite{GJ1}, \cite{Nistor1}, \cite{Nistor2}. \\
Our Green-Julg theorem involves crossed products of pro-$ G $-algebras. We remark that the construction 
of crossed products for $ G $-algebras can immediately be extended to 
pro-$ G $-algebras.
\begin{theorem}\label{Green-Julg} Let $ G $ be a compact group and let $ A $ be a
pro-$ G $-algebra. Then there is a natural isomorphism
\begin{equation*}
HP^G_*(\mathbb{C},A) \cong HP_*(A \rtimes G).
\end{equation*}
\end{theorem}
For the proof of theorem \ref{Green-Julg} we need some preparations. Throughout this section we 
assume that the Haar measure on the compact group $ G $ is normalized and we denote 
by $ H = \D(G) $ the smooth group algebra of $ G $. There are $ H $-bimodule splittings 
$ \sigma_n: H \rightarrow H^{\cotimes n} $ for 
the iterated multiplication given by 
$$
\sigma_n(f)(s_1,\dots, s_n) = f(s_1 \cdots s_n). 
$$
Using this fact it is not hard to show that $ H $ is projective as an $ H $-bimodule and quasifree as a bornological 
algebra. 
\begin{prop}\label{GJ1} Let $ G $ be a compact group and let $ R $ be a unital quasifree 
pro-$ G $-algebra. Then the pro-algebra $ R \rtimes G $ is quasifree. 
\end{prop}
\proof We have to construct a splitting homomorphism $ w: R \rtimes G \rightarrow 
\mathcal{T}(R \rtimes G) $ for the canonical projection. Since $ R $ is assumed to be quasifree there 
exists an equivariant lifting homomorphism $ u: R \rightarrow \mathcal{T}R $ for 
the projection $ \tau_R: \mathcal{T}R \rightarrow R $. After taking crossed 
products we obtain a homomorphism $ u \rtimes G: R \rtimes G \rightarrow \mathcal{T}R \rtimes G $ lifting the 
homomorphism $ \tau_R \rtimes G $. Consider the equivariant linear map $ h: R \rtimes G \rightarrow \mathcal{T}R \rtimes G $ 
obtained by tensoring $ \sigma_R $ with the identity on $ H $. It is straightforward to check that 
$ h $ is a lonilcur. Hence according to proposition \ref{PeriodicTensorAlg} we obtain a 
homomorphism $ [[h]]: \mathcal{T}(R \rtimes G) \rightarrow \mathcal{T}R \rtimes G $ such that 
$ [[h]] \sigma_{R \rtimes G} = h $. We obtain a linear splitting 
$ \sigma: \mathcal{T}R \rtimes G \rightarrow \mathcal{T}(R \rtimes G) $ for $ [[h]] $ by setting 
\begin{align*}
\sigma(x_0 &dx_1 \cdots dx_{2n} \rtimes f)(r_0, \dots, r_{2n}) = \\
&\sigma_{2n + 1}(f)(r_0,\dots, r_{2n})\, x_0 d(r_0^{-1} \cdot x_1) d((r_0 r_1)^{-1} \cdot x_2) \cdots d((r_0 \cdots r_{2n - 1})^{-1} \cdot x_{2n}).  
\end{align*}
This implies that the homomorphism $ [[h]] $ fits into an extension
\begin{equation*} 
    \xymatrix{\mathcal{J} \;\; \ar@{>->}[r] &
         \mathcal{T}(R \rtimes G) \ar@{->>}[r] \ar@{<.}[l]&
           \mathcal{T}R \rtimes G   
     }
\end{equation*}
where the kernel $ \mathcal{J} $ of $ [[h]] $ is locally nilpotent. Hence this extension is 
a universal locally nilpotent extension of $ \mathcal{T} R \rtimes G $. 
Consider the homomorphism $ \iota: H \rightarrow R \rtimes G $ given by $ \iota(f) = 1_R \rtimes f $. 
We compose $ \iota $ with $ u \rtimes G $ to obtain a homomorphism 
$ (u \rtimes G) \iota: H \rightarrow \mathcal{T}R \rtimes G $. 
Since $ G $ is compact the smooth group algebra $ H $ is quasifree. 
By theorem \ref{qf} we obtain a homomorphism $ \phi: H \rightarrow \mathcal{T}(R \rtimes G) $ such that 
$ [[h]] \phi = (u \rtimes G) \iota $. 
In this way the algebra $ \mathcal{T}(R \rtimes G) $ becomes an $ H $-bimodule. 
We shall now construct another linear lifting $ \lambda $ of the homomorphism 
$ \tau_{R \rtimes G}: \mathcal{T}(R \rtimes G) \rightarrow R \rtimes G $. 
Consider first the map $ l: R \rtimes G \rightarrow H \cotimes (R \rtimes G) \cotimes H $  
given by  
$$
l(x \rtimes f)(r,s,t) = \sigma_3(f)(r,s,t)\, r^{-1} \cdot x. 
$$
By construction $ l $ is an $ H $-bimodule map splitting the canonical multiplication 
map $ H \cotimes (R \rtimes G) \cotimes H \rightarrow R \rtimes G $. If we compose 
$ l $ with $ \phi \cotimes \sigma_{R \rtimes G} \cotimes \phi $ and apply multiplication 
in $ \mathcal{T}(R \rtimes G) $ we obtain an $ H $-bimodule map $ \lambda: R\rtimes G \rightarrow \mathcal{T}(R \rtimes G) $. 
One computes $ \tau_{R\rtimes G} \lambda = \id $ which implies in particular that  
$ \lambda $ is a lonilcur. By proposition \ref{PeriodicTensorAlg} we obtain a 
homomorphism $ [[\lambda]]: \mathcal{T}(R \rtimes G) \rightarrow \mathcal{T}(R \rtimes G) $ 
such that $ [[\lambda]] \sigma_{R \rtimes G} = \lambda $. 
Since $ \lambda $ is an $ H $-bimodule map it follows that $ [[\lambda]] $ descends to 
a homomorphism $ v: \mathcal{T}R \rtimes G \rightarrow \mathcal{T}(R \rtimes G) $  
satisfying $ v [[h]] = [[\lambda]] $. We compute 
$$
(\tau_R \rtimes G) [[h]] \sigma_{R \rtimes G} = (\tau_R \rtimes G) h = \id = 
\tau_{R \rtimes G} \lambda = \tau_{R \rtimes  G} [[\lambda]] \sigma_{R \rtimes G} 
= \tau_{R \rtimes G} v [[h]] \sigma_{R \rtimes G}
$$ 
and again by proposition \ref{PeriodicTensorAlg} we deduce $ (\tau_R \rtimes G) [[h]] = \tau_{R \rtimes G} v [[h]] $. 
Composition with the splitting $ \sigma: \mathcal{T}R \rtimes G \rightarrow \mathcal{T}(R \rtimes G) $ from above 
yields $ \tau_R \rtimes G  = \tau_{R \rtimes G} v $. Now we set $ w = v (u \rtimes G) $ and compute 
$$
\tau_{R \rtimes G} w = (\tau_R \rtimes G) (u \rtimes G) = \id. 
$$
Hence $ w $ is a splitting homomorphism for $ \tau_{R \rtimes G} $. \qed \\
Let us assume that $ R $ is a unital pro-$ G $-algebra and write $ B = R \rtimes G $. 
Since $ R $ is unital there exists a natural homomorphism $ H \rightarrow B $ and we always view $ B $ 
as an $ H $-bimodule in this way. Our next goal is to 
define a relative version of the $ X $-complex of $ B $ which can be compared to the equivariant 
$ X $-complex of $ R $. \\
Consider the linear map $ \lambda_0: B \rightarrow B $ defined by 
$$
\lambda_0(f)(t) = \int_G s \cdot f(s^{-1}t s) ds. 
$$
This map vanishes on the the space of commutators $ [B,H] $ and defines a linear splitting 
for the extension 
\begin{equation} \label{GJeq2}
    \xymatrix{
     [B,H] \;\; \ar@{>->}[r] &
         B \ar@{->>}[r] \ar@{<.}[l]&
           B/[B,H].
     }
\end{equation}
If we define $ K^0 = [B,H] $ and $ X^0(B)_H = B/[B,H] $ we can rewrite this as 
\begin{equation} \label{GJeq2a}
    \xymatrix{
     K^0 \;\; \ar@{>->}[r] &
         X^0(B) \ar@{->>}[r] \ar@{<.}[l]&
           X^0(B)_H.   
     }
\end{equation}
The space $ X^0(B)_H $ is the even part of the relative $ X $-complex. \\
Now consider the extension 
$$
   \xymatrix{
     \Omega^1(H)\;\; \ar@{>->}[r] &
         H^+ \cotimes H^+ \ar@{->>}[r] &
           H^+   
     }
$$
of $ H $-bimodules. This extension has a left $ H $-linear splitting, hence 
tensoring from the left with $ H $ over itself we obtain an extension  
\begin{equation} \label{GJeq1}
    \xymatrix{
     H \cotimes_H \Omega^1(H)\;\; \ar@{>->}[r] &
         H \cotimes H^+ \ar@{->>}[r] \ar@{<.}[l]&
           H   
     }
\end{equation}
of $ H $-bimodules. Remark that the map $ \sigma_2: H \rightarrow H \cotimes H $ from above 
yields an $ H $-bimodule splitting for extension (\ref{GJeq1}).  
We tensor extension (\ref{GJeq1}) over $ H $ with $ B $ on the left and 
with $ B^+ $ on the right to obtain the split extension 
\begin{equation} \label{GJeq3}
    \xymatrix{
     B \cotimes_H \Omega^1(H) \cotimes_H B^+ \;\; \ar@{>->}[r] &
         B \cotimes B^+ \ar@{->>}[r] \ar@{<.}[l]&
           B \cotimes_H B^+   
     }
\end{equation}
of $ B $-bimodules. Since $ R $ is unital we have a left $ B $-linear splitting 
$ \lambda_B: B \rightarrow B \cotimes B $ of the multiplication defined by 
$ \lambda_B(f)(s,t) = f(st) \cotimes 1_R $ 
where we identify $ B \cotimes B \cong R \cotimes R \cotimes H \cotimes H $ with a flip of the tensor 
factors. This yields split extensions of $ B $-bimodules
\begin{equation} \label{GJeq4}
    \xymatrix{
     B \cotimes_B \Omega^1(B) \;\; \ar@{>->}[r] &
         B \cotimes B^+ \ar@{->>}[r] \ar@{<.}[l]&
           B    
     }
\end{equation}
and 
\begin{equation} \label{GJeq5}
 \xymatrix{
     B \cotimes_B \Omega^1(B)_H \;\; \ar@{>->}[r] &
         B \cotimes_H B^+ \ar@{->>}[r] \ar@{<.}[l]&
           B   
     }
\end{equation}
where $ \Omega^1(B)_H $ is the kernel of the multiplication map $  B^+ \cotimes_H B^+ \rightarrow B^+ $. 
Assembling the extensions (\ref{GJeq3}), (\ref{GJeq4}) and (\ref{GJeq5}) we obtain a commutative diagram 
\begin{equation} \label{GJeq6}
    \xymatrix{
     B \cotimes_H \Omega^1(H) \cotimes_H B^+ \ar@{->}[r] \ar@{->}[d]^{\id} & 
       B \cotimes_B \Omega^1(B)  \ar@{->}[r] \ar@{->}[d]\ar@{<.}[l]&
           B \cotimes_B \Omega^1(B)_H  \ar@{->}[d] \\
 B \cotimes_H \Omega^1(H) \cotimes_H B^+ \ar@{->}[r] \ar@{->}[d] & 
         B \cotimes B^+ \ar@{->}[r] \ar@{->}[d]\ar@{<.}[l]&
           B \cotimes_H B^+ \ar@{->}[d]  \\
0 \;\ar@{->}[r]  & 
         B \ar@{->}[r]^{\id}\ar@{<.}[l]& B 
     }
\end{equation}
of $ B $-bimodules with split exact rows and columns. Observe that there are 
natural $ B $-bimodule maps $ B \cotimes_B \Omega^1(B) \rightarrow \Omega^1(B) $ 
and $ B \cotimes_B \Omega^1(B)_H \rightarrow \Omega^1(B)_H $. If we set 
$ X^1(B)_H = \Omega^1(B)_H/[-,B] $ we obtain a commutative diagram of pro-vector spaces
\begin{equation} \label{GJeq7}
    \xymatrix{
     (B \cotimes_B \Omega^1(B))/[-,B] \ar@{->}[d] \ar@{->>}[r] &
           (B \cotimes_B \Omega^1(B)_H)/[-,B]  \ar@{->}[d] \\
         X^1(B) \ar@{->}[r] & X^1(B)_H     
     }
\end{equation}
by taking commutator quotients with respect to $ B $ where the upper horizontal arrow has a linear splitting 
according to diagram (\ref{GJeq6}). We want to show that 
the vertical arrows in diagram (\ref{GJeq7})  are isomorphisms. Let $ \tau: B \cotimes B \rightarrow B \cotimes B $ 
be the flip of the tensor factors. Moreover let $ j: B \cotimes B \rightarrow (B \cotimes_B \Omega^1(B))/[-,B] $ be 
the map given by $ j(x_0 \otimes x_1) = x_0 \otimes dx_1 $. We define a linear map 
$ \rho: \Omega^1(B) \rightarrow B \cotimes_B \Omega^1(B)/[-,B] $ 
by setting 
$$ 
\rho(dx_1) = j \lambda_B(x_1) + j \tau \lambda_B(x_1), \qquad \rho(x_0 dx_1) = x_0 \otimes dx_1. 
$$
Using the Leibniz rule and the fact that $ \lambda_B $ is left $ B $-linear it is not hard to show that $ \rho $ 
descends to a map $ \rho: X^1(B) \rightarrow  B \cotimes_B \Omega^1(B)/[-,B] $. 
Once this is established it is easy to see that this map 
provides an inverse to the canonical map $ B \cotimes_B \Omega^1(B)/[-,B] \rightarrow X^1(B) $. A similar 
argument shows that the map $ B \cotimes_B \Omega^1(B)_H/[-,B] \rightarrow X^1(B)_H $ is an isomorphism. \\
If we define $ K^1 = (B \cotimes_H \Omega^1(H) \cotimes_H B^+)/[-,B] $ we now obtain an extension 
\begin{equation} \label{GJeq8}
 \xymatrix{
     K^1 \;\; \ar@{>->}[r] &
         X^1(B) \ar@{->>}[r] \ar@{<.}[l]&
           X^1(B)_H   
     }
\end{equation}
of pro-vector spaces using the first row in diagram (\ref{GJeq6}). \\
The differentials in the $ X $-complex $ X(B) $ descend to differentials in $ X(B)_H $. 
Hence diagrams \ref{GJeq2a} and \ref{GJeq8} yield an extension  
\begin{equation} \label{GJeq9}
 \xymatrix{
     K \;\; \ar@{>->}[r] &
         X(B) \ar@{->>}[r] \ar@{<.}[l]&
           X(B)_H   
     }
\end{equation}
of complexes with linear splitting. The complex $ X(B)_H $ will be called the relative 
$ X $-complex of $ B $ with respect to $ H $.  
\begin{prop} \label{SepHomEq}
The canonical chain map $ X(B) \rightarrow X(B)_H $ is a homotopy equivalence.
\end{prop}
\proof It suffices to show that 
the complex $ K $ is contractible. Consider the map $ \alpha: [B,H] \rightarrow (B \cotimes B^+)/[-,B] $ 
given by $ \alpha(x) = x \otimes 1 $. Since composition of $ \alpha $ with the natural map 
$ (B \cotimes B^+)/[-,B] \rightarrow (B \cotimes_H B^+)/[-,B] $ is zero we can view $ \alpha $  as a map from 
$ K^0 $ to $ K^1 $. It is straightforward to check that $ \alpha $ is inverse to the boundary 
$ b: K^1 \rightarrow K^0 $. This yields the claim. \qed \\
If $ R $ is a pro-$ G $-algebra we denote by $ X_G(R)^G $ 
the invariant part of the equivariant $ X $-complex of $ R $. Note that $ X_G(R)^G $ is in fact a pro-supercomplex. 
\begin{prop}\label{GJ3} Let $ G $ be a compact group and let $ R $ be a unital pro-$ G $-algebra. There is a natural isomorphism 
\begin{equation*}
X_G(R)^G \cong X(R \rtimes G)_H
\end{equation*}
of pro-supercomplexes where $ X(R \rtimes G)_H $ denotes the relative $ X $-complex. 
\end{prop}
\proof Since $ G $ is compact we can identify $ X_G(R)^G $ with the 
$ G $-coinvariants of $ X_G(R) $ by averaging over $ G $. We will denote 
the space of $ G $-coinvariants of $ X_G(R) $ by $ X_G(R)_G $. \\
In the sequel we identify elements of $ \mathcal{O}_G $ with elements in the group algebra $ \D(G) $ in the 
evident way.  
The action of $ s \in G $ on $ f \in \mathcal{O}_G $ corresponds to the adjoint action of $ s $ on
$ f $ in the group algebra $ \D(G) $. \\
We define a map $ \alpha: X_G(R)_G \rightarrow X(R \rtimes G)_H $ by 
\begin{align*}
&\alpha_0(f \otimes x)(s) = f(s) x \\
&\alpha_1(f \otimes xdy)(s,t) = f(st) x d(s^{-1} \cdot y) \\
&\alpha_1(f \otimes dy)(t) = f(t) dy
\end{align*}
where we view $ \alpha_1(f \otimes xdy) \in (R \rtimes G) \cotimes (R \rtimes G) $ as a function on $ G \times G $ with values in $ R \times R $. 
Moreover we define a map $ \beta: X(R \rtimes G)_H \rightarrow X_G(R)_G $ by 
\begin{align*}
&\beta_0(x \rtimes f) = f(r) x \\
&\beta_1((x \rtimes f) d(y \rtimes g))(r) = f(r) g(r) x d(r \cdot y) \\
&\beta_1(d(y \rtimes g))(r) = g(r) dy.
\end{align*}
Some straightforward computations show that these maps are well-defined and it is easy to see that $ \alpha $ and $ \beta $ are 
inverse to each other. We only show that $ \alpha $ is a chain map. One computes
\begin{equation*}
(d\alpha_0)(f \otimes x)(s) = f(s) dx = \alpha_1(f \otimes dx)(s) = 
(\alpha_1 d)(f \otimes x)(s)
\end{equation*}
and 
\begin{align*}
&(b \alpha_1)(f \otimes xdy)(t) = f(t) xy - \int_G f(r^{-1} t r) (t^{-1} r \cdot y) (r \cdot x) dr \\
&= f(t) xy - f(t) (t^{-1} \cdot y) x = (\alpha_0 b)(f \otimes xdy)(t). 
\end{align*}
This finishes the proof of proposition \ref{GJ3}. \qed \\
Now we come back to the proof of theorem \ref{Green-Julg}. 
Using the long exact sequences obtained in theorem \ref{Excision} both for 
equivariant cyclic homology and ordinary cyclic homology it suffices to prove the 
assertion for an augmented pro-$ G $-algebra of the form $ A^+ $. \\
On the one hand we have to compute the equivariant periodic cyclic homology of 
$ A^+ $. Due to proposition \ref{Unitalqf} we can use the universal locally nilpotent 
extension 
\begin{equation} \label{GJeq11}
 \xymatrix{
     \mathcal{J}A \;\; \ar@{>->}[r] &
         (\mathcal{T}A)^+ \ar@{->>}[r]^{\tau_A^+} &
           A^+   
     }
\end{equation}
to do this. Since the group $ G $ is compact and the $ G $-algebra 
$ \mathbb{C} $ is quasifree the equivariant periodic cyclic homology of $ A $ is 
consequently the homology of 
\begin{align*}
\SHom_G&(X_G(\mathbb{C}),X_G((\mathcal{T}A)^+)
= \SHom_G(\mathcal{O}_G[0],X_G((\mathcal{T}A)^+) \\
&= \Hom_G(\mathbb{C}[0],X_G((\mathcal{T}A)^+) = X_G((\mathcal{T}A)^+)^G. 
\end{align*}
On the other hand we have to calculate the cyclic homology of 
the crossed product $ A^+ \rtimes G $. Taking crossed products in extension (\ref{GJeq11}) we obtain 
an extension 
\begin{equation} \label{GJext}
 \xymatrix{
     \mathcal{J}A \rtimes G\;\; \ar@{>->}[r] &
         (\mathcal{T}A)^+ \rtimes G\ar@{->>}[r] &
           A^+ \rtimes G  
     }
\end{equation}
of pro-algebras. It is easy to check that the pro-$ G $-algebra $ \mathcal{J}A \rtimes G $ 
is locally nilpotent. Proposition \ref{GJ1} shows that $ (\mathcal{T}A)^+ \rtimes G $  
is quasifree and hence (\ref{GJext}) is in fact a universal locally nilpotent extension of $ A^+ \rtimes G $. 
This means that $ HP_*(A^+ \rtimes G) $ can be computed 
using $ X((\mathcal{T}A)^+ \rtimes G) $. 
Consider the relative $ X $-complex $ X((\mathcal{T}A)^+ \rtimes G)_H $ described above. Due to 
proposition \ref{SepHomEq} the pro-supercomplexes 
$ X((\mathcal{T}A)^+ \rtimes G) $ and 
$ X((\mathcal{T}A)^+ \rtimes G)_H $ are homotopy equivalent. 
From proposition \ref{GJ3} we obtain a natural isomorphism 
\begin{equation*}
X((\mathcal{T}A)^+ \rtimes G)_H \cong X_G((\mathcal{T}A)^+)^G. 
\end{equation*}
Hence we see that both theories agree. Since all constructions are natural in $ A $ this finishes the 
proof of theorem \ref{Green-Julg}. 

\section{The dual Green-Julg theorem}
\label{secDGJ}

In this section we study equivariant periodic cyclic cohomology in the case of  discrete groups. The main
result is the following dual version of the Green-Julg theorem \ref{Green-Julg}.
\begin{theorem} \label{DGJ} Let $ G $ be a discrete group and let 
$ A $ be a pro-$ G $-algebra. Then there is a natural isomorphism
\begin{equation*}
HP^G_*(A,\mathbb{C}) \cong HP^*(A \rtimes G).
\end{equation*}
\end{theorem}
This theorem yields in particular a description of $ HP^G_*(\mathbb{C},\mathbb{C}) $.
By the work of Burghelea ~\cite{Burghelea} it follows that the group cohomology 
of $ G $ with complex coefficients constitutes a direct factor of $ HP^G_*(\mathbb{C},\mathbb{C}) $.  
We remark that the isomorphism in theorem \ref{DGJ} is compatible with 
natural decompositions of $ HP^G_*(A,\mathbb{C}) $ and $ HP^*(A \rtimes G) $
over the conjugacy classes of $ G $. \\
The proof of theorem \ref{DGJ} is divided into two
parts. In the first part we obtain a simpler description of 
$ HP^G_*(A, B) $ for arbitrary pro-$ G $-algebras $ A $ and $ B $. 
For this we do not have to assume that $ G $ is discrete. \\
Let $ G $ be any locally compact group and let $ B $ be a pro-$ G $-algebra. Consider the map 
$ \Tr: \Omega_G(B \cotimes \mathcal{K}_G) \rightarrow \Omega_G(B) $ given on $ n $-forms by
\begin{align*}
\Tr(f(s) &\otimes (x_0 \otimes  k_0) d(x_1 \otimes k_1) \cdots d(x_n \otimes k_n)) \\
&= f(s) \otimes x_0 dx_1 \cdots dx_n \int k_0(r_0,r_1) k_1(r_1,r_2) \cdots k_n(r_n, sr_0) dr_0 \cdots dr_n 
\end{align*}
and
\begin{align*}
\Tr(f(s) &\otimes d(x_1 \otimes k_1) \cdots d(x_n \otimes k_n)) \\
&= f(s) \otimes dx_1 \cdots dx_n \int k_1(r_1,r_2) \cdots k_n(r_n, sr_1) dr_1 \cdots dr_n. 
\end{align*}
One checks that $ \Tr $ is a covariant map and that it commutes with the Hochschild boundary $ b $. By definition it 
commutes with the operator $ d $ and it follows that $ \Tr $ is a map of paramixed complexes. We remark that $ \Tr $ is 
closely related to the trace map that occured in the proof of the stability theorem \ref{StabLemma}. 
\begin{prop} \label{KG-Hochschild} Let $ G $ be a locally compact group and let $ B $ be a unital pro-$ G $-algebra. The map  
$ \Tr: \Omega_G(B \cotimes \mathcal{K}_G) \rightarrow \Omega_G(B) $ is a linear 
homotopy equivalence with respect to the equivariant Hochschild boundary. 
\end{prop}
\proof As in ordinary Hochschild homology we may view the equivariant Hochschild complex $ \Omega_G(C) $ of any pro-$ G $-algebra $ C $ as the total 
complex of a double complex with two columns. This is induced by the decomposition 
$ \Omega^n_G(C) = \mathcal{O}_G \hat{\otimes} C^{\hat{\otimes} n + 1} \oplus \mathcal{O}_G \hat{\otimes} C^{\hat{\otimes} n} $. One checks 
easily that the second columns of this double complex is simply the Bar-complex of $ C $ tensored with $ \mathcal{O}_G $ whereas the first column is 
equipped with the equivariant Hochschild boundary. \\
We apply this description to the $ G $-algebras $ B \cotimes \mathcal{K}_G $ and $ B $. In order to prove the proposition it suffices to show that the 
columns of the corresponding bicomplexes are linearly homotopy equivalent.\\
Choose a smooth function $ \chi \in \D(G) $ such that 
$$
\int_G \chi^2(t) dt = 1
$$
and consider the bounded linear map $ \sigma: \mathcal{K}_G \rightarrow \mathcal{K}_G \cotimes \mathcal{K}_G $ 
defined by 
$$
\sigma(k)(r_1,t_1,r_2,t_2) = k(r_1, t_2) \chi(t_1) \chi(r_2). 
$$
It is easy to check that $ \sigma $ is a $ \mathcal{K}_G $-bimodule map that splits the multiplication $ \mathcal{K}_G \cotimes \mathcal{K}_G 
\rightarrow \mathcal{K}_G $. We remark that the map $ \sigma $ can be used to show that $ \mathcal{K}_G $ is a quasifree algebra. However, we emphasize 
that this algebra is usually far from being \emph{equivariantly} quasifree. \\
Let us consider the second column in the bicomplex associated to $ B \cotimes \mathcal{K}_G $. We define a contracting 
homotopy for this complex by inserting the map $ \lambda: B \cotimes \mathcal{K}_G \rightarrow B \cotimes B \cotimes \mathcal{K}_G \cotimes \mathcal{K}_G  
\cong (B \cotimes \mathcal{K}_G)^{\cotimes 2} $ defined by $ \lambda(x \otimes k) = 1 \otimes x \otimes \sigma(k) $ in the first 
tensor factor. Similarly, the second column of the bicomplex associated to $ B $ is linearly contractible since $ B $ is unital. Hence the Bar-complexes of 
$ B \cotimes \mathcal{K}_G $ and $ B $ are linearly homotopy equivalent. \\
Now consider the first columns. We view $ \mathcal{O}_G \cotimes B \cotimes \mathcal{K}_G $ as a bimodule over $ B \cotimes \mathcal{K}_G $ in two different ways. 
Both bimodules $ M $ and $ N $ have the obvious right action by multiplication. The left action on $ M $ is given by 
\begin{align*}
(x \otimes k) * (f &\otimes y \otimes l)(s, r, t) = f(s) \otimes (s^{-1} \cdot x \otimes s^{-1} \cdot k)(y \otimes l)(r,t) \\
&= f(s) \otimes (s^{-1} \cdot x)y \int_G k(sr, sp)l(p,t) dp
\end{align*}
whereas the left action in $ N $ is 
$$ 
(x \otimes k) \cdot (f \otimes y \otimes l)(s,r,t) = f(s) \otimes (s^{-1} \cdot x) y \int_G k(r,p) l(p,t) dp.
$$ 
The crucial point is that there is a bimodule isomorphism $ \phi: N \rightarrow M $ given by 
$$ 
\phi(f \otimes x \otimes k)(s,r,t) = f(s) \otimes x \otimes k(sr,t). 
$$
Using the map $ \phi $ we obtain a linear isomorphism of complexes between the 
first columns of $ \Omega_G(B \cotimes \mathcal{K}_G) $ and $ \Omega_G(B \cotimes \mathcal{K}) $ where 
$ \mathcal{K} $ is the algebra $ \mathcal{K}_G $ equipped with the trivial $ G $-action. Under this isomorphism $ \Tr $ corresponds to the 
trace map $ \tau: \Omega_G(B \cotimes \mathcal{K}) \rightarrow \Omega_G(B) $ 
given by 
\begin{align*}
\tau(f(s) &\otimes (x_0 \otimes  k_0) d(x_1 \otimes k_1) \cdots d(x_n \otimes k_n)) \\
&= f(s) \otimes x_0 dx_1 \cdots dx_n \int k_0(r_0,r_1) k_1(r_1,r_2) \cdots k_n(r_n, r_0) dr_0 \cdots dr_n 
\end{align*}
on the first column. Let us show that this map is a linear homotopy equivalence on the  first columns of the bicomplexes associated to 
$ \Omega_G(B \cotimes \mathcal{K}) $ and $ \Omega_G(B) $. The function $ \chi \in \D(G) $ chosen above determines an idempotent 
$ p = \chi \otimes \chi $ in $ \mathcal{K} $. This idempotent induces an equivariant homomorphism $ \iota: B \rightarrow B \cotimes \mathcal{K} $ 
by defining $ \iota(x) = x \otimes p $ and a corresponding chain map 
$ \Omega_G(\iota): \Omega_G(B) \rightarrow \Omega_G(B \cotimes \mathcal{K}) $. One immediately checks the relation $ \tau \Omega_G(\iota) = \id $ on $ \Omega_G(B) $. 
As in the proof of Morita invariance in ordinary Hochschild homology we construct a presimplicial homotopy between $ \Omega_G(\iota) \tau $ and the identity as 
follows \cite{Loday}. For $ j = 0, \dots, n $ we define on the first column of $ \Omega_G(B \cotimes \mathcal{K}) $ the operator 
\begin{align*}
h_j(x_0 &\otimes |p_0 \ket \bra q_0| \otimes \cdots x_n \otimes |p_n \ket \bra q_n|) 
= x_0 \otimes |p_0 \ket  \bra \chi| \otimes x_1 \otimes |\chi\ket \bra \chi| \otimes \cdots \\
&\cdots \otimes x_j \otimes |\chi \ket \bra \chi| \otimes 
1 \otimes |\chi \ket \bra q_j| \otimes x_{j + 1} \otimes |p_{j + 1} \ket \bra q_{j + 1}| \otimes \cdots \otimes x_n \otimes |p_n \ket \bra q_n|.
\end{align*}
It is straightforward to verify that this yields indeed a presimplicial homotopy between $ \Omega_G(\iota) \tau $ and $ \id $ for the equivariant Hochschild 
operator on the first column of $ \Omega_G(B \cotimes \mathcal{K}) $. \qed \\
Since the map $ \Tr: \Omega_G(\mathcal{K}_G) \rightarrow \Omega_G(\mathbb{C}) $ is a linearly split surjection we obtain a 
linearly split exact sequence of paramixed complexes 
\begin{equation*} 
 \xymatrix{
    K \;\; \ar@{>->}[r] &
         \Omega_G(\mathcal{K}_G)\ar@{->>}[r] & \Omega_G(\mathbb{C})  }
\end{equation*}
where $ K $ is the kernel of $ \Tr $. From proposition \ref{KG-Hochschild} we deduce that 
$ K $ is linearly contractible with respect to the Hochschild boundary. \\
Recall from section \ref{secX} the definition of the Hodge tower of a paramixed complex and consider the 
$ n $-th level $ \theta^n K $ of the Hodge tower of $ K $. The Hodge filtration yields a finite decreasing filtration of $ \theta^n K $. 
Since $ K $ is contractible with respect to $ b $ it follows that the paracomplex 
\begin{equation*}
F^p \theta^n K/F^{p + 1} \theta^n K = 
   \xymatrix{
      {b(K_{p + 1})\;} \ar@<1ex>@{->}[r]^-{B} &
          {\; K_{p + 1}/b(K_{p + 2})} 
            \ar@<1ex>@{->}[l]^-{b} 
               } 
\end{equation*}
is covariantly contractible for all $ p $. \\
If $ P $ is a relatively projective paracomplex of covariant pro-modules the Hodge filtration of $ \theta^n K $ induces a finite decreasing 
filtration of the supercomplex $ \SHom_G(P,\theta^n K) $. Since this filtration is bounded the associated spectral 
sequence converges and one gets  
\begin{equation*}
H_*(\SHom_G(P,\theta^n K)) = 0
\end{equation*}
for all $ n $ by our previous argument. 
\begin{lemma} \label{Milnor} With the notation as above put
$ C_n = \SHom_G(P, \theta^n K) $. Then there exists an exact
sequence
\begin{equation*}
\xymatrix{
 {H_0(\SHom_G(P, \theta K))\;} \ar@{->}[r] \ar@{<-}[d] &
     H_0(\prod_{n \in \mathbb{N}} C_n) \ar@{->}[r] &
       H_0(\prod_{n \in \mathbb{N}} C_n)\ar@{->}[d] \\
   {H_1(\prod_{n \in \mathbb{N}} C_n)\;} \ar@{<-}[r] &
    {H_1(\prod_{n \in \mathbb{N}} C_n)}  \ar@{<-}[r] &
     {H_1(\SHom_G(P, \theta K))} \\
}
\end{equation*}
\end{lemma}
\proof First remark that each $ C_n $ is indeed a complex.
We let $ C $ be the corresponding inverse system of complexes.
Using Milnor's description of $ \varprojlim^1 $ we obtain an exact
sequence of supercomplexes
\begin{equation*}
 \xymatrix{
    \varprojlim_n  C_n \;\; \ar@{>->}[r] & \prod_{n \in \mathbb{N}} C_n \ar@{->}[r]^{\id - \sigma} &
        \prod_{n \in \mathbb{N}} C_n \ar@{->>}[r] & \varprojlim^1_n C_n  }
\end{equation*}
where $ \sigma $ denotes the structure maps in $ (C_n)_{n \in \mathbb{N}} $. 
Since all structure maps in $ \theta K $ are linearly split surjections and 
$ P $ is relatively projective the structure maps in the inverse system 
$ (C_n)_{n \in\mathbb{N}} $ are surjective. This implies 
$ \varprojlim^1 C_n = 0 $. Therefore the exact sequence above
reduces to a short exact sequence 
\begin{equation*}
\xymatrix{
    \varprojlim_n  C_n \;\; \ar@{>->}[r] & \prod_{n \in \mathbb{N}} C_n \ar@{->>}[r] &
        \prod_{n \in \mathbb{N}} C_n }
\end{equation*}
of supercomplexes. The associated long exact sequence in homology yields the 
claim. \qed 
\begin{theorem} \label{secvarsimp}
Let $ G $ be a locally compact group. Then there exists a natural isomorphism 
\begin{equation*}
HP^G_*(A, \mathbb{C}) \cong H_*(\SHom_G(X_G(\mathcal{T}(A \hat{\otimes} \mathcal{K}_G)), X_G(\mathbb{C})).
\end{equation*}
for every pro-$ G $-algebra $ A $. 
\end{theorem}
\proof According to theorem \ref{homotopyeq} we have a natural isomorphism 
\begin{equation*}
HP^G_*(A,\mathbb{C}) \cong 
H_*(\SHom_G(X_G(\mathcal{T}(A \hat{\otimes} \mathcal{K}_G)),
\theta \Omega_G(\mathcal{K}_G))) 
\end{equation*}
for every  pro-$ G $-algebra $ A $. Moreover the paracomplex $ P = \theta \Omega_G(A \hat{\otimes} \mathcal{K}_G) $ 
is relatively projective due to corollary \ref{localprj}. Consider the linearly split extension of paracomplexes 
\begin{equation*}
\xymatrix{
    \theta K \;\; \ar@{>->}[r] & \theta \Omega_G(\mathcal{K}_G) \ar@{->>}[r] &
        \theta \Omega_G(\mathbb{C}). }
\end{equation*}
This extension induces a short exact sequence of supercomplexes 
\begin{equation*}
\xymatrix{
    \SHom_G(P,\theta K) \;\; \ar@{>->}[r] & \SHom_G(P,\theta \Omega_G(\mathcal{K}_G)) \ar@{->>}[r] &
        \SHom_G(P,\theta \Omega_G(\mathbb{C})). }
\end{equation*}
The supercomplex $ \SHom_G(P,\theta K) $ is acyclic according to lemma \ref{Milnor}. 
Hence the map $ \Tr: \Omega_G(\mathcal{K}_G) \rightarrow \Omega_G(\mathbb{C}) $ induces an isomorphism  
\begin{equation*}
HP^G_*(A, \mathbb{C}) \cong H_*(\SHom_G(X_G(\mathcal{T}(A \hat{\otimes} \mathcal{K}_G)), \theta \Omega_G(\mathbb{C})).
\end{equation*}
Using theorem \ref{homotopyeq} we can pass to the $ X $-complex $ X_G(\mathcal{T}\mathbb{C}) $ in the second variable again. 
Since the $ G $-algebra $ \mathbb{C} $ is quasifree composition with the chain map 
$ X_G(\mathcal{T}\mathbb{C}) \rightarrow X_G(\mathbb{C}) $ induced by the projection $ \mathcal{T}\mathbb{C} \rightarrow \mathbb{C} $
is a homotopy equivalence. This yields the assertion. \qed \\ 
If $ G $ is discrete this description of $ HP^G_*(A, \mathbb{C}) $ can be simplified further. It is easy to check that in this case the map 
$ \mathcal{O}_G \rightarrow \mathbb{C} $ induced by integration of functions with respect to the counting measure yields an isomorphism
\begin{equation*}
\SHom_G(M, \mathcal{O}_G) \cong \Hom_G(M, \mathbb{C}) \cong \Hom(M_G,\mathbb{C})
\end{equation*}
for every covariant module $ M $ where $ M_G $ denotes the quotient of $ M $ obtained by taking $ G $-coinvariants. 
Let us denote by $ \Omega_G(A \hat{\otimes} \mathcal{K}_G)_G $ the mixed complex obtained by 
taking coinvariants in $ \Omega_G(A \hat{\otimes} \mathcal{K}_G) $. 
Using the previous observation, lemma \ref{XC} and theorem \ref{homotopyeq} we see that theorem \ref{secvarsimp} implies the following result. 
\begin{theorem} Let $ G $ be a discrete group and let $ A $ be a pro-$ G $-algebra. There is a natural isomorphism
\begin{equation*}
HP^G_*(A,\mathbb{C}) \cong
HP^*(\Omega_G(A\hat{\otimes} \mathcal{K}_G)_G)
\end{equation*}
where $ HP^*(\Omega_G(A\hat{\otimes} \mathcal{K}_G)_G) $ denotes the 
periodic cyclic cohomology of the mixed complex 
$ \Omega_G(A\hat{\otimes} \mathcal{K}_G)_G $.
\end{theorem}
For the remaining part of this section $ G $ will be discrete.
In order to complete the proof of theorem \ref{DGJ} we shall show that the 
mixed complexes $ \Omega_G(A\hat{\otimes} \mathcal{K}_G)_G $ and $ \Omega(A \rtimes G) $ 
have isomorphic periodic cyclic cohomologies. We view $ s \in G $ as element of 
$ \mathbb{C}G $ or $ \mathcal{O}_G $ in the canonical way. 
Moreover we write $ T = \sum_{r,s} T_{rs} [r,s] $ for an element 
$ \sum_{r,s} T_{rs}\; r \otimes s $ in $ \mathcal{K}_G $ in the sequel and occasionally 
omit tensor signs in order to improve legibility. \\ 
We define the map $ \phi: \Omega(A \rtimes G) \rightarrow
\Omega_G(A \hat{\otimes} \mathcal{K}_G)_G $ on $ n $-forms by
\begin{align*}
\phi((a_0 \rtimes s_0) d(a_1 \rtimes s_1) \cdots & d(a_n \rtimes s_n))
= s_0 \cdots s_n \otimes a_0[e,s_0]d(s_0 \cdot a_1)[s_0,s_0s_1] \cdots \\
&\qquad \cdots d(s_0 \cdots s_{n - 1} \cdot a_n)[s_0 \cdots s_{n - 1},s_0 \cdots s_n] 
\end{align*}
for $ a_0 \rtimes s_0 \in A \rtimes G $ and
\begin{align*}
\phi(d(a_1 \rtimes s_1) \cdots & d(a_n \rtimes s_n))
= s_1 \dots s_n \otimes d a_1[e,s_1] d(s_1 \cdot a_2)[s_1,s_1s_2] \cdots \\
&\qquad \cdots d(s_1 \cdots s_{n - 1} \cdot a_n)[s_1\cdots s_{n - 1},s_1 \cdots s_n].
\end{align*}
The map $ \tau: \Omega_G(A \hat{\otimes} \mathcal{K}_G)_G
\rightarrow \Omega(A \rtimes G) $ is defined by
\begin{align*}
&\tau(s \otimes (a_0 \otimes T^0) d(a_1 \otimes T^1) \cdots d(a_n \otimes T^n)) \\
&= \sum_{r_0,\dots,r_n \in G} (r_0^{-1} \cdot a_0 \rtimes T^0_{r_0r_1} r_0^{-1} r_1)
d(r_1^{-1} \cdot a_1 \rtimes T^1_{r_1r_2} r_1^{-1} r_2) \cdots \\
&\qquad \cdots d(r_n^{-1} \cdot a_n \rtimes T^n_{r_n, s r_0} r_n^{-1} s r_0)
\end{align*}
for $ a_0 \otimes T^0 \in A \hat{\otimes} \mathcal{K}_G $ and
\begin{align*}
&\tau(s \otimes d(a_1 \otimes T^1) \cdots d(a_n \otimes T^n)) \\
&= \sum_{r_1,\dots,r_n \in G}
d(r_1^{-1} \cdot a_1 \rtimes T^1_{r_1r_2} r_1^{-1} r_2) \cdots 
d(r_n^{-1} \cdot a_n \rtimes T^n_{r_n, s r_1} r_n^{-1} s r_1).
\end{align*}
Observe that the sums occuring here are finite since only finitely many 
entries in the matrices $ T^j $ are nonzero. 
\begin{prop} The bounded linear maps $ \phi: \Omega(A \rtimes G) \rightarrow
\Omega_G(A \hat{\otimes} \mathcal{K}_G)_G $ and
$ \tau: \Omega_G(A \hat{\otimes} \mathcal{K}_G)_G
\rightarrow \Omega(A \rtimes G) $ are maps of mixed complexes
and we have $ \tau \phi = \id $.
\end{prop}
\proof The formulas given above clearly define bounded linear maps. 
Remark that $ \tau $ is well-defined since it vanishes on 
coinvariants.
It is immediate from the definitions that $ \phi $ and $ \tau $
commute with $ d $. A direct calculation shows that
both maps also commute with the Hochschild operators. 
This implies that $ \phi $ and $ \tau $ are maps of mixed complexes. 
Furthermore one computes easily that $ \tau \phi $ is equal to the identity on 
$ \Omega(A \rtimes G) $. This yields the claim. \qed\\
We calculate explicitly 
\begin{align*}
(&\phi \tau)(s \otimes(a_0 \otimes T^0) d(a_1 \otimes T^1) \cdots d(a_n \otimes T^n)) \\
&= \phi\biggl(\;\sum_{r_0,\dots,r_n \in G} (r_0^{-1} \cdot a_0 \rtimes  T^0_{r_0r_1} r_0^{-1} r_1)
d(r_0^{-1} \cdot a_1 \rtimes T^1_{r_1r_2} r_1^{-1} r_2) \cdots \\
&\qquad \cdots d(r_n^{-1} \cdot a_n \rtimes T^n_{r_n, s r_0} r_n^{-1} s r_0) \biggr) \\
&= \sum_{r_0,\dots,r_n \in G} r_0^{-1} s r_0 \otimes (r_0^{-1} \cdot a_0 \otimes T^0_{r_0r_1} [e,r_0^{-1} r_1])
d(r_0^{-1} \cdot a_1 \otimes T^1_{r_1 r_2} [r_0^{-1} r_1, r_0^{-1} r_2]) \\
& \qquad \qquad \cdots d(r_0^{-1} \cdot a_n \otimes T^n_{r_n, s r_0} [r_0^{-1} r_n, r_0^{-1} s r_0]) \\
&= \sum_{r_0,\dots,r_n \in G} s \otimes (a_0 \otimes T^0_{r_0r_1} [r_0, r_1])
d(a_1 \otimes T^1_{r_1 r_2} [r_1, r_2]) \cdots
d(a_n \otimes T^n_{r_n, s r_0} [r_n, s r_0]).
\end{align*}
In the same way one obtains
\begin{align*}
(\phi \tau)(&s \otimes d(a_1 \otimes T^1) \cdots d(a_n \otimes T^n)) \\
&= \sum_{r_1,\dots,r_n \in G} s \otimes d(a_1 \otimes T^1_{r_1 r_2} [r_1, r_2]) \cdots
d(a_n \otimes T^n_{r_n, s r_1} [r_n, s r_1]).
\end{align*}
\begin{prop}\label{DGJhomotopy} Let $ G $ be a discrete group and assume that $ A $ is a unital pro-$ G $-algebra. 
Then the map $ \phi \tau: \Omega_G(A \hat{\otimes} \mathcal{K}_G)_G \rightarrow
\Omega_G(A \hat{\otimes} \mathcal{K}_G)_G $ is homotopic to the identity 
with respect to the Hochschild boundary.
\end{prop}
\proof We construct a chain homotopy 
connecting $ \id $ and $ \phi \tau $ on the Hochschild complex
associated to the mixed complex $ \Omega_G(A \hat{\otimes} \mathcal{K}_G)_G $. \\
Let us associate to an element of the form
$ s \otimes a_0[r_0,s_0]da_1[r_1,s_1] \cdots da_n[r_n,s_n] $ a certain number $ M $.
If $ s_j = r_{j + 1} $ for all $ j = 0, \dots, n - 1 $ and
$ s^{-1}s_n = r_0 $ we set $ M = \infty $. If at least one of
these conditions is not fulfilled, we let $ M $ be the smallest
number $ i $ such that $ s_i \neq r_{i + 1} $ (or $ M = n $ if
all $ s_j = r_{j + 1} $ for $ j = 0,\dots,n - 1 $ and
$ s^{-1} s_n \neq r_0 $).
In a similar way we proceed with elements of the form
$ s \otimes da_1[r_1,s_1] \cdots da_n[r_n,s_n] $. Here the 
first condition disappears and the last condition becomes $ s^{-1}s_n = r_1 $. 
The number $ M $ is then defined as before. \\
We construct bounded linear maps
$ h: \Omega^n_G(A \hat{\otimes} \mathcal{K}_G)_G
\rightarrow \Omega^{n + 1}_G(A \hat{\otimes} \mathcal{K}_G)_G $ for all $ n $ as
follows. For an element
$ s \otimes a_0[r_0,s_0]da_1[r_1,s_1] \cdots da_n[r_n,s_n] $
we set
\begin{align*}
&h(s \otimes a_0[r_0,s_0]da_1[r_1,s_1] \cdots da_n[r_n,s_n]) \\
&= (-1)^M s \otimes a_0[r_0,s_0]da_1[r_1,s_1]\cdots da_M[r_M,s_M] d1_A[s_M,s_M] \cdots 
da_n[r_n,s_n] 
\end{align*}
if $ M < \infty $ and
\begin{equation*}
h(s \otimes a_0[r_0,s_0]da_1[r_1,s_1] \cdots da_n[r_n,s_n]) = 0
\end{equation*}
if $ M = \infty $. Here $ 1_A $ denotes the unit of $ A $.\\
For elements of the form
$ s \otimes da_1[r_1,s_1] \cdots da_n[r_n,s_n] $
we have to distinguish four cases. The first case is $ s^{-1} s_n = r_1 $ and 
$ M < \infty $. In this case we set
\begin{align*}
&h(s \otimes da_1[r_1,s_1] \cdots da_n[r_n,s_n]) \\
&= (-1)^M s \otimes da_1[r_1,s_1]\cdots da_M[r_M,s_M] d1_A[s_M,s_M] \cdots 
da_n[r_n,s_n] 
\end{align*}
as before. The second case is $ s^{-1} s_n \neq r_1 $
and $ M = n $. We set
\begin{align*}
h(s &\otimes da_1[r_1,s_1] \cdots da_n[r_n,s_n]) \\
&= (-1)^M s \otimes da_1[r_1,s_1]\cdots da_n[r_n,s_n] d1_A[s_n,s_n] \\
&+ (-1)^{M + n} s \otimes d1_A[s^{-1}s_n,s^{-1} s_n] da_1[r_1,s_1] \cdots da_n[r_n,s_n].
\end{align*}
The third case is $ s^{-1} s_n \neq r_1 $ and $ M < n $.
We set
\begin{align*}
&h(s \otimes da_1[r_1,s_1] \cdots da_n[r_n,s_n]) \\
&= (-1)^M s \otimes da_1[r_1,s_1]\cdots da_M[r_M,s_M] d1_A[s_M,s_M] \cdots da_n[r_n,s_n] \\
&\qquad + (-1)^{M + n} s \otimes (s^{-1} \cdot a_n)[s^{-1}r_n, s^{-1} s_n] d1_A[s^{-1} s_n,s^{-1} s_n] da_1[r_1,s_1] \cdots \\
&\qquad \cdots da_M[r_M,s_M] d1_A[s_M,s_M] \cdots da_{n - 1}[r_{n - 1},s_{n - 1}].
\end{align*}
Finally if $ M = \infty $ we set
\begin{equation*}
h(s \otimes da_1[r_1,s_1] \cdots da_n[r_n,s_n]) = 0.
\end{equation*}
Remark that in all cases coinvariants are mapped to coinvariants
and hence $ h $ is well-defined. \\
A lengthy but straightforward computation shows $ bh + hb = \id - \phi \tau $. \qed
\begin{prop}\label{DGJcor}
Let $ G $ be a discrete group and let $ A $ be any pro-$ G $-algebra. The periodic cyclic cohomologies of $ \Omega(A \rtimes G) $ and
$ \Omega_G(A \hat{\otimes} \mathcal{K}_G)_G $ are isomorphic. Inverse
isomorphisms are induced by the maps $ \phi $ and $ \tau $.
\end{prop}
\proof This follows after dualizing from proposition
\ref{DGJhomotopy} using excision, the SBI-sequence
and the fact that periodic cyclic cohomology is the direct limit of the
cyclic cohomology groups. \qed \\
This finishes the proof of theorem \ref{DGJ}.


\begin{thebibliography}{99}
\bibitem{AM} Artin, M., Mazur, B., \a'Etale Homotopy,
Lecture Notes in Mathematics 100, Springer, 1969
\bibitem{Blackadar} Blackadar, B., $ K $-theory for operator algebras, second 
edition, Mathematical Sciences Research Institute Publications 5, Cambridge 
University Press, 1998
\bibitem{Blanc} Blanc, P., Cohomologie diff\'erentiable et changement de 
groupes, Ast\'erisque 124 - 125 (1985), 113 - 130
\bibitem{Block} Block, J., Excision in cyclic homology of topological algebras,
Harvard university thesis, 1987
\bibitem{BG} Block, J., Getzler, E., Equivariant cyclic homology and
equivariant differential forms, Ann. Sci. \'Ecole. Norm. Sup. 27 (1994), 493 - 527
\bibitem{Brylinski1} Brylinski, J.-L., Algebras associated with
group actions and their homology, Brown university preprint, 1986
\bibitem{Brylinski2} Brylinski, J.-L., Cyclic homology and
equivariant theories, Ann. Inst. Fourier 37 (1987), 15 - 28
\bibitem{Bues1} Bues, M., Equivariant differential forms and 
crossed products, Harvard university thesis, 1996
\bibitem{Bues2} Bues, M., Group actions and quasifreeness, preprint, 1998
\bibitem{Burghelea} Burghelea, D., The cyclic homology of the group rings, 
Comment. Math. Helv. 60 (1985), 354 - 365
\bibitem{Cartan1} Cartan, H., Notions d'alg\a`ebre diff\a'erentielle;  
applications aux groupes de Lie et aux vari\'et\'es o\a`u op\a`ere un
groupe de Lie, Colloque de topologie, C.B.R.M. Brussels (1950), 15 - 27
\bibitem{Cartan2} Cartan, H., La transgression dans un groupe de 
Lie et dans un espace fibr\'e principal, Colloque de topologie, C.B.R.M. Brussels (1950), 57 - 71
\bibitem{Connes1} Connes, A., Noncommutative differential geometry,
Publ. Math. IHES 39 (1985), 257 - 360
\bibitem{Connes2} Connes, A., Noncommutative Geometry, Academic Press,
1994
\bibitem{Crainic} Crainic, M., Cyclic homology of smooth groupoids: The general case, 
$ K $-theory 17 (1999), 319 - 362
\bibitem{CQ1} Cuntz, J., Quillen, D., Algebra extensions and
nonsingularity, J. Amer. Math. Soc. 8 (1995), 251 - 289
\bibitem{CQ2} Cuntz, J., Quillen, D., Cyclic homology and
nonsingularity, J. Amer. Math. Soc. 8 (1995), 373 - 442
\bibitem{CQ3} Cuntz, J., Quillen, D., Operators on noncommutative 
differential forms and cyclic homology, in: Geometry, Topology and Physics, 
77 - 111, Internat. Press, 1995
\bibitem{CQ4} Cuntz, J., Quillen, D., Excision in bivariant
periodic cyclic cohomology, Invent. Math. 127 (1997), 67 - 98
\bibitem{FT} Feigin, B. L., Tsygan, B. L., Additive $ K $-theory,
Lecture Notes in Mathematics 1289, Springer, 1987, 67 - 209
\bibitem{GJ1} Getzler, E., Jones, J. D. S., The cyclic homology of
crossed product algebras, J. Reine Angew. Math. 445 (1993), 161 - 174
\bibitem{Green} Green, P., Equivariant $ K $-theory and crossed product 
$ C^* $-algebras, Proc. Sympos. Pure Math. 38, Amer. Math. Soc., Providence, 
1982, 337 - 338
\bibitem{Grothendieck} Grothendieck, A., Produits tensoriel topologiques
et espaces nucl\'eaires, Mem. Amer. Math. Soc. 16, 1955
\bibitem{H-L1} Hogbe-Nlend, H., Compl\'etion, tenseurs et nucl\'earit\'e 
en bornologie, J. Math. Pures Appl. 49 (1970), 193 - 288
\bibitem{H-L2} Hogbe-Nlend, H., Bornologies and functional analysis, 
North-Holland Publishing Co., 1977
\bibitem{Julg} Julg, P., $ K $-th\'eorie \'equivariante et produits crois\'es,
C. R. Acad. Sci. Paris 292 (1981), 629 - 632
\bibitem{Kasparov1} Kasparov, G. G., The operator $ K $-functor and 
extensions of $ C^* $-algebras, Izv. Akad. Nauk SSSR Ser. Mat. 44 (1980), 571 - 636
\bibitem{Kasparov2} Kasparov, G. G., Equivariant $ KK $-theory and 
the Novikov conjecture, Invent. Math. 91 (1988), 147 - 201
\bibitem{Kassel} Kassel, C., Homologie cyclique, caract\`ere de Chern et lemme de 
perturbation, J. Reine Angew. Math. 408 (1990), 159 - 180
\bibitem{KKL1} Klimek, S., Kondracki, W., Lesniewski, A.,
Equivariant entire cyclic cohomology, I. Finite groups, $ K $-Theory 4 (1991), 
201 - 218
\bibitem{KKL2} Klimek, S., Lesniewski, A., Chern character in
equivariant entire cyclic cohomology, $ K $-Theory 4 (1991),
219 - 226
\bibitem{Loday} Loday, J.-L., Cyclic Homology, Grundlehren der
Mathematischen Wissenschaften 301, Springer, 1992
\bibitem{Meyer} Meyer, R., Analytic cyclic cohomology, Preprintreihe 
SFB 478, Geometrische Strukturen in der Mathematik, Heft 61, M\"unster, 
\bibitem{Meyersmoothrep} Meyer, R., Smooth group representations on bornological vector spaces, Bull. Sci. Math. 128 (2004), 127 - 166 
\bibitem{Meyerborntop} Meyer, R., Bornological versus topological analysis in metrizable spaces, to appear in 
Conference Proceedings Banach Algebras 2003, Contemporary Mathematics
\bibitem{Nistor1} Nistor, V., Group cohomology and the cyclic cohomology
of crossed products, Invent. Math. 99 (1990), 411 - 424
\bibitem{Nistor2} Nistor, V., Cyclic cohomology of crossed products by 
algebraic groups, Invent. Math. 112 (1993), 615 - 638
\bibitem{Segal} Segal, G., Equivariant $ K $-theory, Publ. Math. IHES 34 (1968), 
129 - 151
\bibitem{Voigtthesis} Voigt, C., Equivariant cyclic homology, Preprintreihe 
SFB 478, Geometrische Strukturen in der Mathematik, Heft 287, M\"unster 

\end{thebibliography}
\end{document}